\documentclass[preprint]{elsarticle}
\usepackage[bookmarks=true,colorlinks=true,linkcolor=blue]{hyperref}

\usepackage{amsmath,amssymb,mathabx,mathtools}
\usepackage{graphicx,esint,ulem}

\biboptions{sort&compress}
\usepackage[usenames,dvipsnames,svgnames,table]{xcolor}

\usepackage{color}
\definecolor{jwbGreen}{rgb}{0, .6, 0}
\definecolor{jbaPurple}{HTML}{6600FF}
\definecolor{purple}{rgb}{.7, 0., .8}

\usepackage{calc}
\usepackage[margin=1.in]{geometry}

\usepackage{fancyvrb}

\usepackage{empheq}
\usepackage[most]{tcolorbox}

\newtcbox{\mymath}[1][]{%
    nobeforeafter, math upper, tcbox raise base,
    enhanced, colframe=blue!60!black,
    colback=blue!20, boxrule=1pt,
    #1}

\usepackage{algorithm,algpseudocode}
\usepackage{algorithmicx}
\algrenewcommand\alglinenumber[1]{\footnotesize #1:} 
\newcommand{\algFontSize}{\footnotesize}

\usepackage{tcolorbox}
\usepackage[english]{babel}
\usepackage{amsmath}
\usepackage{amsfonts}
\usepackage{fullpage}
\usepackage{graphicx}
\usepackage{setspace}
\usepackage{float}

\usepackage{tabu}

\usepackage{multirow}

\usepackage{mdframed}

\allowdisplaybreaks

\usepackage{tikz}

\newtheorem{theorem}{Theorem}

\newenvironment{proof}[1][Proof]{\begin{trivlist}
\item[\hskip \labelsep {\bfseries #1.}]}{\end{trivlist}}

\newcommand{\citeCount}[1]{}

\newcommand{\bogus}[1]{{}}
\newcommand{\strutt}{\rule{0pt}{10pt}}

\newcommand{\mni}{\medskip\noindent}

\newcommand{\shadedBoxWithShadow}[3]{\begin{empheq}[box={\mymath[colback=#2!20,colframe=#2!60!black,drop lifted shadow, sharp corners]}]{#1} #3\end{empheq}}

\newcommand{\ECR}{{\rm ECR}}
\newcommand{\CR}{{\rm CR}}

\newtheorem{recipe}{Recipe}

\newlength{\ycbTop}
\newlength{\ycbMid}%

\newcommand{\p}{\partial}

\newcommand{\Np}{{N_p}}

\newcommand{\tableFont}{\footnotesize}

\newcommand{\num}[2]{#1e#2} 

\newcommand{\f}[2]{\frac{#1}{#2}}

\def\ba#1\ea{\begin{align}#1\end{align}}

\def\bas#1\eas{\begin{align*}#1\end{align*}}

\def\bat#1\eat{\begin{alignat}{3}#1\end{alignat}}

\def\bats#1\eats{\begin{alignat*}{3}#1\end{alignat*}}

\newcommand{\bse}{\begin{subequations}}
\newcommand{\ese}{\end{subequations}}

\newcommand{\Nt}{N_t}

\newcommand{\Dzt}{D_{0t}}
\newcommand{\Dpt}{D_{+t}}
\newcommand{\Dmt}{D_{-t}}
\newcommand{\dt}{\Delta t}
\newcommand{\dx}{\Delta x}

\newcommand{\dr}{{\Delta r}}

\newcommand{\eqdef}{\overset{{\rm def}}{=}}

\newcommand{\bv}{\mathbf{ b}}

\newcommand{\ev}{\mathbf{ e}}
\newcommand{\fv}{\mathbf{ f}}

\newcommand{\iv}{\mathbf{ i}}
\newcommand{\jv}{\mathbf{ j}}

\newcommand{\rv}{\mathbf{ r}}

\newcommand{\vv}{\mathbf{ v}}
\newcommand{\wv}{\mathbf{ w}}
\newcommand{\xv}{\mathbf{ x}}

\newcommand{\Gv}{\mathbf{ G}}

\newcommand{\half}{{1\over2}}

\newcommand{\Real}{{\mathbb R}}

\newcommand{\zerov}{\mathbf{0}}

\newcommand{\Ac}{{\mathcal A}}
\newcommand{\Bc}{{\mathcal B}}

\newcommand{\Dc}{{\mathcal D}}
\newcommand{\Ec}{{\mathcal E}}

\newcommand{\Gc}{{\mathcal G}}

\newcommand{\Lc}{{\mathcal L}}

\newcommand{\Nc}{{\mathcal N}}

\newcommand{\Wc}{{\mathcal W}}

\newcommand{\grad}{\nabla}

\newcommand{\nd}{n_d}

\newcommand{\eps}{\epsilon}

\newcommand{\ds}{{\Delta s}}

\newcommand{\Dpx}{D_{+x}}
\newcommand{\Dmx}{D_{-x}}

\newcommand{\ssf}{\scriptscriptstyle}

\newcommand{\Tbar}{{\widebar T}}
\newcommand{\Tf}{\Tbar}

\newcommand{\sinc}{{\rm sinc}}
\newcommand{\Uhat}{\hat{U}}

\newcommand{\omegaHat}{\hat{\omega}}
\newcommand{\omegaTilde}{\tilde{\omega}}
\newcommand{\omegaExplicit}{\omega_{e}}
\newcommand{\omegaImplicit}{\omega_{i}}
\newcommand{\Ttilde}{\tilde{T}}

\newcommand{\vHat}{\hat{v}}
\newcommand{\fHat}{\hat{f}}

\newcommand{\wHat}{\hat{w}}
\newcommand{\What}{\hat{W}}
\newcommand{\Vhat}{\hat{V}}
\newcommand{\uHat}{\hat{u}}

\newcommand{\Gcd}{\Gc_d}
\newcommand{\NITS}{N_{\ssf ITS}}
\newcommand{\Nits}{\NITS}

\newcommand{\lambdaTilde}{\tilde{\lambda}}
\newcommand{\lambdaExplicit}{\lambda^e_{h,m}}
\newcommand{\lambdaImplicit}{\lambda^i_{h,m}}

\newcommand{\sincd}{\sinc_d}

\newcommand{\ACRe}{\mu_{{\ssf E},h}}
\newcommand{\ACRi}{\mu_{{\ssf I},h}}

\newcommand{\PPW}{{\rm PPW}}

\newcommand{\alphag}{a_{g}}
\newcommand{\betag}{b_{g}}

\newcommand{\Ep}{\Ec}

\newcommand{\kTilde}{\tilde{k}}
\newcommand{\kappaTilde}{\tilde{\kappa}}
\newcommand{\up}{u^f}
\newcommand{\Up}{U^f}
\newcommand{\Ef}{E^f}

\newcommand{\Eamp}{\Ec_{\ssf A}}
\newcommand{\Ephase}{\Ec_{\phi}}

\newcommand{\deltak}{\delta k}

\newlength{\tfwidth}
\newlength{\tfheight}
\newlength{\tfxa}
\newlength{\tfxb}
\newlength{\tfya}
\newlength{\tfyb}
\newcommand{\trimFigWithBox}[6]{%
\setlength\fboxsep{0pt}%
\setlength\fboxrule{1.0pt}
\fbox{\includegraphics[width=#2, clip, trim=#3 #4 #5 #6]{#1}}%
}
\newcommand{\trimFigNoBox}[6]{%
\setlength\fboxsep{1pt}
\setlength\fboxrule{0.0pt}
\fbox{\includegraphics[width=#2, clip, trim=#3 #4 #5 #6]{#1}}%
}
\newcommand{\trimFigHeightWithBox}[6]{%
\setlength\fboxsep{0pt}%
\setlength\fboxrule{1.0pt}
\fbox{\includegraphics[height=#2, clip, trim=#3 #4 #5 #6]{#1}}%
}
\newcommand{\trimFigHeightNoBox}[6]{%
\setlength\fboxsep{1pt}
\setlength\fboxrule{0.0pt}
\fbox{\includegraphics[height=#2, clip, trim=#3 #4 #5 #6]{#1}}%
}

\newcommand{\trimFig}[6]{%
\setlength{\tfwidth}{(#2+#2*\real{#3})+#2*\real{#4}}
\setlength{\tfheight}{(#2+#2*\real{#5})+#2*\real{#6}}%
\setlength{\tfxa}{\tfwidth*\real{#3}}%
\setlength{\tfxb}{\tfwidth*\real{#4}}%
\setlength{\tfya}{\tfheight*\real{#5}}%
\setlength{\tfyb}{\tfheight*\real{#6}}%
\trimFigNoBox{#1}{#2}{\tfxa}{\tfya}{\tfxb}{\tfyb}%
}

\newsavebox\figBox

\newcommand{\trimw}[6]{%
\sbox\figBox{\includegraphics{#1}}
\setlength{\tfwidth}{\the\wd\figBox}
\setlength{\tfheight}{\the\ht\figBox}
\setlength{\tfxa}{\tfwidth*\real{#3}}%
\setlength{\tfxb}{\tfwidth*\real{#4}}%
\setlength{\tfya}{\tfheight*\real{#5}}%
\setlength{\tfyb}{\tfheight*\real{#6}}%
\trimFigNoBox{#1}{#2}{\tfxa}{\tfya}{\tfxb}{\tfyb}%
}

\newcommand{\trimwb}[6]{%

\sbox\figBox{\includegraphics{#1}}
\setlength{\tfwidth}{\the\wd\figBox}
\setlength{\tfheight}{\the\ht\figBox}
\setlength{\tfxa}{\tfwidth*\real{#3}}%
\setlength{\tfxb}{\tfwidth*\real{#4}}%
\setlength{\tfya}{\tfheight*\real{#5}}%
\setlength{\tfyb}{\tfheight*\real{#6}}%
\trimFigWithBox{#1}{#2}{\tfxa}{\tfya}{\tfxb}{\tfyb}%
}

\newcommand{\trimh}[6]{%
\sbox\figBox{\includegraphics{#1}}
\setlength{\tfwidth}{\the\wd\figBox}
\setlength{\tfheight}{\the\ht\figBox}
\setlength{\tfxa}{\tfwidth*\real{#3}}%
\setlength{\tfxb}{\tfwidth*\real{#4}}%
\setlength{\tfya}{\tfheight*\real{#5}}%
\setlength{\tfyb}{\tfheight*\real{#6}}%
\trimFigHeightNoBox{#1}{#2}{\tfxa}{\tfya}{\tfxb}{\tfyb}%
}

\newcommand{\trimhb}[6]{%

\sbox\figBox{\includegraphics{#1}}
\setlength{\tfwidth}{\the\wd\figBox}
\setlength{\tfheight}{\the\ht\figBox}
\setlength{\tfxa}{\tfwidth*\real{#3}}%
\setlength{\tfxb}{\tfwidth*\real{#4}}%
\setlength{\tfya}{\tfheight*\real{#5}}%
\setlength{\tfyb}{\tfheight*\real{#6}}%
\trimFigHeightWithBox{#1}{#2}{\tfxa}{\tfya}{\tfxb}{\tfyb}%
}
\usepackage{xargs}

\newcommandx{\figByHeight}[9][5=0, 6=0, 7=0, 8=0,9=]{
\draw (#1,#2) node[anchor=south west,xshift=-16pt,yshift=-4pt] {\trimh{#3}{#4}{#5}{#6}{#7}{#8}};}
\newcommandx{\figByHeightb}[9][5=0, 6=0, 7=0, 8=0,9=]{
\draw (#1,#2) node[anchor=south west,xshift=-16pt,yshift=-4pt] {\trimhb{#3}{#4}{#5}{#6}{#7}{#8}};}

\newcommandx{\figByHeightWithLabel}[9][5=0, 6=0, 7=0, 8=0,9=]{
\draw (#1,#2) node[anchor=south west,xshift=-16pt,yshift=-4pt] {\trimh{#3}{#4}{#5}{#6}{#7}{#8}} node[draw=white,fill=white,inner sep=1pt,anchor=south west] {#9};}
\newcommandx{\figByHeightWithLabelb}[9][5=0, 6=0, 7=0, 8=0,9=]{
\draw (#1,#2) node[anchor=south west,xshift=-16pt,yshift=-4pt] {\trimhb{#3}{#4}{#5}{#6}{#7}{#8}} node[draw=white,fill=white,inner sep=1pt,anchor=south west] {#9};}

\newcommandx{\figByWidth}[9][5=0, 6=0, 7=0, 8=0,9=]{
\draw (#1,#2) node[anchor=south west,xshift=-16pt,yshift=-4pt] {\trimw{#3}{#4}{#5}{#6}{#7}{#8}};}
\newcommandx{\figByWidthb}[9][5=0, 6=0, 7=0, 8=0,9=]{
\draw (#1,#2) node[anchor=south west,xshift=-16pt,yshift=-4pt] {\trimwb{#3}{#4}{#5}{#6}{#7}{#8}};}

\newcommandx{\figByWidthWithLabel}[9][5=0, 6=0, 7=0, 8=0,9=]{
\draw (#1,#2) node[anchor=south west,xshift=-16pt,yshift=-4pt] {\trimw{#3}{#4}{#5}{#6}{#7}{#8}} node[draw=white,fill=white,inner sep=1pt,anchor=south west] {#9};}
\newcommandx{\figByWidthWithLabelb}[9][5=0, 6=0, 7=0, 8=0,9=]{
\draw (#1,#2) node[anchor=south west,xshift=-16pt,yshift=-4pt] {\trimwb{#3}{#4}{#5}{#6}{#7}{#8}} node[draw=white,fill=white,inner sep=1pt,anchor=south west] {#9};}

\usepackage{todonotes}

\definecolor{jwbGreen}{rgb}{0, .6, 0}
\definecolor{jbaPurple}{HTML}{6600FF}
\definecolor{purple}{rgb}{.7, 0., .8}
\definecolor{pinegreen}{rgb}{0.0, 0.47, 0.44}

\begin{document}

\begin{frontmatter}
 \title{An Optimal $O(N)$ Helmholtz Solver for Complex Geometry using WaveHoltz and Overset Grids}

\author[vtu]{Daniel Appel\"o\fnref{DanielThanks}}
\ead{appelo@vt.edu}

\address[vtu]{Department of Mathematics, Virginia Tech, Blacksburg, VA 24061 USA}


\author[rpi]{Jeffrey W.~Banks\fnref{DOE}}
\ead{banksj3@rpi.edu}

\author[rpi]{William D.~Henshaw\corref{cor}\fnref{NSFgrants}}
\ead{henshw@rpi.edu}

\author[rpi]{Donald~W.~Schwendeman\fnref{NSFgrants}}
\ead{schwed@rpi.edu}

\address[rpi]{Department of Mathematical Sciences, Rensselaer Polytechnic Institute, Troy, NY 12180, USA}

\cortext[cor]{Corresponding author}

\fntext[DanielThanks]{Research supported by National Science Foundation under grant DMS-2345225, and Virginia Tech.}

\fntext[NSFgrants]{Research supported by the National Science Foundation under grants DMS-1519934 and DMS-1818926.}

\begin{abstract} 
 We develop efficient and high-order accurate solvers for the Helmholtz equation on complex geometry.
 The schemes are based on the WaveHoltz algorithm which computes solutions of the Helmholtz equation by time-filtering
 solutions of the wave equation. The approach avoids the need to invert an indefinite matrix which can cause convergence
 difficulties for many iterative solvers for indefinite Helmholtz problems. Complex geometry is treated with overset grids
 which use Cartesian grids throughout most of the domain together with curvilinear grids near boundaries. The basic
 WaveHoltz fixed-point iteration is accelerated using GMRES and 
 also by a \textit{deflation} technique using a set of precomputed eigenmodes.
 The solution of the wave equation is solved efficiently with implicit time-stepping using as few as five time-steps
 per period, independent of the mesh size. The time-domain solver is adjusted to remove dispersion errors in time and
 this enables the use of such large time-steps without degrading the accuracy. When multigrid is used to solve the
 implicit time-stepping equations, the cost of the resulting WaveHoltz scheme scales linearly with the total number of grid points $N$ (at fixed
 frequency) and is thus optimal in CPU-time and memory usage as the mesh is refined. A simple rule-of-thumb formula is
 provided to estimate the number of points-per-wavelength required for a $p$-th order accurate scheme which accounts
 for pollution (dispersion) errors. Numerical results are given for problems in two and three space dimensions, to
 second and fourth-order accuracy, and they show the potential of the approach to solve a wide range of large-scale
 problems.
\end{abstract}


\end{frontmatter}

\tableofcontents

\clearpage
\section{Introduction} \label{sec:introduction}

Helmholtz problems commonly arise in applications of engineering and applied sciences involving systems exhibiting
time-harmonic behavior, e.g.~electromagnetics, acoustics, elasticity, quantum mechanics, and other dispersive and
non-dispersive wave propagation problems. 
In this article we develop efficient and accurate Helmholtz solvers for complex geometry.
The schemes are based on the WaveHoltz algorithm~\cite{WaveHoltz1}, which computes solutions to a Helmholtz equation 
by time-filtering solutions to an associated wave equation. WaveHoltz avoids the need to invert an indefinite matrix which can cause convergence 
difficulties for many iterative approaches~\cite{ModernSolversForHelmholtz2017}.
Complex geometry is treated with overset grids which use Cartesian grids throughout
most of the domain together with curvilinear grids near boundaries.
The basic WaveHoltz fixed-point iteration is accelerated using a matrix free GMRES method.
The method is also accelerated using a \textit{deflation} technique whereby selected precomputed eigenmodes, corresponding to the slowest converging components, are removed from the iteration.
 The solution of the wave equation is solved efficiently with
implicit time-stepping using as few as five time-steps per period, independent of the mesh size.
The time-domain solver is adjusted to remove dispersion errors in time and this enables the use of large time-steps without degrading the accuracy.
When multigrid is used to solve the implicit time-stepping equations, the cost of the resulting WaveHoltz scheme 
for a fixed frequency scales linearly with~$N$, the number of grid points (or degrees of freedom in an equivalent finite element method), and is thus optimal
in CPU-time and memory usage as the mesh is refined.
Figure~\ref{fig:multiObjectsGridAndSolution} shows sample results using the new algorithm. The overset grid for a domain with multiple bodies is shown on the left, 
contours of the computed Helmholtz solution are plotted in the middle, and the convergence of the WaveHoltz fixed-point iteration (blue line) and a
GMRES accelerated iteration (red line) are shown on the right. Further details are provided in subsequent sections.

{
\newcommand{\drawContour}[7]{%
\begin{scope}[#1]
\draw(0.0,0) node[anchor=south west,xshift=-4pt,yshift=+0pt] {\trimfiga{#2}{\figWidtha}};
\begin{scope}[xshift=.1cm,yshift=14pt]
  \draw (\xcb,\ycb) node[anchor=south west,xshift=0.25cm,yshift=.5cm,rotate=-90] {\trimfigcb{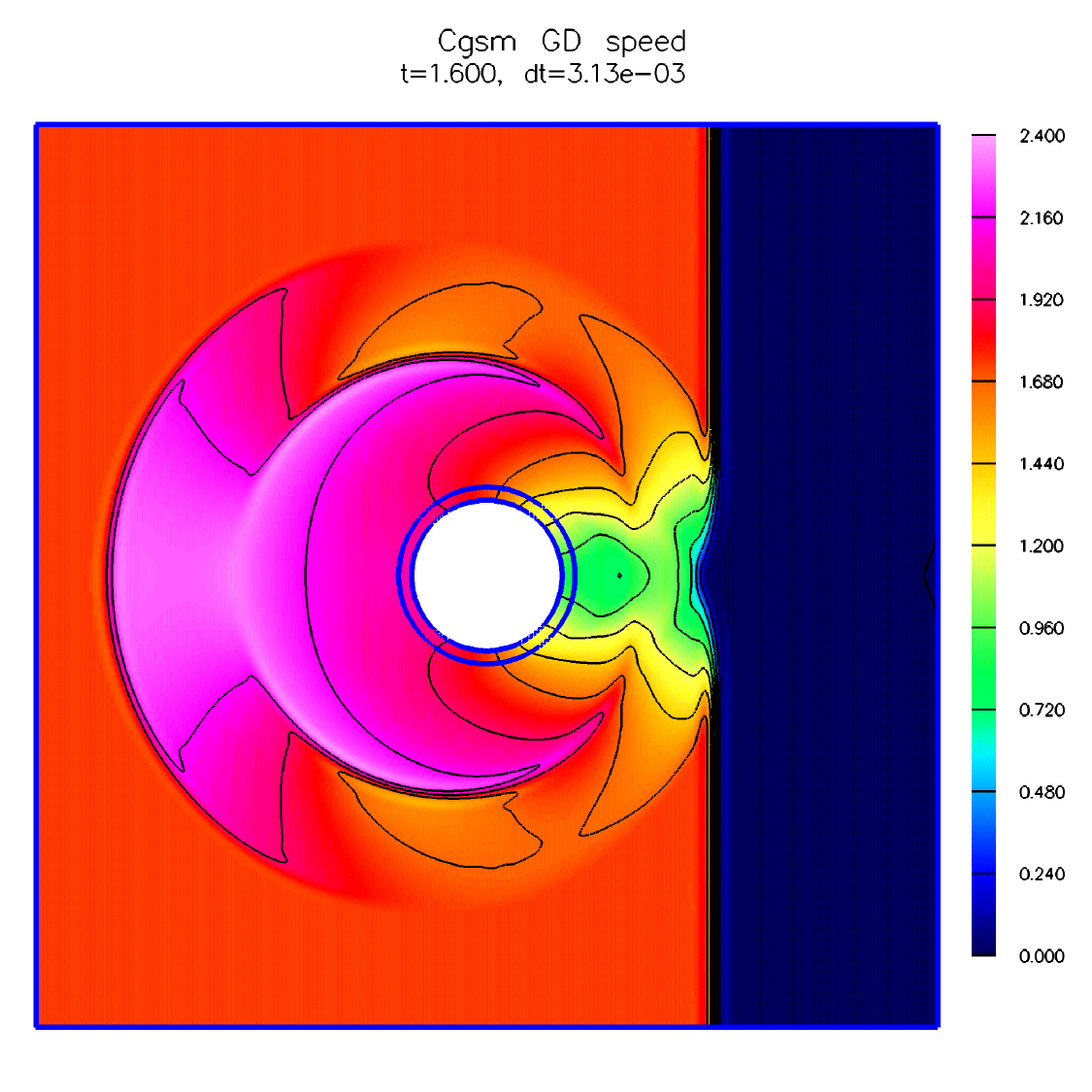}{\cbWidth}{\cbHeight}};
  \draw (.8,0) node[anchor=north,xshift=+3pt,yshift=+2pt] {\scriptsize $#6$};
  \draw (4.8,0) node[anchor=north,xshift=+0pt,yshift=+2pt] {\scriptsize $#7$};
\end{scope}
\end{scope}
}
\newcommand{\cbWidth}{.2cm}
\newcommand{\cbHeight}{4cm}
\newcommand{\xcb}{.5cm}
\newcommand{\ycb}{-.2cm}
\setlength{\ycbTop}{\ycb+\cbHeight}
\setlength{\ycbMid}{\ycb+\cbHeight*\real{.5}}
\newcommand{\trimfigcb}[3]{\includegraphics[width=#2, height=#3, clip, trim=17cm 2.35cm 1.65cm 2.35cm]{#1}}
\newcommand{\figWidtha}{5.25cm}
\newcommand{\trimfiga}[2]{\trimw{#1}{#2}{.11}{.115}{.11}{.17}}
\newcommand{\figSize}{5.5cm}
\begin{figure}[htb]
\begin{center}
\begin{tikzpicture}
   \useasboundingbox (0.2,.45) rectangle (16,4.5);  

   \begin{scope}[yshift=0cm]
     \figByWidthb{   0}{.2}{multiObjectsGridG2}{4.87cm}[0.1][0.1][0.16][0.16]

     \figByWidth{10.7}{-.1}{multiObjectsFreq10Np4G4Order4Deflate128}{\figSize}[0][0][0][0]  
   \end{scope}
    \begin{scope}[xshift=4.8cm,yshift=-.5cm]
      \drawContour{xshift=0cm,yshift=0cm}{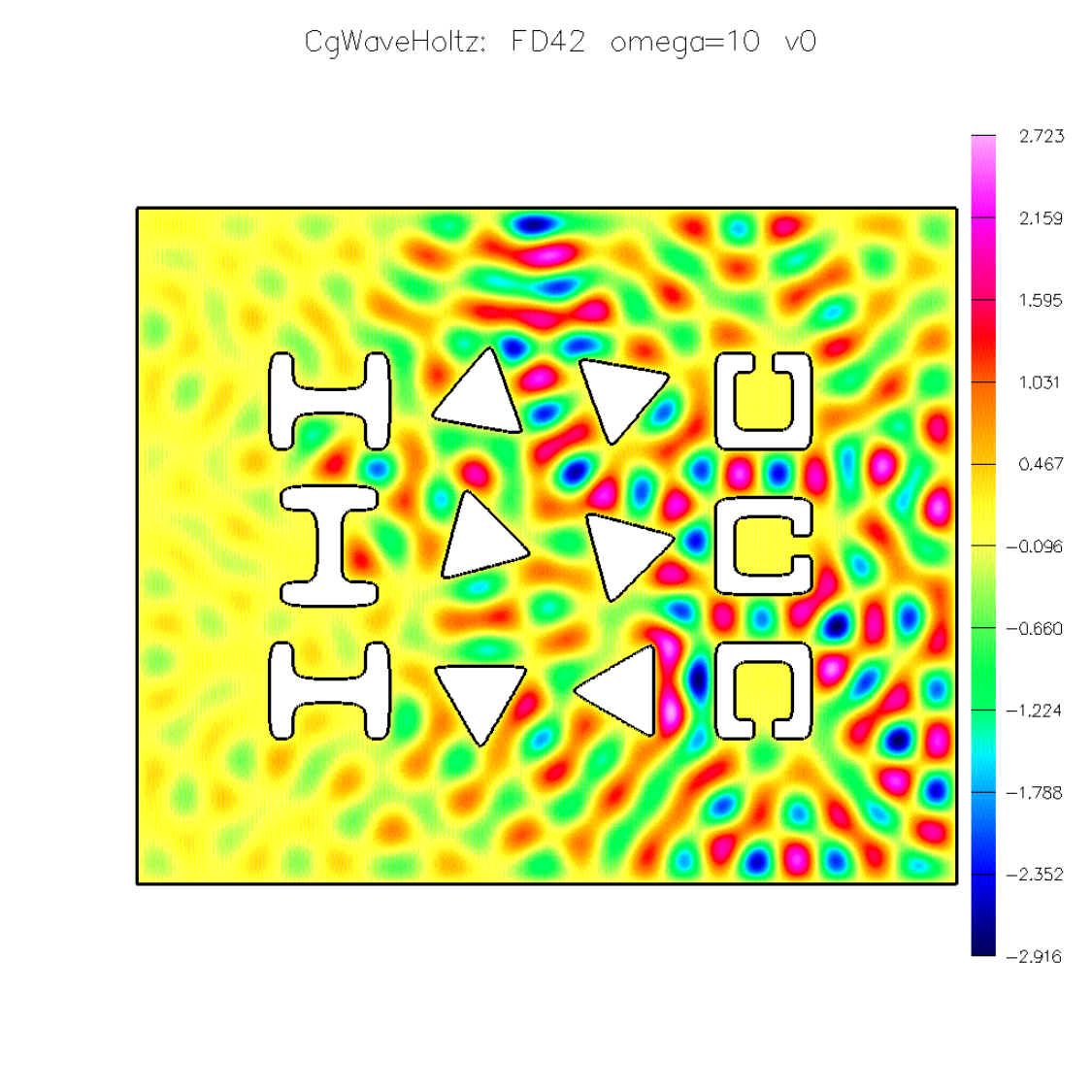}{$v$}{$v$}{$t=0.3$}{$-2.9$}{$2.7$}
     \draw(1.5,1) node[draw,fill=white,anchor=west,xshift=2pt,yshift=2pt,inner sep=3pt] {\scriptsize Gaussian source};  
     \draw[thick,black,->] (2.8,1.27) -- (3.15,2.4);  
   \end{scope}

\end{tikzpicture}
\end{center}
\caption{Gaussian source amongst multiple bodes. 
Left: overset grid (coarse version) consisting of a background blue grid and body fitted grids around each object. Middle: computed WaveHoltz solution using implicit time-stepping with $10$ time-steps per period and deflation. 
Right: WaveHoltz convergence history for the fixed-point iteration and GMRES accelerated iteration (further details are provided in subsequent sections).
    }
\label{fig:multiObjectsGridAndSolution}
\end{figure}
}

Developing efficient solvers for the Helmholtz equation is an important and challenging topic with wide applicability in the applied sciences.
On one hand there are attractive direct methods for the solution of linear system associated with Helmholtz discretizations, e.g.~the Hierarchically Semi-Separable (HSS) parallel multifrontal sparse solver by deHoop and co-authors~\cite{2011_deHoop_Xia}, and the spectral collocation solver by Gillman, Barnett and Martinsson~\cite{Gillman2015}. However, direct solvers have limitations in terms of computational time and memory use as the size of the problem grows, and thus these solvers are not generally a viable option for very large-scale problems. Therefore, for large problems, it is usually necessary to resort to iterative methods.  However, iterative methods often have competing requirements of keeping the number of iterations bounded as the frequency increases while at the same time keeping the memory use, startup costs and time per iteration small as the number of grid points $N$ increases. An ideal iterative solver would have $O(N)$ computational cost and memory use,
 and bounded iteration counts as the frequency increases. Furthermore, for many applications, solutions are needed over a wide range of frequencies, and for changing geometry or material parameters.  Meeting these challenges is notoriously difficult, and it has been the subject of much research (see Ernst and Gander~\cite{ernst2012difficult}, Erlangga~\cite{erlangga2008advances}, or the collection of papers in \cite{ModernSolversForHelmholtz2017} for more details). 

{
\newcommand{\figSize}{6cm}

\begin{figure}[htb]
\begin{center}
\begin{tikzpicture}[scale=1]
  \useasboundingbox (0,.6) rectangle (14,8.2);  

  \begin{scope}[yshift=0cm]
    \figByWidth{0.0}{5.15}{pointsPerWaveLengthSine100WaveLengths}{14cm}[0][0][0][0];
    \figByWidth{4}{0}{pointsPerWaveLengthSineCurves}{\figSize}[0][0][0][0];
     \begin{scope}[xshift=6.9cm,yshift=6.3cm]
       \draw[black, very thick] (0,0) rectangle (.2,1.5);  
       \draw[black, thick,->] (0,0) -- (-2.6,-1.3);
       \draw[black, thick,->] (.2,0) -- (3.0,-1.3);
    \end{scope}

  \end{scope}  

\end{tikzpicture}
\end{center}
\caption{Pollution errors and points-per-wavelength (PPW) comparison. The plots show the grid point distribution for a $p^{\rm th}$-order accurate scheme to reach a
 relative error tolerance of $\eps=10^{-2}$, for a domain that is $N_\Lambda=100$ wavelengths long.
 The sine curves are plotted showing the grid points and are offset vertically for clarity.
 The second-order accurate scheme (red curve) requires a massive number of grid points to manage the pollution (dispersion) errors.
   The formula for $\PPW$ is given by~\eqref{eq:ppwPollution} in Recipe~\ref{recipe:PPW}.
    } 
\label{fig:pointsPerWavelengthSine}
\end{figure}
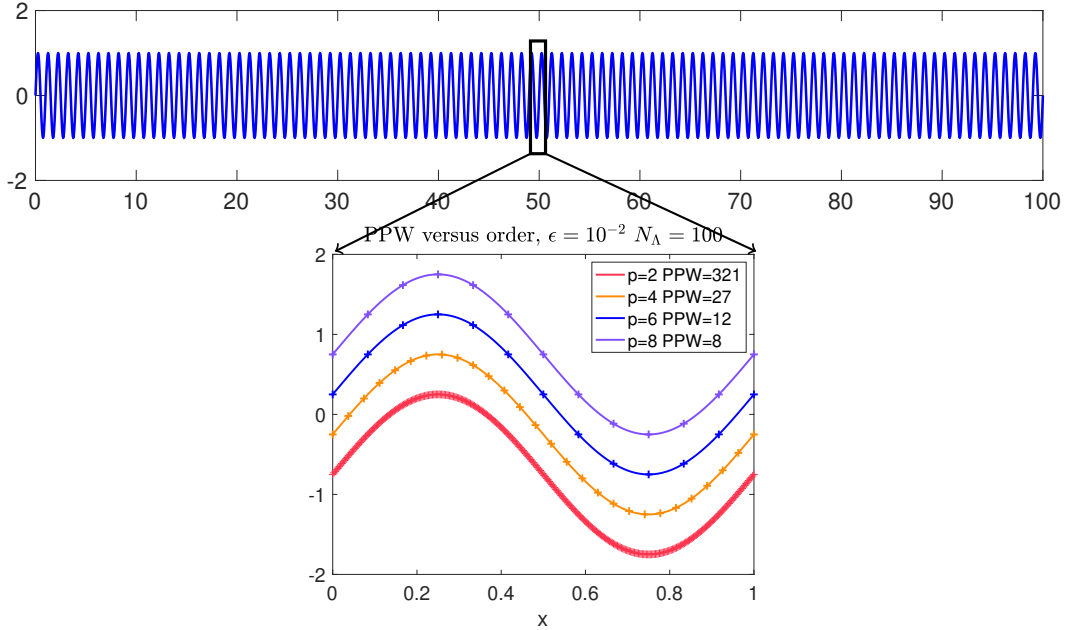
}

A major challenge for solving Helmholtz problems is the highly indefinite character of the discretized system of equation. 
This indefinite nature of the linear systems causes traditional iterative algorithms, such as preconditioned GMRES or multigrid methods, 
to either converge very slowly or not at all~\cite{ernst2012difficult}. 
A second major challenge is the resolution requirements to manage pollution (dispersion) errors at high frequencies~\cite{BaylissGoldsteinTurkel1985,IlenburgBabuska1995}.
The suppression of pollution errors is a serious problem when the domain is large compared to the wavelength $\Lambda=2\pi/k$ 
associated with the wave-number $k=\omega/c$ in definition of the  Helmholtz problem.
The model problem analysis given in Section~\ref{sec:pollutionErrors} provides a rule of thumb that the number of points-per-wavelength (PPW) 
for a $p^{\rm}$-order accurate scheme should be proportional to $(N_\Lambda/\eps)^{1/p}$ where $N_\Lambda$ is the size of the domain in wavelengths
and $\eps$ is the relative error tolerance.
Figure~\ref{fig:pointsPerWavelengthSine} shows the estimated grid resolution requirements for $N_\Lambda=100$ and $\eps=10^{-2}$ (only one wavelength is plotted). A second-order accurate scheme requires a massive $321$ points-per-wavelength (PPW).
A fourth-order accurate scheme requires a much more reasonable $PPW=27$, while sixth and eight-order schemes are even better.
Methods with higher orders of accuracy are thus attractive in terms of reducing the number of grid points and therefore the computational cost\footnote{Assuming the computational cost of the higher-order accurate scheme does not increase too fast with order.}.

In recent years there have been important advances in iterative methods for Helmholtz problems. For example, schemes based on sweeping preconditioners combined with domain-decomposition methods using sophisticated transmission conditions have shown promise~\cite{engquist2011sweepingH,engquist2011sweepingPML,stolk2013rapidly,Helmholtz_Chen_Xiang_2013,PoulsonEngquistLiYing2013,vion2014double,leonardo_1,ZEPEDANUNEZ2016347,2018arXiv180108655Z,GANDER2000261,Gander_Zhang_SIAM_REV}. Similarly, schemes using preconditioners based on a complex shifted Laplacian coupled perhaps with multigrid~\cite{BaylissGoldsteinTurkel1985,LahayeVuik2017,ErlanggaRamosNabben2017} have also been found to be effective. In fact, some of these schemes are able to achieve near $O(N)$ scaling and/or frequency independent iteration counts in some cases. However, it is fair to say that  these methods all have limitations in terms of startup costs, memory use, parallel scalability, and applicability to very large scale problems.


\newcommand{\jj}{j}
\section{Preliminaries} \label{sec:preliminaries}

In this section we introduce the Helmholtz boundary-value problem of interest and the related 
initial-boundary-value problem for the time-domain wave equation with periodic forcing.
We then describe the WaveHoltz fixed-point iteration (FPI).

\subsection{Governing equations}  \label{sec:governing}
  
Consider the problem of finding numerical approximations to solutions $u(\xv)$ of the Helmholtz boundary-value problem (BVP),
\bse
\label{eq:HelmholtzBVP}
\bat
  &  \Lc u + \omega^2 \, u = f(\xv), \qquad&&  \xv\in\Omega, \\
  &  \Bc u = g(\xv) ,                   \qquad && \xv\in\partial\Omega, 
\eat
\ese
on a domain $\Omega\in\Real^{\nd}$ in $\nd$ dimensions with boundary $\partial\Omega$.
Here $\Lc$ is an elliptic operator and $\Bc$ denotes the boundary condition operator.
The frequency $\omega$ is real and we take $\omega\ge 0$ without loss of generality.
The forcing functions $f(\xv)$ and $g(\xv)$ are assumed to be given.
The elliptic operator~$\Lc$, along with boundary operator~$\Bc$, is assumed to be self-adjoint.
For the purposes of this article we take $\Lc= c^2 \Delta$, $c>0$, along with Dirichlet or Neumann boundary conditions,
but we note that the WaveHoltz approach can be extended to 
more general elliptic operators $\Lc$, variable coefficients, and more general boundary conditions~\cite{peng2021emwaveholtz,ElWaveHoltz,WaveHoltzSemi}.

Solutions to the Helmholtz problem~\eqref{eq:HelmholtzBVP} can be found 
by finding time-periodic solutions 
$w(\xv,t)=u(\xv)\cos(\omega t)$
to the associated initial-boundary-value-problem (IBVP) for the wave equation\footnote{For the real-valued Helmholtz solutions found here, the choice of $\cos(\omega t)$ implies the initial conditions in~\eqref{eq:waveEquation}.} given by
\bse
\label{eq:waveEquation}
\bat
   & \p_t^2 w = \Lc w - f(\xv) \cos(\omega t)   ,   \qquad&&  \xv\in\Omega, \\
   & \Bc w = g(\xv)\cos(\omega t) ,                \qquad && \xv\in\partial\Omega, \\
   & w(\xv,0) = w_0(\xv),                          \qquad && \xv\in\Omega, \\
   & \p_t w(\xv,0) = 0,                          \qquad && \xv\in\Omega.
\eat
\ese
The WaveHoltz algorithm defines a procedure for finding the 
initial condition $w_0(\xv)$ in~\eqref{eq:waveEquation} so that the solution $w(\xv,t)$ is time periodic with period $T=2\pi/\omega$.  Once $w_0(\xv)$ is found, then the solution of the corresponding Helmholtz problem is simply $u(\xv)=w_0(\xv)$.

\subsection{The WaveHoltz fixed-point iteration}

\renewcommand{\algFontSize}{\small}
\begin{algorithm}
\algFontSize 
\caption{WaveHoltz Algorithm - Basic Fixed-Point Iteration.}
\begin{algorithmic}[1]

  \Function{WaveHoltz}{$\omega$,$f$,$g$,$\Np$}  
    \State // Final time is $\Tbar=\Np T$, with $T=2\pi/\omega$. Filter uses $\alpha=1/2$. 
    \State $k=0$ \Comment WaveHoltz~iteration counter.
    \State $v^{(k)}=0$   \Comment Assign initial guess for Helmholtz iterate 
    \While{ not converged} \Comment Start WaveHoltz~iterations.

      \State $w^{(k)}(\xv,0) = v^{(k)}(\xv)$ \Comment Initial condition for wave equation solve.
      \State $w^{(k)}(\xv,0:\Tf)$ = \Call{solveWaveEquation}{$w^{(k)}(\xv,0)$,$f$,$g$} \Comment Solve for $\wv(\xv,t)$, $t\in[0,\Tf]$. 
      \State $\displaystyle v^{(k+1)}(\xv) = \f{2}{\Tbar} \int_{0}^{\Tbar} \left( \cos(\omega t) - \f{\alpha}{2} \right) \, w^{(k)}(\xv,t; v^{(k)}) \, dt$
           \Comment Time filter the wave equation solution.
      \State $k = k+1$
    \EndWhile    \Comment End WaveHoltz iterations.
    \State $u(\xv) = v^{(k)}(\xv)$ \Comment Approximate Helmholtz solution.
 \EndFunction
\end{algorithmic} 
\label{alg:waveHoltz}
\end{algorithm}

The basic features of the WaveHoltz algorithm can be described at the continuous level. Details of the discrete approximations
are left to later sections.
Let $v^{(k)}(\xv)$, $k=0,1,2,\ldots,$ denote the $k\sp{{\rm th}}$ iterate in the WaveHoltz algorithm which is an approximate solution to the Helmholtz BVP~\eqref{eq:HelmholtzBVP}.
The basic WaveHoltz fixed-point iteration that generates $\{v^{(k)}\}$ is given in Algorithm~\ref{alg:waveHoltz}.
The input to the algorithm is the frequency $\omega$, the forcing functions $f(\xv)$ and $g(\xv)$, and $\Np$,
 the number of periods over which to integrate (taking $\Np>1$ can sometimes be advantageous as discussed in later sections).
The final time is thus given by $\Tf=\Np T$, where $T$ is the period defined above.
After setting an initial guess for $v^{(0)}(\xv)$ (here taken as zero but any guess could be used) the iteration commences.
At each stage in the iteration the current value of $v^{(k)}(\xv)$ is used as the initial condition for solving
the wave equation IBVP~\eqref{eq:waveEquation}.
Given the solution over time to the wave equation, $w^{(k)}(\xv,t) = w(\xv,t;v^{(k)})$, the new iterate $v^{(k+1)}(\xv)$ is computed using the time filter,
\ba
  v^{(k+1)}(\xv) = \f{2}{\Tbar} \int_{0}^{\Tbar} \left( \cos(\omega t) - \f{\alpha}{2} \right) \, w^{(k)}(\xv,t;v^{(k)}) \, dt  , \label{eq:waveHoltzTimeFilter}
\ea
where $\alpha$ is a constant, often taken to be $1/2$.
Under suitable conditions, $v^{(k+1)}(\xv)$ converges to the Helmholtz solution $u(\xv)$, see Section~\ref{sec:convergenceAnalysis}.


\section{Solving the wave equation in complex geometry using overset grids} \label{sec:numericalApproach}

We have developed two new computer programs to implement the WaveHoltz algorithm on overset grids.
The first program, called CgWaveHoltz, implements the WaveHoltz algorithm.
CgWaveHoltz in turn uses the second program, CgWave, which solves the scalar wave equation.
These programs are built using the Overture framework.\footnote{www.overtureFramework.org and sourceforge.net/projects/overtureframework (for Overture), and sourceforge.net/projects/cgwave (for CgWave and CgWaveHoltz).}
The software is open source and freely available. 
A brief overview of the numerical scheme used by CgWave is given in this section, while further details of our approach to solving wave propagation problems on overset grids
can be found in~\cite{max2006b,smog2012,adegdm2019,adegdmi2020}, for example.

\subsection{Discretizing PDEs on overset grid}

To provide context for our numerical discretizations and to establish some notation, we give a brief overview of the overset grid approach.
As illustrated in Figure~\ref{fig:overlap3dCartoon}, an overset grid, denoted as~$\Gc$, 
consists of a set of component grids $\{G_g\}$, $g=1,\ldots,{\mathcal N}$, that cover the PDE domain~$\Omega$.
The primary motivation for our use of overset grids is to enable efficient finite difference schemes on structured grids, while simultaneously treating complex geometry with high-order accuracy up to and including the boundaries.
In three dimensions, each component grid, $G_g$, is a logically rectangular, curvilinear grid
defined by a smooth mapping from a unit cube parameter space~$\rv$ to physical
space~$\xv$,
\begin{equation}
  \xv = \Gv_g(\rv),\qquad \rv\in[0,1]^3,\qquad \xv\in\Real^3.
\label{eq:gridMapping}
\end{equation}
All grid points in $\Gc$ are classified as discretization, interpolation or unused points~\cite{CGNS}.
The overlapping grid generator {\bf Ogen}~\cite{ogen} from the {\it Overture} framework is used to
construct the overlapping grid information.
In a typical overset grid, one or more boundary-fitted curvilinear grids represent each boundary.
The remainder of the domain is covered by one or more Cartesian grids.
{\bf Ogen} cuts holes in the appropriate component grids by using physical boundaries to distinguish between
the interior and exterior to the domain.
Grid points outside the domain are classified as unused points.
For instance, the \lq\lq cylinder'' grid displayed in the upper right image of Figure~\ref{fig:overlap3dCartoon} cuts a hole in the Cartesian \lq\lq box'' grid so that the latter grid has many unused points (those not being plotted in the lower right image).
{\bf Ogen} also provides the interpolation information for all interpolation points in the overlap region between component grids.

{
 \newcommand{\figWidth}{5.5cm}  
 \newcommand{\figWidtha}{3.75cm}   
 \newcommand{\labelSize}{\small}
%
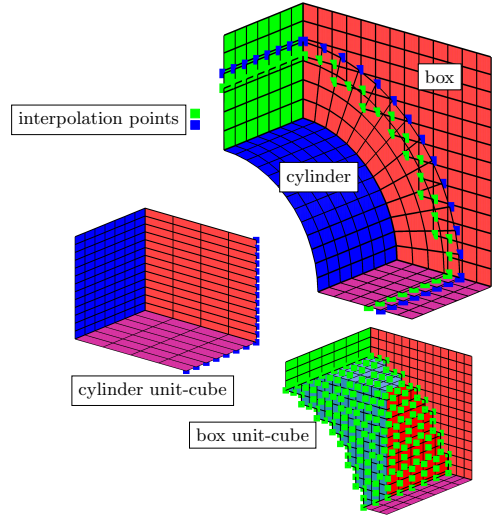
\begin{figure}[hbt]
\newcommand{\trimfiga}[2]{\trimFig{#1}{#2}{.05}{.025}{.1}{.02}}
\newcommand{\trimfigb}[2]{\trimFig{#1}{#2}{.055}{.05}{.1}{.04}}
 \newcommand{\trimfig}[2]{\trimFig{#1}{#2}{.05}{.05}{.0}{.0}}
 \begin{center}
\resizebox{6cm}{!}{
\begin{tikzpicture}[scale=.5]
  \useasboundingbox (-1,-2.8) rectangle (14.5,13.8);  
  \draw(9,8.5) node[]{\trimfig{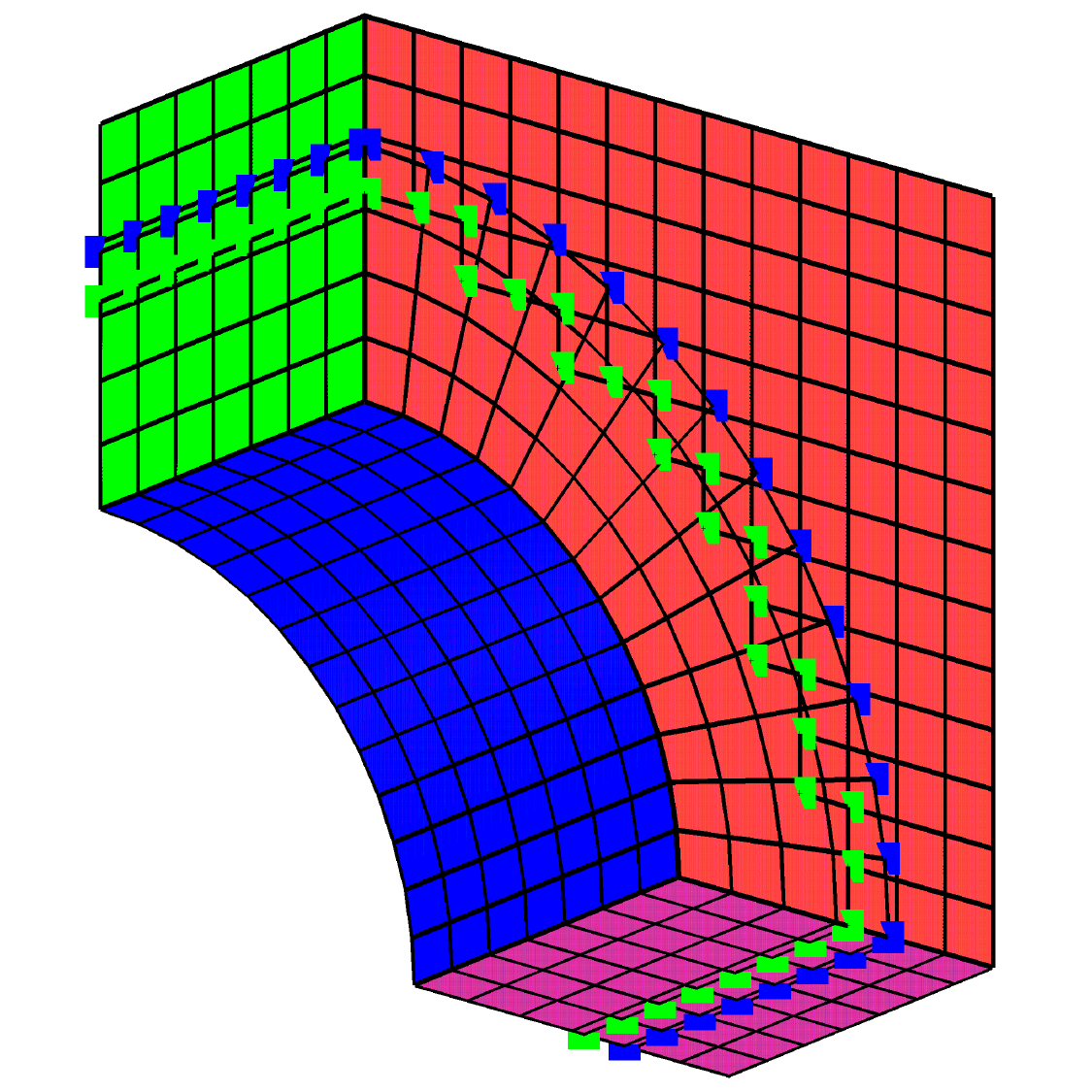}{\figWidth}};
  \draw(2.5,3.8) node[]{\trimfigb{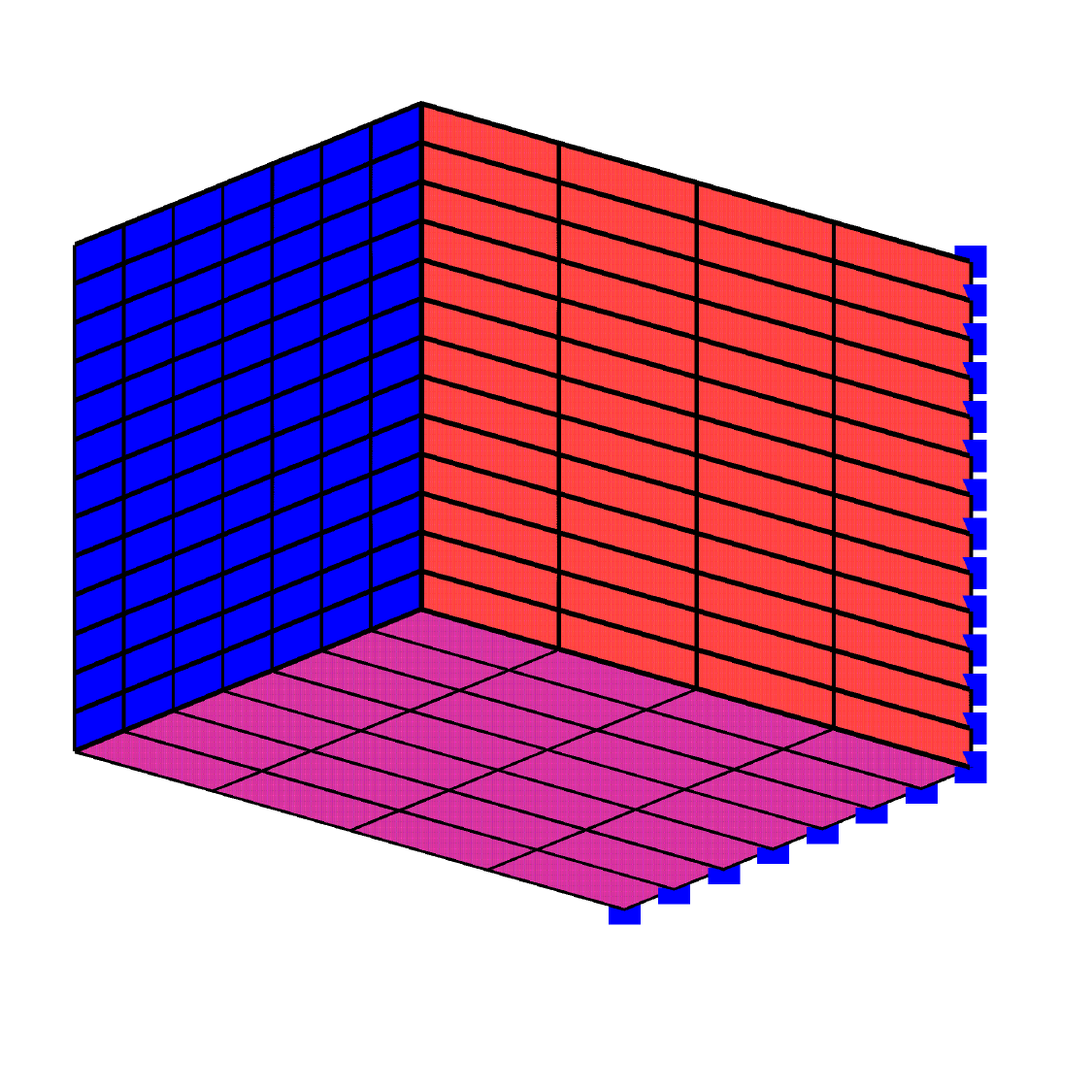}{\figWidtha}};
  \draw(9.6,-.8) node[]{\trimfiga{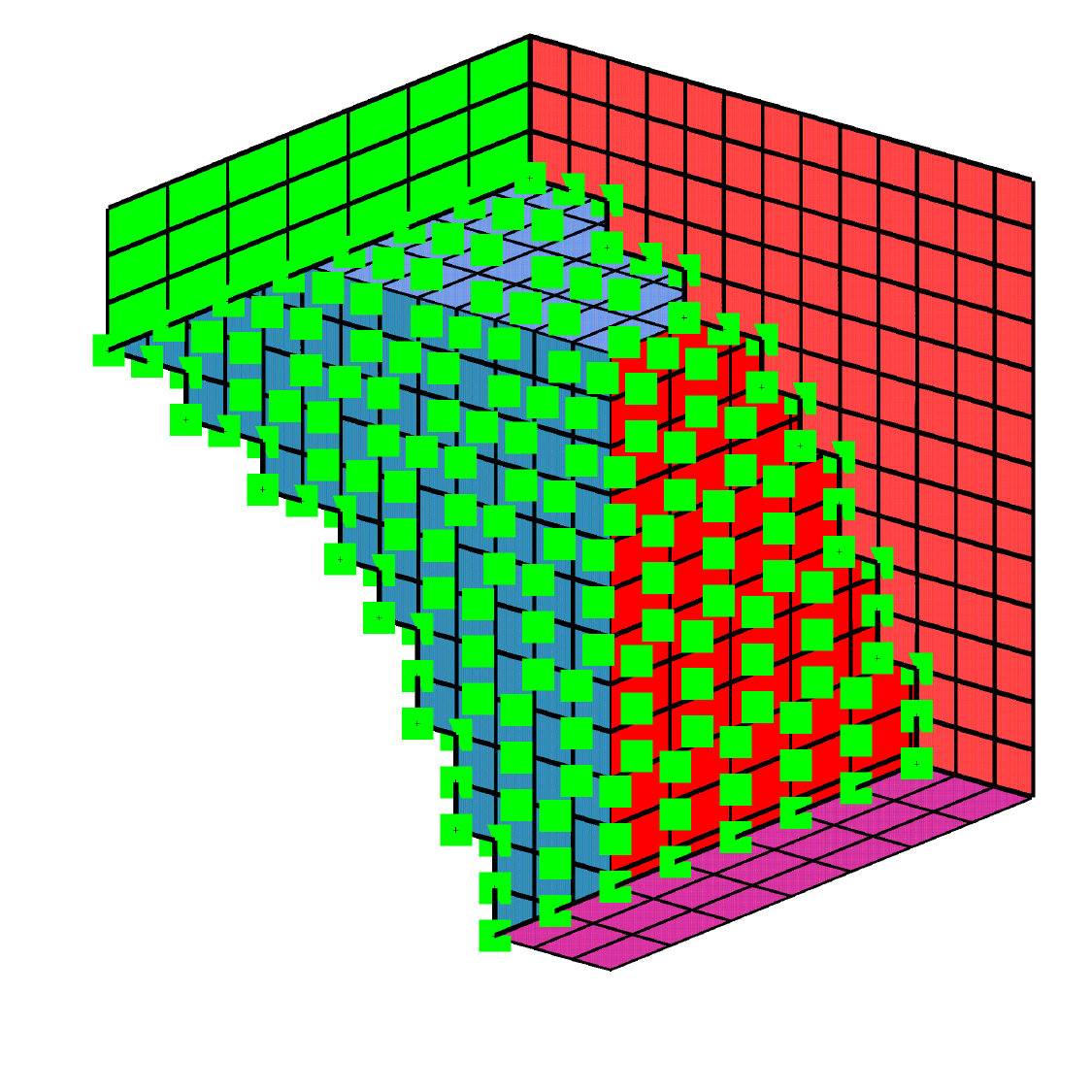}{\figWidtha}};
\draw(12.5,11.5) node[draw,fill=white,anchor=east,xshift=2pt,yshift=-1pt]{\labelSize box};
\draw(8.9,8) node[draw,fill=white,anchor=east,xshift=2pt,yshift=-1pt]{\labelSize cylinder};
\draw(3.,-1.)node[draw,fill=white,anchor=west,xshift=2pt,yshift=-1pt]{\labelSize box unit-cube};
\draw(-1.1,.55)node[draw,fill=white,anchor=west,xshift=2pt,yshift=-1pt]{\labelSize cylinder unit-cube};
\begin{scope}[xshift=-.7cm]
\draw(3.6,10) node[draw,fill=white,anchor=east,xshift=2pt,yshift=-1pt]{\labelSize interpolation points};
\draw[fill=green,draw=green](3.9,10) rectangle ++(.3,.3);
\draw[fill=blue,draw=blue](3.9,9.6) rectangle ++(.3,.3);
\end{scope}
%
\end{tikzpicture}
}
\end{center}
\caption{Top: a three-dimensional overlapping grid for a quarter-cylinder in a box.
    Bottom left and right: component grids for the cylindrical and box grids in the unit cube parameter space.
Interpolation points at the grid overlap are marked and
color-coded for each component grid. }
\label{fig:overlap3dCartoon}
\end{figure}
}

The interpolation between grids is defined using tensor-product Lagrange interpolation in the parameter
space of the mapping $\Gv_g$.  The unit square coordinates $\rv$ of a given point $\xv$ on one grid are located
in the donor-grid parameter space. The interpolation is performed in the Cartesian-grid parameter space and
is thus straightforward~\cite{CGNS}. For a $p^{\rm th}$-order accurate scheme ($p$ even), $p/2$ layers of interpolation points are
required to support the stencil width of $p+1$. An interpolation stencil of width $p+1$ is used, as required for
$p^{\rm th}$-order accuracy on typical grids~\cite{CGNS}. Note that wider interpolation stencils may used for upwind schemes~\cite{ssmx2023} but then the interpolation stencil can still have a width of $p+1$.

Forming approximations to derivatives on a Cartesian grid is straightforward.
Let $\xv_{g,\iv}$ denote the grid points on a grid $g$, where $\iv=[i_1,i_2,i_3]$ is a multi-index with components $i_m=0,1,\ldots,N_m$, where $N_m$ is the number
of grid cells in the $m\sp{{\rm th}}$ direction.
Let $\dx_m=1/N_m$  denote the grid spacing so that $\xv_{g,\iv}=(i_1 \dx_1, i_2\dx_2, i_3\dx_3)$.
Let $W_\iv \approx w(\xv_{g,\iv})$ and define the standard divided difference operators, 
\ba
   D_{+ x_m} W_{\iv} \eqdef \f{u_{\iv+\ev_m} - W_\iv}{\dx_m},        \qquad 
   D_{- x_m} W_{\iv} \eqdef \f{W_{\iv}       - u_{\iv-\ev_m}}{\dx_m}, \qquad 
   D_{0 x_m} W_{\iv} \eqdef \f{u_{\iv+\ev_m} - u_{\iv-\ev_m}}{2\dx_m} ,
\ea
where $\ev_m$ is the unit vector in the $m\sp{{\rm th}}$ direction 
(e.g.~$\ev_2=[0,1,0]$).  Second-order accurate approximations to the Laplacian and gradient in Cartesian coordinates are then
\ba
  & \Delta_h \eqdef \sum_{m=1}^{\nd} D_{+x_m}D_{-x_m} , \qquad 
   \grad_h \eqdef [ D_{0,x_1}, D_{0,x_2}, D_{0,x_3} ]^T .
\ea
High-order accurate approximations use higher-order accurate difference approximations~\cite{max2006b}.

Forming approximations to derivatives on a curvilinear grid is a bit more complicated and there are several ways
to approach this; here we use the \textit{mapping method}.
Given a mapping $\xv = \Gv_g(\rv)$ and its inverse metric derivatives, $\p r_\ell/\p x_m$, $1\le\ell,m\le\nd$, 
the derivatives of a function $w(\xv)=w(\Gv(\rv))=W(\rv)$ are first written in parameter space using the chain rule, for example,
\ba
   \f{\p w}{\p x_m } = \sum_{\ell=1}^{\nd} \f{\p r_\ell}{\p x_m} \f{\p W}{\p r_\ell } . \label{eq:firstDerivCurv}
\ea
Derivatives of $W$ with respect to $r_\ell$ are then approximated with standard finite differences.
Let~$\rv_\iv$ denote grid points on the unit cube, where $i_k=0,1,\ldots,N_k$.
Let $\dr_k=1/N_k$ denote the grid spacing in the $k\sp{{\rm th}}$ direction so that $\rv_{\iv} = (i_1 \dr_1, i_2\dr_2, i_3\dr_3)$. 
Let $W_\iv \approx W(\rv_\iv)$ and define the difference operators,
\ba
   D_{+ r_\ell} W_{\iv} \eqdef \f{U_{\iv+\ev_\ell} - W_\iv}{\dr_\ell},        \qquad 
   D_{- r_\ell} W_{\iv} \eqdef \f{W_{\iv}       - U_{\iv-\ev_\ell}}{\dr_\ell}, \quad 
   D_{0 r_\ell} W_{\iv} \eqdef \f{U_{\iv+\ev_\ell} - U_{\iv-\ev_\ell}}{2\dr_\ell}.
   \label{eq:dividedDiffr}
\ea
Second-order accurate approximations to the first derivatives in~\eqref{eq:firstDerivCurv} are
\ba
   D_{x_m,h} W_\iv \eqdef \sum_{\ell=1}^{\nd} \left.\f{\p r_\ell}{\p x_m}\right\vert_{\iv} D_{0,r_\ell} W_\iv,
\ea
where we assume the metric terms $\p r_\ell/\p x_m$ are known at grid points from the mapping. 
Second derivatives are formed using the chain rule,
\ba
  \f{\p^2 w}{\p x_m \p x_n } = \sum_{k=1}^{\nd} \sum_{l=1}^{\nd} \f{\p r_k}{\p x_m} \f{\p r_l}{\p x_n} \f{\p^2  W}{\p r_k\p r_l }
     +  \sum_{k=1}^{\nd} \left\{ \sum_{l=1}^{\nd} \f{\p r_l}{\p x_n} \f{\p}{\p r_l}  \f{\p r_k}{\p x_m} \right\}  \f{\p W}{\p r_k } .
\ea
The second derivatives are then approximated using finite differences in the parameter space.
We normally do not assume the second derivatives of the mapping are known (to avoid the extra storage) and these 
are computed using finite differences of the metrics.
As an example, second-order accurate  approximations are
\bse
\bat
    \left.\f{\p^2  W}{\p r_k\p r_l }\right\vert_{\rv_\iv} & \approx  D_{+r_k} D_{-r_l} W_\iv , \qquad&& \text{for $k =  l$}, \\
    \left.\f{\p^2  W}{\p r_k\p r_l }\right\vert_{\rv_\iv} & \approx  D_{0r_k} D_{0r_l} W_\iv , \qquad&& \text{for $k\ne l$}, \\
    \left.\f{\p}{\p r_l}\left( \f{\p r_k}{\p x_m}\right)\right\vert_{\rv_\iv} & \approx D_{0r_l} \left( \left.\f{\p r_k}{\p x_m}\right\vert_\iv \right) .
\eat
\ese
Higher order approximations for wave equations on overset grids are discussed further in~\cite{max2006b}.

\subsection{Discretizing the wave equation} \label{sec:discretizingTheWaveEquation}

Now consider solving the IBVP for the wave equation given by~\eqref{eq:waveEquation}.
Let $W_{g,\iv}^n \approx w(\xv_{g,\iv},t^n)$ denote the discrete approximation in space and time on grid~$g$. 
We consider both explicit and implicit methods in time. 
While high-order accurate methods in both space and time are available (see~\cite{lcbc2022} for example), 
we use only second-order accurate schemes in time since we can correct for time discretization errors in the WaveHoltz algorithm using the approaches described in~\ref{sec:timeCorrections}.  
The spatial approximations, on the other hand, are $p^{\rm th}$ order accurate, where $p=2$ and $4$ for the purposes of this paper (although higher-order accurate discretizations are possible). 
The explicit time-stepping scheme takes the form
\ba
  \Dpt\Dmt W_{g,\iv}^n = L_{ph}  W_{g,\iv}^n  + F(\xv_{g,\iv},t^n),
  \label{eq:explicitScheme}
\ea
where $L_{ph}$ denotes a $p^{\rm th}$ order accurate approximation to $\Lc=c^2 \Delta$ and $F(\xv,t)$ is a time-periodic forcing function
whose form is nominally $f(\xv)\cos(\omega t)$ but with adjustments for correcting time-discretization errors as described in subsequent sections.
The implicit time-stepping scheme is given by 
\ba
  \Dpt\Dmt W_{g,\iv}^n = L_{ph}  \Big[ \frac{1}{2} W_{g,\iv}^{n+1} + \frac{1}{2} W_{g,\iv}^{n-1} \Big] 
      + F(\xv_{g,\iv},t^n) .
  \label{eq:implicitScheme}
\ea
These schemes are augmented with appropriate approximations to the initial conditions and boundary conditions,
 and with suitable choices for $F(\xv,t)$ (see Section~\ref{sec:fullyDiscreteConvergence} and~\ref{sec:timeCorrections}).
For stability on overset grids, upwind dissipation would normally be included; the form of this dissipation
is described in~\cite{ssmx2023}. We note, however, that upwind dissipation is not generally needed with the WaveHoltz 
algorithm; this is discussed further in Section~\ref{sec:numerics}.

On a Cartesian grid, the time-step restriction for the explicit scheme~\eqref{eq:explicitScheme}
can be found from a von Neumann analysis, and takes the form
 \ba
    c^2 \,\dt^2\, \left( \sum_{m=1}^{\nd} \f{1}{\dx_m^2 } \right) < C_{2,p} , \label{eq:dtForCartesian}
 \ea  
 where $C_{2,p}$ is a constant that depends on the order of accuracy in time, i.e.~$2$, and the
 order of accuracy in space~$p$.  For example, it is found that~$C_{2,2}=1$ and $C_{2,4}=\sqrt{3}/2 \approx 0.866$.
 The time-step restriction for curvilinear grids is chosen by freezing coefficients and using a von Neumann analysis, and then the time step is chosen to satisfy all frozen coefficient problems.

\renewcommand{\algFontSize}{\small}
\begin{algorithm}
\algFontSize 
\caption{Overset grid algorithm for the wave equation (explicit time-stepping).}
\begin{algorithmic}[1]

  \Function{CgWave}{}  
    \State $W^0$ = assignInitialConditions
    \State $W^1$ = takeFirstStep($W^0$)
    \For{ $n=1,2,\ldots,N_t-1$} \Comment Start time-stepping
      \State $t^n= n\dt$  \Comment current time

      \For{ $g=1,2,\ldots,N_g$}
         \State $W_g^{n+1}$ = advanceGrid( $W_g^n$, $W_g^{n-1}$ ) \Comment Advance component grid $g$
         \State $W_g^{n+1}$ = applyBoundaryConditions( $W_g^{n+1}$, $t^n+\dt$  )
      \EndFor 
      \State $W^{n+1}$=interpolate( $W^{n+1}$ ) \Comment overset grid interpolation
    \EndFor    \Comment End time-stepping loop
 \EndFunction
\end{algorithmic} 
\label{alg:cgWaveExplicit}
\end{algorithm}

Algorithm~\ref{alg:cgWaveExplicit} gives the basic form of the scheme for explicit time-stepping of the wave equation on an overset grid.
At each time-step the solution is advanced independently on each component grid.
After all component grid solutions have been updated to the new time and the boundary conditions
applied, the solution is interpolated to update the solution on interpolation points.
With implicit time-stepping (see Section \ref{sec:implicitTimeStepping}), a sparse linear system of equations is formed representing the 
interior equations on all component grids, as well as equations for the boundary conditions and the interpolation between grids.  These equations can
be solved with a direct sparse solver (for smaller size problems) or iterative methods such
as Krylov methods,\footnote{We obtain good results using a bi-conjugate gradient stabilized scheme with an ILU preconditioner.} algebraic multigrid, or domain decomposition methods. 
Geometric multigrid methods for overset grids can also be used to solve these implicit time-stepping equations as discussed in Section~\ref{sec:optimal}.

\section{Convergence of the WaveHoltz fixed-point iteration}  \label{sec:convergenceAnalysis}

The convergence properties of the WaveHoltz algorithm can be studied through
an eigenfunction analysis. 
Section~\ref{sec:continuousAnalysis} reviews the known convergence results for the continuous problem, and this provides useful background for the new convergence analysis of the fully discrete problem in Section~\ref{sec:fullyDiscreteConvergence}.


\subsection{WaveHoltz convergence: continuous in space and time}  \label{sec:continuousAnalysis}

We first summarize results of the convergence analysis given in~\cite{WaveHoltz1} for the continuous in time and space problem.
The eigenvalue problem associated with the BVP in~\eqref{eq:HelmholtzBVP} is given by
\bse
\label{eq:eigBVP}
\bat
  &  \Lc \phi_m = - \lambda_m^2 \, \phi_m, \qquad&&  \xv\in\Omega, \qquad \\
  &  \Bc \phi_m = 0,  \qquad && \xv\in\partial\Omega, 
\eat
\ese
and since the elliptic operator $\Lc$ with boundary operator $\Bc$ is self-adjoint the eigenvalues $\lambda_m^2$ are real and there exists a complete set of orthogonal eigenfunctions $\phi_m(\xv)$, $m=1,2,\ldots$.  We further assume that the eigenvalues are non-negative and take $\lambda_m \ge 0$ without loss of generality.  The eigenfunctions are normalized so that
\ba
   ( \phi_l, \phi_m )_\Omega = \delta_{l,m},
\ea
where $(\,\cdot\,,\,\cdot\,)_\Omega$ is the usual $L_2$ inner product on $\Omega$ and $\delta_{l,m}$ is the Kronecker delta.

Consider solving the Helmholtz problem~\eqref{eq:HelmholtzBVP} with homogeneous boundary conditions.
Let the solution~$u(\xv)$ and forcing function~$f(\xv)$ be written in terms of the eigenfunction expansions
\ba
  & u(\xv) = \sum_{m=1}^\infty \uHat_{m} \phi_m(\xv), \qquad
    f(\xv) = \sum_{m=1}^\infty \fHat_{m} \phi_m(\xv),
\label{eq:eigenExpansion}
\ea
where $\uHat_{m}$ and $\fHat_{m}$ are generalized Fourier coefficients.
Substituting~\eqref{eq:eigenExpansion} into~\eqref{eq:HelmholtzBVP} leads to the following formula
for the Fourier coefficients of the Helmholtz solution
\ba
  \uHat_{m} = \frac{ \fHat_{m} }{\omega^2 - \lambda_m^2},\qquad m=1,2,3,\ldots\;.
\ea
To study the behavior of the WaveHoltz iteration, we also write $v^{(k)}(\xv)$ and $w^{(k)}(\xv,t)$ in terms of the eigenfunction expansions,
\ba
   v^{(k)}(\xv) = \sum_{m=0}^\infty \vHat^{(k)}_{m} \phi_m(\xv), \qquad
   w^{(k)}(\xv,t) = \sum_{m=0}^\infty \wHat^{(k)}_{m}(t) \, \phi_m(\xv),
\ea
where $\vHat^{(k)}_{m}$ and $\wHat^{(k)}_{m}(t)$ are coefficients in the expansions at the $k\sp{{\rm th}}$ iterate.  Substituting these expressions into the wave equation IBVP~\eqref{eq:waveEquation}
leads to an initial-value problem for each coefficient $\wHat_m^{(k)}(t)$ and whose solution is given by
\ba
  \wHat_m^{(k)}(t) =    \big(\vHat^{(k)}_{m} - \uHat_{m} \big)\, \cos(\lambda_m t) 
                     +  \uHat_{m} \cos(\omega t) .
      \label{eq:wHatFormula}
\ea
Substituting the eigenfunction expansions and the expression for $\wHat_m^{(k)}(t)$  in~\eqref{eq:wHatFormula} 
into the WaveHoltz time filter~\eqref{eq:waveHoltzTimeFilter}
leads to a fixed-point iteration given by
\bse
\label{eq:waveHoltzFourierFPI}
\ba
  \vHat^{(k+1)}_m &=   \big( \vHat^{(k)}_m - \uHat_m \big) \, \beta(\lambda_m) + \uHat_m \, \beta(\omega) , \\
                  &=  \beta(\lambda_m) \,\vHat^{(k)}_m  + \big( \beta(\omega)- \beta(\lambda_m)\big) \uHat_m  , \qquad k=0,1,2,\ldots,
\ea
\ese
where 
\ba
  \beta(\lambda) 
      & \eqdef \f{2}{\Tbar} \, \int_0^\Tbar \left( \cos(\omega t) - \f{\alpha}{2} \right)\, \cos(\lambda t) \, dt ,
  \label{eq:filterFunction}
\ea
is the WaveHoltz filter function.  In addition to the principal dependence on $\lambda$, the filter function also depends on the frequency $\omega$, the final time $\Tbar$ and the filter parameter~$\alpha$, and when appropriate we indicate the dependence on these parameters as $\beta=\beta(\lambda;\,\omega,\Tbar,\alpha)$. 
 It is readily shown that $\beta(\lambda)=1$ when $\lambda=\omega$.
 Thus, if the WaveHoltz fixed-point iteration (FPI) in~\eqref{eq:waveHoltzFourierFPI} converges, and $\lim_{k\rightarrow\infty}\vHat^{(k)}_m = \vHat^*_m$,
 then from~\eqref{eq:waveHoltzFourierFPI}
 \ba
     \vHat^*_m = \f{\beta(\omega)- \beta(\lambda_m)}{1-\beta(\lambda_m)} \, \uHat_m  = \uHat_m,
 \ea 
 and $\vHat^*_m$ are the coefficients in the expansion of the solution of the Helmholtz problem. Also from~\eqref{eq:waveHoltzFourierFPI}
 it is seen that the asymptotic convergence rate depends on $|\beta(\lambda_m)|$.

In order to assess the convergence of the WaveHoltz FPI, it is important to consider the behavior of the WaveHoltz filter function.  To do this, we note that filter function can be written as
\bse
\ba
   \beta(\lambda;\,\omega,\Tbar,\alpha)  &= \sinc\big((\omega-\lambda)\Tbar\big) + \sinc\big((\omega+\lambda)\Tbar\big) - \alpha \, \sinc(\lambda \Tbar),
   \label{eq:filterThreeSincs} 
\ea  
where $\sinc(x)\equiv\sin(x)/x$, or as
\ba    
  \beta(\lambda;\,\omega,\Tbar,\alpha)   &=  \f{2}{\Tbar} \sin( \lambda \Tbar) \, \left( \f{\lambda}{\lambda^2-\omega^2} - \f{ \alpha}{\lambda} \right) .
\ea
\ese
For the typical choice $\alpha=1/2$ it can be shown that $\vert\beta(\lambda)\vert$ has a global maximum equal to $1$ when $\lambda=\omega$ and that $\vert\beta(\lambda)\vert<1$ otherwise (assuming $\lambda\ge0$).  Figure~\ref{fig:waveHoltzBeta} shows plots of $\beta$ versus $\lambda/\omega$ for $\alpha=1/2$ and for $\Np=1,2,3$.
The asymptotic convergence rate $\mu$ of the the WaveHoltz algorithm is generally determined by the value of $|\beta(\lambda_m)|$ for the  eigenvalue $\lambda_m$ closest to $\omega$ (assuming $\lambda_m\ne\omega$).
As $\Np$ increases the main peak near $\lambda=\omega$ narrows and thus the $\mu$, in general, decreases for increasing~$\Np$.

{
\newcommand{\figw}{7cm}
\newcommand{\figh}{6cm}
\begin{figure}[htb]
\begin{center}
\begin{tikzpicture}
  \useasboundingbox (0,.65) rectangle (\figw,.95*\figh);  
  \begin{scope}[yshift=0*\figh]
    \figByWidth{0}{0}{waveHoltzBetaFunction}{\figw}[0.][0.][0.][0.]
  \end{scope}  
\end{tikzpicture}
\end{center}
\caption{
  WaveHoltz filter function $\beta$ for $\Np=1$, $\Np=2$, and $\Np=3$ periods per time-interval.
   }
\label{fig:waveHoltzBeta}
\end{figure}
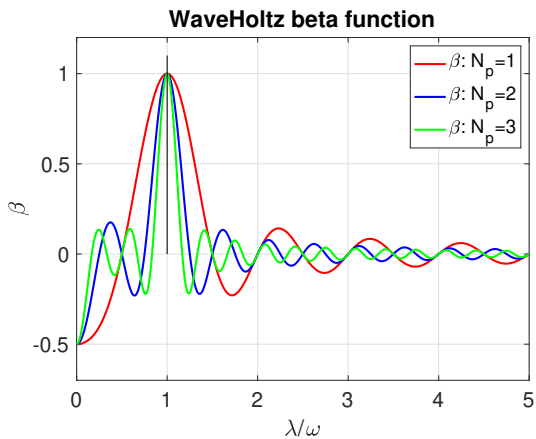
}

These results are summarized in the following theorem.
\begin{theorem}[WaveHoltz FPI Convergence Rate] \label{th:WaveHoltzConvergenceRate}
  Assume $\lambda_m\ne\omega$ are the eigenvalues of the problem in~\eqref{eq:eigBVP} so that $\vert\beta(\lambda_m)\vert<1$ for all $\lambda_m$.  The WaveHoltz fixed-point iteration has asymptotic convergence rate $\mu$ given by
\ba
    \mu = \max_{\lambda_m} |  \beta(\lambda_m) | .
\ea
\end{theorem}
\begin{proof}
  The proof follows from~\eqref{eq:waveHoltzFourierFPI} and the assumptions in the statement of the theorem, see~\cite{WaveHoltz1} for details.
\end{proof}

\subsection{Fully discrete convergence analysis} \label{sec:fullyDiscreteConvergence}

The convergence of the WaveHoltz algorithm for a fully discretized problem is now presented.
Consider a discrete approximation to the wave equation using either explicit or
implicit time-stepping with second-order accuracy in time and $p^{\rm th}$ order accuracy in space.
The WaveHoltz filter function is approximated with a trapezoidal quadrature in time\footnote{Which is \textit{spectrally} accurate for periodic functions.}.
We study a model problem discretized on a single grid.
In this section we take $\Np=1$ for simplicity; the results for $\Np>1$ are similar.

\subsubsection{Explicit time-stepping}

The explicit time-stepping scheme for the wave equation with modified frequency $\omegaExplicit$ (chosen to adjust for time-discretization errors as described below)
takes the form
\bse
\label{eq:WaveSchemeexplicit}
\bat
  & \Dpt\Dmt W_\jv^n = L_{ph} W_\jv^n - f(\xv_\jv) \, \cos(\omegaExplicit t^n),  \qquad&& \jv \in \Omega_h , \quad n=0,1,2,\ldots \label{eq:WaveSchemeexplicitA} \\
  & W_\jv^0 = V_\jv,                                                             \qquad&& \jv \in \bar{\Omega}_h  , \\
  & \Dzt W_\jv^0 = 0,                                                            \qquad&& \jv \in \bar{\Omega}_h  ,    \label{eq:WaveSchemeexplicitC} \\
  & \Bc_{ph} W_\jv^n = 0 ,                                                       \qquad&& \jv \in \p\Omega_h , \quad n=1,2,\ldots
\eat
\ese
where $\Omega_h$ denotes the set of grid points, $\jv$, where the interior equation is applied,
$\bar{\Omega}_h$ denotes the set of all grid points, and $\p\Omega_h$ denotes the set of points
where the boundary conditions are applied.  Here, $\Bc_{ph}$ denotes the discrete boundary condition operator
and the modified frequency $\omegaExplicit$ for explicit time-stepping is
\ba
\omegaExplicit \eqdef {2\over\Delta t}\sin\sp{-1}\left({\omega\Delta t\over2}\right).
\ea
This form for $\omegaExplicit$
is chosen to correct for the time discretization as described in~\ref{sec:timeCorrectionsExplicit} following Recipe~\ref{recipe:OmegaTildeExplicit}.
Note that~\eqref{eq:WaveSchemeexplicitC} can be combined with~\eqref{eq:WaveSchemeexplicitA}, for $n=0$, to eliminate $W_\jv^{-1}$ and this leads 
to an expression for the first time-step,
\ba
    W_\jv^1 = W_\jv^0 + \f{\dt^2}{2} \Big( L_{ph} W_\jv^n - f(\xv_\jv)  \Big). 
\ea

Let us assume that the discrete eigenvalue problem
\bse
\label{eq:discreteEigenvalueBVP}
\bat
  &  L_{ph} \Phi_{m,\jv} = -\lambda_{h,m}^2 \, \Phi_{m,\jv} ,     \qquad&& \jv \in \Omega_h , \\
  & \Bc_{ph} \Phi_{m,\jv} = 0   ,                           \qquad&&  \jv \in \p\Omega_h ,
\eat
\ese
has a complete set of linearly independent eigenvectors $\Phi_{m,\jv}$ 
with corresponding real-valued 
 eigenvalues $\lambda_{h,m}$ for $m=1,2,\ldots,N_{a}$,
where $N_a$ is the total number of approximate eigenmodes. 
Expanding $W_\jv^n$, $V_\jv$, and $f(\xv_\jv)$ in eigenvector expansions with coefficients
$\What_m$, $\Vhat_m$ and $\fHat_m$, respectively,
leads to a discrete ODE for each generalized Fourier coefficient given by
\bse
\label{eq:waveFourierExplicit}
\bat
  & \Dpt\Dmt \What_m^n = L_{ph} \What_m^n - \fHat_m \, \cos(\omegaExplicit t^n),  \quad n=0,1,2,\ldots , \\
  & \What_m^0 = \Vhat_m,                  \\
  & \Dzt \What_m^0 = 0 .      
\eat
\ese
The solution to~\eqref{eq:waveFourierExplicit} takes a similar form to the continuous case (see~\eqref{eq:wHatFormula}) and is 
\bse
\ba
   & \What_m^n = \Uhat_m \cos(\omegaExplicit t^n) + \Vhat_m \cos(\lambdaExplicit t^n) ,  
       \label{eq:WhatExplicit} \\
   & \Uhat_m \eqdef \f{\fHat_m}{\omegaExplicit^2 - \lambda_{h,m}^2} ,  \\
   & \lambdaExplicit \eqdef \f{2}{\dt} \sin^{-1}\left(   \f{\lambda_{h,m}\dt}{2} \right), \label{eq:discreteExplicitLambdaTilde}
\ea
\ese
where the particular solution $\Uhat_m$ is the Fourier coefficient for the solution of the discretized Helmholtz problem.
The WaveHoltz time filter~\eqref{eq:waveHoltzTimeFilter} is approximated using the trapezoidal rule.
Applying this approximate time filter to $\What_m^n$ leads to the update for $\Vhat_m^{(k)}$,
\ba
  \Vhat_m^{(k+1)} = 
    \f{2}{T} \sum_{n=0}^{\Nt} \left( \cos(\omegaExplicit t^n) - \f{\alpha_d }{2}\right) \, \What_m^n \, \sigma_n \, \dt, 
    \label{eq:updateDiscreteExplicit}
\ea
where $\Nt$ is the number of time-steps, and
$\sigma_n$ are quadrature weights given by $\sigma_0=\sigma_{\Nt}=\half$ and $\sigma_n=1$ otherwise. 
The coefficient $\alpha_d = \alpha_d(\omegaExplicit\dt)$ in the discrete filter~\eqref{eq:updateDiscreteExplicit} is an adjusted value for~$\alpha$, derived in~\ref{sec:discreteFilterFunction}, and given by 
\ba
    \alpha_d = \alpha_d(\omegaExplicit\dt)  \eqdef \f{\tan(\omegaExplicit\dt/2)}{\tan(\omegaExplicit\dt)}.   \label{eq:alphad}
\ea
For $\omegaExplicit\dt$ going to zero, $\alpha_d$ approaches the usual choice $\alpha=1/2$ for the continuous filter.
Substituting~\eqref{eq:WhatExplicit} into~\eqref{eq:updateDiscreteExplicit} gives
\ba
   \Vhat_m^{(k+1)} =  \beta_d(\lambdaExplicit;\,\omegaExplicit,\Ttilde,\alpha_d) \, \Vhat_m^{(k)}
       + \left( \beta_d(\omegaExplicit;\,\omegaExplicit,\Ttilde,\alpha_d) - \beta_d(\lambdaExplicit;\,\omegaExplicit,\Ttilde,\alpha_d)\right) \Uhat_m , 
       \label{eq:waveHoltzFourierFPIexplicit}
\ea
where $\beta_d$ is a discrete filter function is given by
\ba
   \beta_d(\lambda;\,\omega,T,\alpha) = \sincd(\omega+\lambda,T) + \sincd(\omega-\lambda,T) - \alpha \, \sincd(\lambda,T) ,
     \label{eq:discreteBetaThreeSincsA}
\ea
see~\ref{sec:discreteFilterFunction} for a derivation.
The function $\sincd(\lambda,T)$ in~\eqref{eq:discreteBetaThreeSincsA} is an approximate $\sinc$ function defined by
\ba
  \sincd(\lambda,T) \eqdef \f{\sin( \lambda \, T)}{T \tan(\lambda\f{\dt}{2})/(\f{\dt}{2})} ,\qquad \dt={T\over\Nt} . \label{eq:sincdA}
\ea
Noting the properties of $\beta_d$ given in~\ref{sec:discreteFilterFunction} leads to the 
following result.
\begin{theorem}[Fully Discrete Explicit WaveHoltz FPI Convergence] \label{th:FullyDiscreteExplicitWaveHoltzConvergenceRate}
  Let $\lambda_{h,m}\ne\omegaExplicit$ be the eigenvalues of the discrete problem in~\eqref{eq:discreteEigenvalueBVP}.  The asymptotic convergence rate $\ACRe$ of the fully discrete WaveHoltz fixed-point iteration 
  with explicit time-stepping is 
\ba
    & \ACRe = \max_{\lambdaExplicit} \,\bigl|  \beta_d(\lambdaExplicit;\,\omegaExplicit,\Ttilde,\alpha_d) \bigr| , 
\ea
where $\beta_d$ is the discrete beta function~\eqref{eq:discreteBetaThreeSincsA}, $\lambdaExplicit$ is defined in~\eqref{eq:discreteExplicitLambdaTilde}, $\omegaExplicit$ is
determined from Recipe~\ref{recipe:OmegaTildeExplicit} in~\ref{sec:timeCorrections}, $\Ttilde=2\pi/\omegaExplicit$, and $\alpha_d$ 
is the adjusted value for $\alpha$ given in~\eqref{eq:alphad}.
\end{theorem}
\begin{proof}
   The proof follows from the iteration~\eqref{eq:waveHoltzFourierFPIexplicit} and the assumptions in the statement of the theorem.
\end{proof}

\subsubsection{Implicit time-stepping}

Now consider the case of implicit time-stepping. The implicit scheme with corrections for 
time discretization errors is 
\bse
\label{eq:WaveSchemeImplicit}
\bat
  & \Dpt\Dmt W_\jv^n = \half L_{ph} \Big( W_\jv^{n+1} + W_\jv^{n-1}  \Big) 
     - f(\xv_\jv) \, \cos(\omegaImplicit\, t^n) \,\cos(\omegaImplicit\dt),  ~&& \jv \in \Omega_h, ~ n=0,1,2,\ldots  ,  \label{eq:WaveSchemeImplicitA} \\
  & W_\jv^0 = V_\jv,                 \qquad&& \jv \in \bar{\Omega}_h  , \\
  & \Dzt W_\jv^0 = 0,                \qquad&& \jv \in \bar{\Omega}_h  ,  \label{eq:WaveSchemeImplicitC} \\
  & \Bc_{ph} W_\jv^n = 0  , \qquad&&  \jv \in \p\Omega_h , ~ n=1,2,\ldots ,
\eat
\ese
where $\omegaImplicit$ is a discrete-correction to the frequency given by
\ba
    \omegaImplicit \eqdef \f{1}{\dt} \cos^{-1}\Big( \f{1}{1 + (\omega\dt)^2/2 } \Big) ,
    \label{eq:omegaTildeImplicitII}
\ea
see~\ref{sec:timeCorrectionsImplicit}.  
As for the explicit scheme, the initial condition~\eqref{eq:WaveSchemeImplicitC} can be combined with~\eqref{eq:WaveSchemeImplicitA}
with $n=0$ to eliminate $W_\jv^{-1}$ and arrive at an implicit update for the first time-step $W_\jv^{1}$.

Following similar steps to the analysis above for explicit time-stepping
leads to a discrete ODE for each generalized Fourier coefficient given by
\bse
\label{eq:waveFourierImplicit}
\bat
  & \Dpt\Dmt \What_m^n = \half L_{ph} \Big( \What_m^{n+1} + \What_m^{n-1} \Big) 
    - \fHat_m \, \cos(\omegaImplicit t^n) \, \cos(\omegaImplicit\dt),  \quad n=0,1,2,\ldots ,\\
  & \What_m^0 = \Vhat_m,                  \\
  & \Dzt \What_m^0 = 0 .       
\eat
\ese
The solution to~\eqref{eq:waveFourierImplicit} also takes a similar form to the continuous case and is 
\bse
\ba
   & \What_m^n = \Uhat_m \cos(\omegaImplicit t^n) + \Vhat_m \cos(\lambdaImplicit t^n) ,  
       \label{eq:WhatImplicit} \\
   & \Uhat_m \eqdef \f{\fHat_m}{\omegaImplicit^2 - \lambda_{h,m}^2} , \\
   & \lambdaImplicit \eqdef \f{1}{\dt} \cos^{-1}\left( \f{1}{1 + (\lambda_{h,m} \dt)^2/2} \right) , \label{eq:lambdaTildeImplicit}
\ea
where $\lambda_{h,m}$ are the eigenvalues of the discrete problem in~\eqref{eq:discreteEigenvalueBVP}.
Applying the discrete time filter as for the explicit time-stepping case leads to the following result.
\ese
\begin{theorem}[Fully Discrete Implicit WaveHoltz FPI Convergence Rate] 
\label{th:FullyDiscreteImplicitWaveHoltzConvergenceRate}
  Let $\lambda_{h,m}\ne\omegaImplicit$ be the eigenvalues of the discrete problem in~\eqref{eq:discreteEigenvalueBVP}.  The asymptotic convergence rate of the fully discrete WaveHoltz fixed point iteration 
  with implicit time-stepping is 
\ba
    &\ACRi = \max_{\lambdaImplicit} |  \beta_d(\lambdaImplicit;\,\omegaImplicit,\Ttilde,\alpha_d) | ,
\ea
where $\beta_d$ is the discrete beta function~\eqref{eq:discreteBetaThreeSincsA}, $\omegaImplicit$ is defined in~\eqref{eq:omegaTildeImplicitII}, $ \lambdaImplicit$ is defined in~\eqref{eq:lambdaTildeImplicit}, and $\Ttilde=2\pi/\omegaImplicit$.
\end{theorem}

{
\newcommand{\figw}{7cm}
\newcommand{\figh}{6cm}
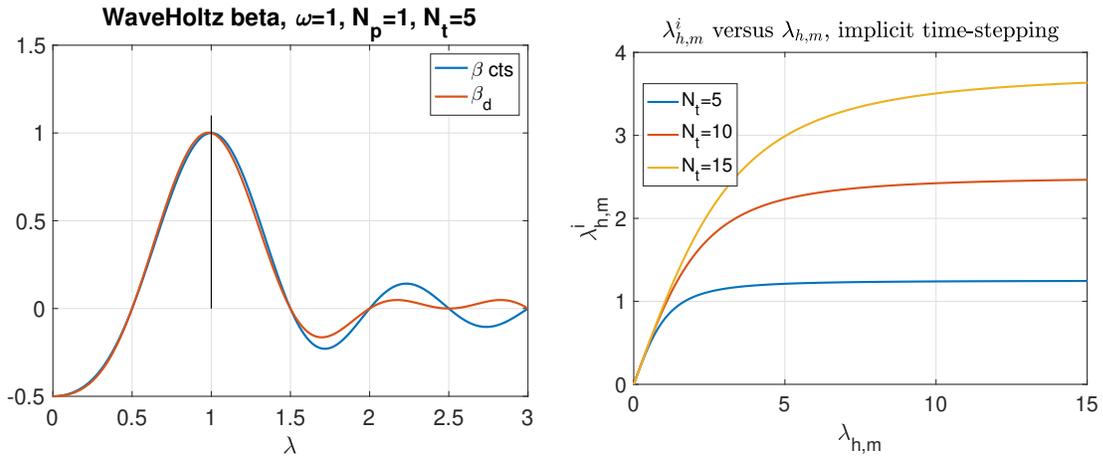
\begin{figure}[htb]
\begin{center}
\begin{tikzpicture}
  \useasboundingbox (0,.75) rectangle (2.2*\figw,1*\figh);  
  \begin{scope}[yshift=0*\figh]
    \figByWidth{0}{0}{betaFunctionNp1Nt5}{\figw}[0.][0.][0.][0.]
    \figByWidth{7.5}{0}{adjustedLambdaImplicit}{\figw}[0.][0.][0.][0.]
  \end{scope}  
\end{tikzpicture}
\end{center}
\caption{
  Left: discrete filter function $\beta_d$ and continuous filter $\beta$ for $\Nt=5$.
  Right: the adjusted $\lambdaImplicit$ for implicit time-stepping versus $\lambda_{h,m}$ for varying number of time-steps $\Nt$, $\dt=T/\dt$, for $\omega=1$.
   }
\label{fig:adjustedLambdaImplicit}
\end{figure}
}

The left graph in Figure~\ref{fig:adjustedLambdaImplicit} compares $\beta_d$ and $\beta$ for $\Nt=5$ time-steps 
per period. 
Even for this large value of $\dt$ the curves are quite similar for $\lambda$ near $\omega$.
The right graph in Figure~\ref{fig:adjustedLambdaImplicit} shows $\lambdaImplicit$ in~\eqref{eq:lambdaTildeImplicit}  
as a function $\lambda_{h,m}$ (treated as a continuous variables) for varying number of implicit time-steps $\Nt$. 
For small values of $\Nt$, the transformation~\eqref{eq:lambdaTildeImplicit} has a significant effect 
with $\lambdaImplicit$ approaching $\omega \Nt/4$ as $\lambda\rightarrow \infty$.
The astute reader will note that $\sinc_d(z,T)$ reaches a maximum of one not only for $z=0$
but also for $z \dt = 2 m \pi$.  Thus $\beta_d(\lambda)$ equals one at additional values of
$\lambda$, for example $\lambda=2\pi/\dt$.
These additional values do not
play a role in the convergence, however, 
since they occur at values of $\lambdaImplicit$ that are outside its range $\lambdaImplicit\in[0,\omega \Nt/4]$.


\section{WaveHoltz iteration: acceleration and performance} \label{sec:accelerations}

Having discussed the convergence of the WaveHoltz fixed-point iteration for both the continuous problem and its space-time approximations, we now consider approaches to accelerate convergence and improve performance.

\newcommand{\DeflationSet}{\Dc} 
\subsection{Deflation: accelerating WaveHoltz by removing some slowly converging eigenmodes} \label{sec:deflation}

The WaveHoltz fixed-point iteration (FPI) can be accelerated using a deflation approach.
Using some precomputed eigenmodes, the components of the solution along the slowest converging eigenmodes can be removed during the
WaveHoltz iteration thus improving the convergence rate. When the iteration has converged the WaveHoltz 
solution can be adjusted to include the components of the Helmholtz solution along the eigenmodes that were deflated.
A drawback of using deflation is that certain eigenmodes must be computed. However, for a given geometry, a set of eigenmodes
can be pre-computed and these can subsequently be used to solve multiple Helmholtz problems.


\renewcommand{\algFontSize}{\small}
\begin{algorithm}
\algFontSize 
\caption{WaveHoltz Algorithm with Deflation.}
\begin{algorithmic}[1]

  \Function{WaveHoltz}{$\omega$,$f$,$g$,$\Np$}  
    \State Set $\Tf=\Tbar$, where $\Tbar=\Np T$ and $T=2\pi/\omega$. 
    \State $k=0$ \Comment WaveHoltz~iteration counter.
    \State $v^{(k)}=0$   \Comment Assign initial guess for Helmholtz iterate 
    \While{ not converged} \Comment Start WaveHoltz~iterations.

      \State $w^{(k)}(\xv,0) = v^{(k)}(\xv)$ \Comment Initial condition for wave equation solve.
      \State $w^{(k)}(\xv,0:\Tf)$= \Call{solveWaveEquation}{$w^{(k)}(\xv,0)$,$f$,$g$} \Comment Solve for $\wv(\xv,t)$, $t\in[0,\Tf]$. 
      \State $\displaystyle v^{(k+1)}(\xv) = \f{2}{\Tbar} \int_{0}^{\Tbar} \left( \cos(\omega t) - \f{\alpha}{2} \right) \, w^{(k)}(\xv,t) \, dt$
           \Comment Time filter the wave equation solution.
      \State $ \displaystyle  v^{(k+1)}(\xv) =  v^{(k+1)}(\xv) - \sum_{\phi_m\in\DeflationSet} ( v^{(k+1)}, \phi_m )_\Omega \, \phi_m(\xv)$ \Comment Deflate.
      \State $k = k+1$
    \EndWhile    \Comment End WaveHoltz iterations.
    \State  $\displaystyle v^{(k)}(\xv) =  v^{(k)}(\xv) + \sum_{\phi_m\in\DeflationSet} \f{(f,\phi_m)_\Omega}{\omega^2-\lambda_m^2} \,  \phi_m(\xv)  $ \Comment Inflate.
    \State $u(\xv) = v^{(k)}(\xv)$; \Comment Approximate Helmholtz solution.
 \EndFunction
\end{algorithmic} 
\label{alg:waveHoltzWithDeflation}
\end{algorithm}

One approach to deflation is to adjust the initial conditions and forcing, iterate until convergence, and then correct the solution. 
A second approach is to not change the forcing but then deflate the WaveHoltz
solution after each iteration. 
The second approach is used here and given in Algorithm~\ref{alg:waveHoltzWithDeflation}.
Let $\DeflationSet$ denote the set of eigenmodes that are deflated.
The components of $v^{(k+1)}(\xv)$ along the eigenmodes $\phi_m \in \DeflationSet$
are removed at the end of each WaveHoltz iteration,
\ba
   v^{(k+1)}(\xv) =  v^{(k+1)}(\xv) - \sum_{\phi_m\in\DeflationSet} ( v^{(k+1)}, \phi_m )_\Omega \, \phi_m(\xv),
\ea
where $(\,\cdot\,,\,\cdot\,)_\Omega$ denotes the usual $L_2$ inner product on $\Omega$.
After the deflated WaveHoltz solution has converged, $v^{(k)}(\xv)$ is corrected by adding back the missing components of the Helmholtz solution,
\ba
  & v^{(k)}(\xv) =  v^{(k)}(\xv) + \sum_{\phi_m\in\DeflationSet} \f{(f,\phi_m)_\Omega}{\omega^2-\lambda_m^2}\,  \phi_m(\xv) .
\ea
Note that here we have assumed that $g=0$; further adjustments would be needed for inhomogeneous boundary conditions.
The convergence rate of the deflated algorithm follows easily.
\begin{theorem}[WaveHoltz FPI convergence with deflation]
  Assume $\lambda_m\ne\omega$ are the eigenvalues of the problem in~\eqref{eq:eigBVP} with corresponding eigenfunctions $\phi_m(\xv)$.  The asymptotic convergence rate $\mu$ of the WaveHoltz fixed-point iteration with deflation is 
\ba
    \mu = \max_{\phi_m \not\in \DeflationSet} |  \beta(\lambda_m) | ,
\ea
where $\beta(\lambda)$ is given by~\eqref{eq:filterFunction} and $\DeflationSet$ denotes the set of deflated eigenmodes.
\end{theorem}
In practice the deflation set $\DeflationSet$ should normally be chosen to include eigenmodes whose eigenvalues are closest to $\omega$.  Ideally one would want to deflate enough eigenmodes so that the WaveHoltz FPI convergence rate is less than some specified value, e.g.~$\beta_0=0.7$.
Note, however, that the number of eigenmodes needed for deflation increases in proportion to $\omega^{2\nd}$ for $\nd$ space dimensions~\cite{WaveHoltz1},
and thus many eigenmodes are generally needed for large $\omega$ in three dimensions.

The implementation of deflation on a overset grid requires the calculation of discrete eigenvalue/eigenvector pairs on an overset grid.  
We perform this calculation using a Krylov-Schur algorithm from SLEPSc~\cite{SLEPc2005}.
The discrete approximation to the eigenvalue problem on an overset grid consists of approximations to the PDE and boundary
conditions together with interpolation equations. This is a generalized eigenvalue problem of the form $A x = \lambda B x$, since the eigenvalue does not
appear in the boundary conditions and interpolation equations. The matrix $B$ has ones on the diagonal for points where the PDE is discretized and zeros for constraint equations.
It is possible, in principle, to eliminate all constraint equations
and reduce the problem to a regular eigenvalue problem of the form $A x = \lambda x$ for a reduced matrix $A$.
For practical reasons, however, it is convenient to retain the constraint equations.
The algorithms in SLEPSc seem to work best if the matrix $B$ in the generalized form is nonsingular.
In the overset grid setting $A$ is nonsingular while $B$ is singular. To resolve this issue, the roles of $A$ and $B$ can be reversed and instead we solve a related generalized eigenvalue problem $B x = (1/\lambda) A x $ for the reciprocals of the eigenvalues.
The eigenvectors returned from SLEPSc are normalized using the discrete inner product. 
For any multiple eigenvalues, an orthonormal basis for the corresponding eigen-space is found.  Following this procedure, we are able to compute the required discrete eigenvalue/eigenvector pairs to carry out the deflation algorithm described in Algorithm~\ref{alg:waveHoltzWithDeflation}.

It should be noted, however, that computation of the eigenmodes using SLEPSc requires the inversion of a large (often indefinite) matrix,
and generally we use a direct sparse solver to do this. 
This can be expensive for large problems. However, the eigenpairs can be computed as a pre-processing step and used to solve multiple Helmholtz problems.  It turns out that the WaveHoltz algorithm can be used to compute eigenpairs
without the need to an indefinite matrix; this will be described in a forthcoming paper.


\subsection{Implicit time-stepping with a very large time-step} \label{sec:implicitTimeStepping}

There is a potential to dramatically improve
the run-time performance of the WaveHoltz algorithm
through the use of implicit time-stepping and a large time-step. 
A key result of using implicit time-stepping 
is that a small number of time-steps per period can be taken, this number being independent
of the mesh spacing or order of accuracy in space. Thus, as the mesh is refined, the total
number of time-steps needed to reach convergence should be independent of the mesh spacing 
(see Section~\ref{sec:optimal} for further details).
This is in contrast to explicit time-stepping where a stability constraint on the time-step size forces the number of time-steps to increase as
$1/h$ when the mesh spacing $h$ decreases.

We adopt a variation of the implicit, high-order accurate modified equation time-stepping schemes for the wave equation
developed in~\cite{carson2024highorderaccurateimplicitexplicittimestepping}.  The present schemes use high-order accuracy in space but only second-order accuracy in time.
Second-order accuracy in time is used since time-discretization errors can be removed (see Section~\ref{sec:timeCorrectionsImplicit})
from the WaveHoltz solution.
The implicit time-stepping scheme for a grid function $W_\iv^n \approx w(\xv_\iv,t^n)$
takes the form 
\ba
  \Dpt\Dmt W_\iv^n = L_{ph}  \left[ \frac{1}{2} W_\iv^{n+1} + \frac{1}{2} W_\iv^{n-1} \right]  + F(\xv_i,t^n) ,
  \label{eq:implicitSchemeI}
\ea
which uses a second-order accurate (\textit{trapezodial}) in time approximation 
and a $p^{\rm th}$-order accurate spatial approximation $L_{ph}$ (see Section~\ref{sec:numericalApproach} for further details on the spatial discretization).
The linear system that needs to be inverted with implicit time-stepping is well suited
to be solved by fast methods such as multigrid~\cite{multigridWithNonstandardCoarsening2023}.

When using a large time-step $\dt$ it is important to correct for time-discretization errors, and,
as derived in~\ref{sec:timeCorrectionsImplicit}, the minimum number of time-steps per (smallest) period,
denoted by $\NITS$,
is then given by
\ba
    \NITS \ge 5.
\ea
The convergence rate of the scheme depends on the value of $\Nits$ since the time-step depends on $\Nits$ and this affects the discrete filter function $\beta_d$.
 Larger values of $\Nits$ may give faster convergence (to a point) but
at a larger computational cost. In practice we have found that a value of $\Nits=10$ is often a good comprise, although this could be problem dependent.
It should also be noted that when using a large implicit time-step, 
the first time-step should also be implicit otherwise the WaveHoltz FPI iteration may fail to converge properly.
See~\cite{carson2024highorderaccurateimplicitexplicittimestepping} for details of the form of the 
implicit first time-step.


\subsection{Krylov methods} \label{sec:krylov}

\newcommand{\OpWH}{\Wc}

The WaveHoltz fixed-point iteration (FPI), with or without deflation, can be accelerated with Krylov methods. 
Consider the FPI for the continuous problem which takes the general form 
\ba
     v^{(k+1)}  = \OpWH( v^{(k)},f ),
\ea
where $\OpWH$ is the affine operator that takes $v^{(k)}$ as initial condition to the wave equation
and returns $v^{(k+1)}$ as the next iterate.  (The dependence of $v^{(k)}$, $v^{(k+1)}$ and $f$ on the independent variable $\xv$ is suppressed.)
This function can be written in the form
\ba
   & v^{(k+1)}  = \OpWH( v^{(k)},f ) = S v^{(k)} + b(f)   \label{eq:waveHoltzMatrixIterationForKrylov}
\ea
where $S$ is a linear operator and 
the function $b$ (independent of $k$) is simply $\OpWH(0,f)$, i.e.~the result of one WaveHoltz iteration starting 
from a zero initial condition, $v^{(0)}=0$.
For the class of problems considered in this article, 
the operator $S$ is self-adjoint with real eigenvalues $\beta(\lambda_m)$ as described in Section~\ref{sec:convergenceAnalysis}.
Thus the WaveHoltz time filter has transformed the eigenvalues of the Helmholtz operator to $\beta(\lambda_m)\in[-\half,1]$ assuming $\alpha=\half$ is used in the filter.

The solution to the discretized fixed-point iteration can be found directly by solving
the linear system
\ba
   A \vv_h \eqdef (I-S_h) \vv_h = \bv_h, \label{eq:WaveHoltzLinearSystem}
\ea
where $S_h$, $\vv_h$ and $\bv_h$ denote the discrete approximations to $S$, $v$, and $b$, respectively.
Forming the matrix explicitly would be expensive for large problems and so instead a matrix-free iterative method,
such as a Krylov method, is used.
A matrix-free method requires a function that can evaluate $A \vv_h$ for any vector $\vv_h$.
From a discrete approximation of~\eqref{eq:waveHoltzMatrixIterationForKrylov}, 
$S_h \vv_h$ can be computed by applying one WaveHoltz iteration with initial condition $\vv_h$ and then subtracting $\bv_h$,
\ba
    S_h \vv_h \eqdef \OpWH_h(\vv_h,\fv_h) - \bv_h,
\ea
where $\OpWH_h$ denotes the discrete version of $\OpWH$.
Whence $A\vv_h$ can be evaluated using
\ba
   A \vv_h = \vv_h - \OpWH_h(\vv_h,\fv_h) + \bv_h. 
\ea

\mni
Note from~\eqref{eq:waveHoltzMatrixIterationForKrylov} that $\OpWH_h(\zerov,\fv_h)=\bv_h$ and $\OpWH_h(\vv_h,\zerov) = S_h\vv_h$, and thus 
$A\vv_h$ can also be evaluated using
\ba
     A \vv_h = \vv_h - \OpWH_h(\vv_h,\zerov).
     \label{eq:MevalII}
\ea
Using~\eqref{eq:MevalII} implies that, for the matrix-vector product, 
the discrete wave equation solve can be performed with zero forcing, which can
provide a computational saving.

{
\newcommand{\figw}{7cm}
\newcommand{\figh}{6cm}
\begin{figure}[htb]
\begin{center}
\begin{tikzpicture}
  \useasboundingbox (0,.65) rectangle (\figw,.95*\figh);  
  \begin{scope}[yshift=0*\figh]
    \figByWidth{0}{0}{waveHoltzOneMinusBetaFunction}{\figw}[0.][0.][0.][0.]
  \end{scope}  
\end{tikzpicture}
\end{center}
\caption{
   Plots of $1-\beta(\lambda)$ for $\Np=1$, $\Np=2$, and $\Np=3$ periods per time-interval.
   The Krylov solvers operate on the matrix $A=I-S_h$ which has eigenvalues $1-\beta(\lambda_{h,m})$ of which representative values are shown with red x's.
   }
\label{fig:waveHoltzOneMinusBeta}
\end{figure}
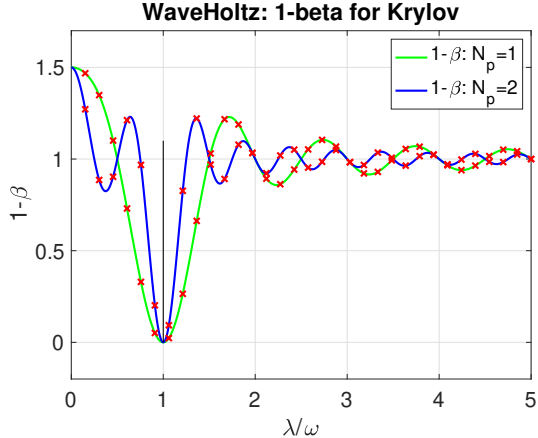
}

We note that for an overset grid, the matrix $A$ is not symmetric in general, and thus Krylov methods
appropriate for non-symmetric matrices such as GMRES should be used.
In fact, the results given in later sections show that GMRES can be quite effective.
The matrix $A=I -S_h$ has eigenvalues $1-\beta(\lambda_{h,m})$. 
Figure~\ref{fig:waveHoltzOneMinusBeta} shows plots of $1-\beta$ for $\Np=1$ and $\Np=2$ along with representative eigenvalues marked as red x's.
In general there will be many eigenvalues of $A$ clustered near $1$ as well as small eigenvalues 
near where $\lambda_{h,m} \approx \omega$. 
GMRES finds the ``best'' solution in the Krylov space spanned by the WaveHoltz iterates $S^k \vv^{(0)}$, and it uses
an Arnoldi process to form an orthonormal basis for this vector space. This Arnoldi process identifies
the most slowly converging eigenvectors through a power-method-type iteration.
When the slowly converging eigenvectors are found, their contribution to the solution can be removed. GMRES is also good at detecting clusters of eigenvalues and so should be effective at identifying the eigenvalues of $A$ near $1$.
Thus if there are just a few slowly converging eigenmodes, or the slowing converging ones have been deflated then GMRES should converge very fast.


\section{Pollution errors and points-per-wavelength rules of thumb} \label{sec:pollutionErrors}

It is well known, see for example~\cite{BaylissGoldsteinTurkel1985,IlenburgBabuska1995}, that for large frequencies $\omega$, or large wave numbers $k = \omega/c$,
discrete solutions to the Helmholtz equation suffer from pollution or dispersion errors 
where, for accuracy, the number of points-per-wavelength (PPW) must increase with increasing $\omega$.
In this section we present a model problem that provides theoretical insight into the source of the pollution errors. We are then able to present a simple rule of thumb
that can be used as a rough estimate for choosing the PPW.

\subsection{Helmholtz model problem and pollution errors}  \label{sec:HelmholtzModelProblem}

A simple way to see the source of pollution (dispersion) errors is to consider a model Helmholtz BVP on the interval $x\in[a,b]$ given by
\bse
\label{eq:pollutionModelProblem}
\bat
    & \p_x^2 u + k^2 \, u = \cos(\kappa x), \qquad && x\in(a,b) , \\
    & u(a)=0, \quad u(b)=0 ,
\eat
\ese
where $k=\omega/c>0$ is a wave number, $\kappa>0$ is a given constant, and $L=b-a$ is the length of the problem domain. 
The forcing term $\cos(\kappa x)$ can be thought of as one term in a Fourier expansion of a more general forcing.
 We assume that $\kappa\ne k$ so that the harmonic forcing is not resonant and that $\sin(kL)\ne0$ so that the BVP is nonsingular and a unique solution exists, i.e.~$k$ is not an eigenvalue of the associated eigenvalue problem.
The solution to~\eqref{eq:pollutionModelProblem} is the sum of a particular solution (of the forced problem)  and a homogeneous solution,
\bse
     \label{eq:modelSolutionContinuous}
\ba
     u(x) = \up(x) + u^h(x),
\ea
where
\ba
  \up(x)    &= {\cos(\kappa x)\over k^2 - \kappa^2}, \label{eq:modelSolutionContinuousP}\\
  u^{h}(x)  &= -\up(b){\sin(k(x-a))\over\sin(kL)}-\up(a){\sin\bigl(k(b-x)\bigr)\over\sin(kL)}  . \label{eq:modelSolutionContinuousH}
\ea
\ese
A second-order accurate discrete approximation of~\eqref{eq:pollutionModelProblem} is
\bse
\label{eq:pollutionModelProblemDiscreteP}
\bat
  & \Dpx\Dmx U_j + k^2 \, U_j = \cos(\kappa x_j), \qquad && j=1,2,\ldots,N-1, \\
  & U_0 = 0 , \quad U_N=0 ,
\eat
\ese
where $x_j=a + j\dx$ and $\dx=L/N$.
The solution to~\eqref{eq:pollutionModelProblemDiscreteP} also takes the form of a particular solution
plus a homogeneous solution, 
\bse
    \label{eq:modelSolutionDiscrete}
\ba
     U_j = U^f_j + U^h_j,
\ea
where
\ba
\Up_j&={\cos(\kappa x_j)\over k^2 - \kappaTilde^2}, \label{eq:modelSolutionDiscreteP}\\
    U_j^{h}&=-\Up_N\,{\sin(\kTilde (x_j-a))\over\sin(\kTilde L)}-\Up_0\,{\sin\bigl(\kTilde(b-x_j)\bigr)\over\sin(\kTilde L)}. \label{eq:modelSolutionDiscreteH}
\ea
\ese
Here, $\kTilde$ and $\kappaTilde$ are related to $k$ and $\kappa$, respectively, through
\ba
    \f{\sin(\kTilde\dx/2)}{\dx/2} = k , \qquad 
    \f{\sin(\kappa\dx/2)}{\dx/2} = \kappaTilde  .
    \label{eq:kandbtilde}
\ea
For later purposes, we note that for $\kappa\dx$ and $k\dx$ approaching zero, we have
\ba
  \kTilde=k\left[1+{1\over24}(k\dx)^2+O\bigl((k\dx)^4\bigr)\right],\qquad 
  \kappaTilde=\kappa\left[1-{1\over24}(\kappa\dx)^2+O\bigl((\kappa\dx)^4\bigr)\right].
  \label{eq:kandbetaexpansions}
\ea

We are interested in the relative
 error between the discrete solution $U_j$ in~\eqref{eq:modelSolutionDiscrete} and the continuous solution $u(x)$ in~\eqref{eq:modelSolutionContinuous} at $x=x_j$.  Define this error as
\ba
  E_j\eqdef {\bigl\vert U_j-u(x_j)\bigr\vert\over \Nc},\qquad \Nc={1\over \bigl\vert k^2-\kappa^2\bigr\vert}\,{1\over\bigl\vert\sin(kL)\bigr\vert}.
\ea
where $\Nc$ scales the error by the size of the homogeneous solution.  Using the triangle inequality, we have
\bse
\ba
E_j\le \Ef_j+E_j^{h},
\ea
where
\ba
\Ef_j\eqdef{\bigl\vert \Up_j-\up(x_j)\bigr\vert\over \Nc},\qquad E_j^{h}\eqdef{\bigl\vert U_j^{h}-u^{h}(x_j)\bigr\vert\over \Nc},
\label{eq:errorContributions}
\ea
\ese
are contributions to the scaled error from the errors in the particular solution and the homogeneous solution. 

\mni
\textbf{Error in the particular solution.}
First consider bounding $\Ef_j$. 
Using the expressions for $\up(x)$ and $\Up_j$ in~\eqref{eq:modelSolutionContinuousP} and~\eqref{eq:modelSolutionDiscreteP}, respectively, we find
\ba
  \Ef_j={1\over\Nc}\left\vert{\cos(\kappa x_j)\over k^2-\kappaTilde^2}-{\cos(\kappa x_j)\over k^2-\kappa^2}\right\vert
       = |\sin(kL)| |\cos(\kappa x_j)| \left\vert{\kappaTilde^2-\kappa^2\over k^2-\kappaTilde^2}\right\vert
       \le \left\vert{\kappaTilde^2-\kappa^2\over k^2-\kappaTilde^2}\right\vert.
\ea
Using the expansion for $\kappaTilde$ in~\eqref{eq:kandbetaexpansions} gives
\bse
\ba
    \Ef_j\le {K_f\kappa^2\over\bigl\vert k^2-\kappa^2\bigr\vert} \, (\kappa\dx)^2  \\
    K_f \eqdef {1\over12}+O\bigl((\kappa\dx)^2\bigr),
\label{eq:polutionPartErrorFinal}
\ea
\ese
assuming $\kappa\dx$ is small.
The particular solution we have chosen thus has a relative error proportional to $(\kappa\dx)^2$.
Even if $\kappa=O(k)$, this error would be controlled provided $k\dx$ is small, which is a usual accuracy requirement based on points-per-wavelength.  Of course, $\Ef_j$ can be large if $\kappa\approx k$, but this is expected for a harmonic forcing near resonance.

\mni
\textbf{Error in the homogeneous solution.}
Due to our choice for the particular solution, it is the error in the discrete homogeneous solution that is the source of the pollution error.  
Substituting the expressions for $u^{h}(x)$ and $U_j^{h}$ in~\eqref{eq:modelSolutionContinuousH} and~\eqref{eq:modelSolutionDiscreteH}, respectively, 
into the expression for $E_j^{h}$ in~\eqref{eq:errorContributions} gives
\bse
\label{eq:polutionHomogError}
\ba      
  & E_j^{h} \le E_j^a + E_j^b,  \\
  \intertext{where}
  & E_j^a \eqdef {1\over\Nc}\left\vert \,\Up_N \,{\sin(\kTilde (x_j-a))\over\sin(\kTilde L)}- \up(b) {\sin(k (x_j-a))\over\sin(kL)}\right\vert  ,\label{eq:polutionHomogErrora} \\
  & E_j^b \eqdef  {1\over\Nc}\left\vert \,\Up_0 \,{\sin\bigl(\kTilde(b-x_j)\bigr)\over\sin(\kTilde L)}- \up(a) {\sin\bigl(k(b-x_j)\bigr)\over\sin(kL)}\right\vert .
        \label{eq:polutionHomogErrorb}
\ea
\ese
Let us obtain a bound for $E_j^a$ (the bound for $E_j^b$ is similar).
Substituting the expressions for $\Up_N$, $\up(b)$, and  $\Nc$ into~\eqref{eq:polutionHomogErrora} gives
\bse
\label{eq:polutionHomogErroraMod}
\ba
  & E_j^{a} = | \cos( b L ) | \,
       \big|        \Ac \sin{(\kTilde(x_j-a))} - \sin(k (x_j-a))    \big|, \\
  & \Ac \eqdef \f{k^2-\kappa^2}{k^2-\kappaTilde^2} \, \f{\sin(kL)}{\sin(\kTilde L)} .  
\ea
\ese
We identify the term $\Ac$ as the ratio of a discrete amplitude to the corresponding continuous one, and this ratio can be written as one plus a correction,
\bse
\ba
   & \Ac = 1 + \Eamp, \\
   & \Eamp \eqdef \f{k^2-\kappa^2}{k^2-\kappaTilde^2} \, \f{\sin(kL)}{\sin(\kTilde L)} -1 ,  \label{eq:ampErrorb}
\ea
\ese
where $\Eamp$ denotes the relative error in the amplitude of the discrete homogeneous solution.
The relative phase error from the expression for $E_j^{a}$ in~\eqref{eq:polutionHomogErroraMod} is identified as 
\ba
  \Ephase \eqdef \big| \sin({\kTilde(x_j-a))} - \sin(k (x_j-a))  \big| . \label{eq:phaseError}
\ea

Consider first the relative amplitude error $\Eamp$. 
Assuming $k\dx \ll 1$ and $\kappa\dx \ll1$ and using~\eqref{eq:kandbetaexpansions}, we have the following expansions
\bse
\ba
   {k^2-\kappa^2\over k^2-\kappaTilde^2}&=1-K_f{\kappa^2\over k^2-\kappa^2}(\kappa\dx)^2,\\
   {\sin(kL)\over\sin(\kTilde L)}&=1-K_h{k L\over \tan(k L)}(k\dx)^2,
\ea
\ese
where $K_f$ is defined in~\eqref{eq:polutionPartErrorFinal} and $K_h$ is given by
\ba
   K_h \eqdef{1\over 24}+O\bigl((k\dx)^2\bigr).
\ea
Substituting these expressions into~\eqref{eq:ampErrorb} leads to the following result.

\mni
\textbf{Summary (amplitude error).}
The amplitude error has the bound 
\shadedBoxWithShadow{align}{blue}{
   \Eamp  \le K_h\,  \f{k L }{|\tan(k L)|} \, (k\dx)^2 + K_f{1\over |(k/\kappa)^2- 1| }(\kappa\dx)^2 .
}
We observe that if  $|\tan(kL)|$ is not too small and if $k/\kappa$ is not too close to one, then the amplitude error can be controlled by 
making $ kL \, (k\dx)^2$ and $(\kappa\dx)^2$ small.
Note the factor $k L$ multiplying $(k\dx)^2$ in this requirement and this can be one source of pollution error if $kL$ is large.
On the other hand, 
suppose that $k$ is close to an eigenvalue $k_m = m\pi/L$ for $m=1,2,3,\ldots$,
\ba 
   k= k_m  +  \deltak,
\ea
then
\ba
   \f{k L }{\tan(k L)} = \f{k_m L + \deltak\, L}{\tan(k_mL+\deltak\, L)} = \f{k_m}{\deltak} + 1 + (k_m L) \, O( \deltak L ) , \qquad \deltak\,L \ll1. 
\ea 
\mni
\textbf{Summary (amplitiude error for $k$ near an eigenvalue).}
When $k$ is close to an eigenvalue the relative error in the amplitude is
\shadedBoxWithShadow{align}{blue}{
    \Eamp =  K_h\, \f{k_m }{\vert\deltak\vert}  \, (k\dx)^2  + O(1) , \qquad \deltak\,L \ll1 , \label{eq:AmpErrSingular}
}
which scales as $k_m/\vert\delta k\vert$, the inverse of the relative distance between $k$ and the eigenvalue $k_m$.
The relative error in the amplitude can thus be large when $k$ is very close to an eigenvalue.

\medskip
Now consider the contribution of the relative phase error given by $\Ephase$ in~\eqref{eq:phaseError}. 
Using
\ba
   \sin(\kTilde (x_j-a))&=\sin(k(x_j-a)) + K_h\, k (x_j-a) \, \cos(k(x_j-a)) \, (k\dx)^2,
\ea
leads to the following.

\mni
\textbf{Summary (phase error).}
\shadedBoxWithShadow{align}{blue}{
   \Ephase \le  K_h\, k \, |x_j-a| \, (k\dx)^2 \le K_h (k L ) \, (k\dx)^2 .   \label{eq:finalPhaseError}
}
Note that the phase error in~\eqref{eq:finalPhaseError} also scales as $k L$; this is another source of pollution error.

To obtain a rule-of-thumb to guide a choice for $\dx$ to suppress pollution errors, we note that the contributions to the relative error in the discrete solution from either the amplitude error or phase error are dominated by terms involving $kL(k\dx)^2$.  While the amplitude error can be large in our model if $\tan(kL)$ is small or $k/\kappa$ is close to one, we ignore these factors for the purposes of this rule-of-thumb.  With this assumption, the dominant contribution to the amplitude error and phase error are similar, and we define
\shadedBoxWithShadow{align}{orange}{
    \Ep_2 \eqdef \f{1}{24} \, k L \, (k\dx)^2  ,  \label{eq:pollutionErrorOrder2}
}
as the approximate relative error of the second-order accurate scheme.
It should be remembered, however, that a grid spacing requirement derived from~\eqref{eq:pollutionErrorOrder2} is just a first guess;
finer grids may be needed such as for problems that are forced close to resonance or problems with $k$ is close to an eigenvalue.

\medskip
Now consider solving the model problem~\eqref{eq:pollutionModelProblem} to order of accuracy $p$, where $p=2,4,6,\ldots$.
We suppress the details of the discrete solution and instead focus on the key ingredients that lead to the form of the phase error as a guide to extend the order-two formula in~\eqref{eq:pollutionErrorOrder2}.
A $p^{\rm th}$-order accurate approximation to the second derivative can be written in the form 
\ba
   \p_x^2 & \approx  \Dpx\Dmx \sum_{\mu=0}^{p/2} b_{\mu} (-\dx^2 \Dpx\Dmx)^\mu .  \label{eq:SecondDerivOrderp}
\ea
The following theorem, proved in~\ref{sec:discreteDispersion}, gives a simple closed form expression for the coefficients $b_\mu$.
\begin{theorem} \label{th:uxxCoeff}
  The coefficients $b_{\mu}$ in the difference approximation~\eqref{eq:SecondDerivOrderp} for the second derivative are
\ba
    b_{\mu}  = \f{2}{ (\mu+1)^2 \, { 2\mu+2 \choose \mu+1} } = \f{2\, (\mu !)^2}{ (2\mu+2)!} , \qquad \mu=0,1,2,\ldots ~.  \label{eq:uxxCoeff}
\ea
\end{theorem}
In particular the first few coefficients are 
\ba
  b_0=1, \quad b_1=\f{1}{12}, \quad b_2=\f{1}{90}, \quad b_3=\f{1}{560}, \quad b_4=\f{1}{3150}.
\ea
To our knowledge this is the first time the explicit formula~\eqref{eq:uxxCoeff} has been presented.

Using~\eqref{eq:SecondDerivOrderp} implies the discrete wave-number $\kTilde$ satisfies, 
\ba
   k^2 = \f{4 \sin^2(\kTilde\dx/2)}{\dx^2} \sum_{\mu=0}^{p/2} b_{\mu} (- 4\sin(\kTilde\dx/2)^2)^\mu. \label{eq:kTildeSymbolOrderp}
\ea
The phase error for a $p^{\rm th}$-order accurate central difference scheme has, to leading order, the same
form as~\eqref{eq:phaseError} except with $\kTilde$ from~\eqref{eq:kTildeSymbolOrderp} satisfying 
\ba
    \kTilde = k \, \Big[ 1 +  \half b_{p/2} \,  (k\dx)^{p} +  O\bigl( (k\dx)^{p+2} \bigr) \Big]. \label{eq:kErrorOrderp}
\ea
See~\ref{sec:discreteDispersion} for the derivation of~\eqref{eq:kErrorOrderp}.

\mni
\textbf{Summary (phase error at order $p$).}
Following a similar argument to the second-order accurate case, $p=2$, the relative error in the phase for a $p^{\rm th}$-order accurate scheme
is approximated by 
\shadedBoxWithShadow{align}{orange}{
   \Ep_p \eqdef  \half b_{p/2} \, k L \, (k\dx)^p  , \quad p=2,4,6,\ldots.  \label{eq:pollutionErrorOrderp}
}

\subsection{Rules of thumb for choosing the point-per-wavelength, $\PPW$}

We now derive a rule-of-thumb estimate, based on~\eqref{eq:pollutionErrorOrderp}, that can be used to estimate approximately how fine the grid spacing should be to manage pollution errors.
Note the following relations between the angular frequency $\omega$, the wave-number $k$, the wave-length $\Lambda$, the
grid spacing $\dx$, and the points-per-wavelength, $\PPW$:
\ba
    \f{\omega}{c}=k, \qquad \Lambda = \f{2\pi}{k}, \qquad \PPW = \f{\Lambda}{\dx} = \f{2\pi}{k \dx} .
\ea
Given a relative error tolerance $\eps$, set $\Ep_p=\eps$, and then re-arrange~\eqref{eq:pollutionErrorOrderp}
to give
\ba
     \f{1}{(k\dx)^p} = \half b_{p/2} \, k L \, \f{1}{\eps}.  \label{eq:EpA}
\ea
Taking the $p\,$-th root of~\eqref{eq:EpA} and multiplying by $2\pi$ gives 
\ba
   \PPW = \f{2\pi}{k \dx} = 2\pi \, (\half b_{p/2})^{1/p} \, \left[ \f{ k L}{\eps} \right]^{1/p}. 
\ea
Define $N_\Lambda$ to be the size of the domain (largest length in any direction) in wave-lengths,
\ba
   N_\Lambda  \eqdef \f{L}{\Lambda},
\ea
and note that $k L = 2\pi L/\Lambda=2\pi N_\Lambda$ leads to the rule of thumb in Recipe~\ref{recipe:PPW}.

{
\newcommand{\figSize}{7cm}

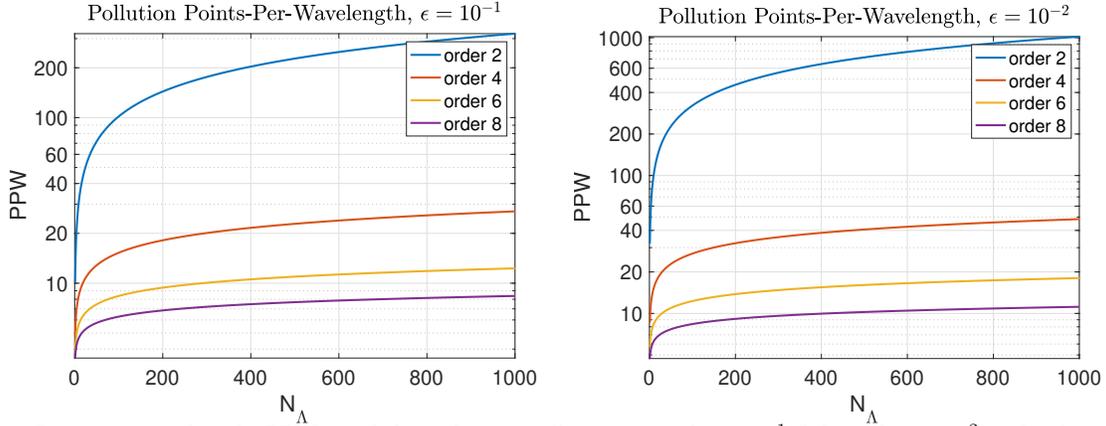
\begin{figure}[htb]
\begin{center}
\begin{tikzpicture}[scale=1]
  \useasboundingbox (0,.8) rectangle (14.5,5.7);  

  \begin{scope}[yshift=0cm]
    \figByWidth{0.0}{0}{pointsPerWaveLengthTol1}{\figSize}[0][0][0][0]
    \figByWidth{7.5}{0}{pointsPerWaveLength}{\figSize}[0][0][0][0]
  \end{scope}  

\end{tikzpicture}
\end{center}
\caption{Points per wavelength ($\PPW$) needed to achieve a pollution error of $\eps=10^{-1}$ (left) and $\eps=10^{-2}$ (right) for a domain of length $N_\Lambda$ wavelengths for schemes of
  different orders of accuracy. 
   The formula for $\PPW$ is given by~\eqref{eq:ppwPollution} in Recipe~\ref{recipe:PPW}.
    } 
\label{fig:pointsPerWavelength}
\end{figure}
}

\begin{table}[htb]\tableFont 
\begin{center}
\begin{tabular}{|c|c|c|c|c|c|} \hline 
  \multicolumn{2}{|c}{} &   \multicolumn{4}{|c|}{Points-per-wavelength}    \\ 
     $\eps$     &  $N_\Lambda$ &  $p=2$   &  $p=4$    & $p=6$  & $p=8$  \\ \hline 
 $10^{-1}$ & $    1$ & $   10$  & $    5$  & $    4$  & $    4$ \\
 $10^{-1}$ & $   10$ & $   32$  & $    9$  & $    6$  & $    5$ \\
 $10^{-1}$ & $  100$ & $  102$  & $   15$  & $    8$  & $    6$ \\
 $10^{-1}$ & $ 1000$ & $  321$  & $   27$  & $   12$  & $    8$ \\ \hline 
 $10^{-2}$ & $    1$ & $   32$  & $    9$  & $    6$  & $    5$ \\
 $10^{-2}$ & $   10$ & $  102$  & $   15$  & $    8$  & $    6$ \\
 $10^{-2}$ & $  100$ & $  321$  & $   27$  & $   12$  & $    8$ \\
 $10^{-2}$ & $ 1000$ & $ 1017$  & $   48$  & $   18$  & $   11$ \\ \hline 
 $10^{-3}$ & $    1$ & $  102$  & $   15$  & $    8$  & $    6$ \\
 $10^{-3}$ & $   10$ & $  321$  & $   27$  & $   12$  & $    8$ \\
 $10^{-3}$ & $  100$ & $ 1017$  & $   48$  & $   18$  & $   11$ \\
 $10^{-3}$ & $ 1000$ & $ 3215$  & $   86$  & $   26$  & $   15$ \\
  \hline
\end{tabular}
\caption{Approximate number of points-per-wavelength required to reach a relative error tolerance $\eps$ for a domain of longest dimension $N_\Lambda$ wavelengths.
 These values come from the formula~\eqref{eq:ppwPollution} for $\PPW$ in Recipe~\ref{recipe:PPW}.
    }
\label{tab:ppw}
\end{center}
\end{table}

\begin{recipe}[Rule of thumb for choosing the points-per-wavelength]
\label{recipe:PPW}
Given a relative error tolerance $\eps$, choose the number of points-per-wavelength for a $p^{\rm th}$-order accurate scheme from
\shadedBoxWithShadow{align}{orange}{
    \PPW_p \eqdef  2 \pi \, (\pi \, b_{p/2})^{1/p} \, \left[ \f{N_\Lambda}{\eps} \right]^{1/p} , \label{eq:ppwPollution}
}
where $N_\Lambda$ is the size of the domain in wave-lengths and $b_{p/2}$, given by~\eqref{eq:uxxCoeff}, is the coefficient in the expansion~\eqref{eq:SecondDerivOrderp} for the order $p$ 
approximation to the second derivative.
\end{recipe}
Note that for $p=2,4,6,8$, the values of the factor $ (\pi\,b_{p/2})^{1/p}$ appearing in~\eqref{eq:ppwPollution} are similar in size,
\ba
  (\pi\, b_1)^{1/2} \approx 0.51,  \quad 
  (\pi\, b_2)^{1/4} \approx 0.43, \quad 
  (\pi\, b_3)^{1/6} \approx 0.42,\quad 
  (\pi\, b_4)^{1/8} \approx 0.42 ~. 
  \label{eq:kErrorCpRoot}
\ea
While the numerical values shown in~\eqref{eq:kErrorCpRoot} appear to settle to $0.42$ as $p$ increases, 
it can be shown from~\eqref{eq:uxxCoeff} that $(\pi b_{p/2})^{1/p}\rightarrow1/2$ as $p\rightarrow\infty$.  This limit, in turn, implies that $\PPW_p$ approaches~$\pi$ for large~$p$.

The rule of thumb in Recipe~\ref{recipe:PPW} is in agreement with the well known result for the accuracy 
of discrete approximations to wave propagation problems~\cite{KreissOliger1972,GustafssonKreissOliger95} that 
the points-per-wavelength should be taken proportional to $(N_{\rm periods}/\eps)^{1/p}$, where $t=TN_{\rm periods}$ is the time in multiples of the period~$T$ over which the wave has traveled, and $\eps$ is the relative error tolerance. 
\textit{The pollution error for the Helmholtz problem thus corresponds to the accumulated dispersion error of a wave that has traveled across the domain.}

\medskip
Figure~\ref{fig:pointsPerWavelength} graphs $\PPW_p$ as a function of $N_\Lambda$ for $\eps=10^{-1}$ and $10^{-2}$ and $p=2,4,6,8$.
All curves show an initial rapid increase in the points-per-wave-length and then asymptote to a slower growth as a function of $N_\Lambda$.
The curves for $p=2$ are seen to quickly grow to very large values of $\PPW$. As the order of accuracy increases the growth in $\PPW$ is 
much more gradual.
Table~\ref{tab:ppw} gives some representative values for different values of $\eps$, $N_\Lambda$ and $p$ (note that the $\PPW_p$ depends
on the ratio $N_\Lambda/\eps$ which explains the repeated values).
For example, for a tolerance of $\eps=10^{-2}$ and a domain $N_\Lambda=100$ wavelengths in size, $\PPW_2=321$, $\PPW_4=27$, $\PPW_6=12$, and $\PPW_8=8$.
The second-order accurate scheme thus requires a massive $321$ points per-wavelength; without an understanding of pollution errors, a grid with this spacing 
would seem, at first glance, to be highly over-resolved.
The fourth-order accurate scheme with $\PPW_4=27$ has much more reasonable resolution requirements, while the sixth ($\PPW_6=12$) and 
eighth order ($\PPW_8=8$) schemes are even better.
It is apparent that for large values of $N_\Lambda/\eps$, the use of high-order accurate schemes is generally advisable.


\section{Numerical Results} \label{sec:numerics}

This section presents numerical results that demonstrate the properties of the WaveHoltz algorithm for overset grids.
The computations are performed at second-order and fourth-order accuracy in space and illustrate the use of explicit and implicit time-stepping (which is performed at second-order accuracy as discussed in Section~\ref{sec:discretizingTheWaveEquation}). 
Although the accuracy of the computed results are important, our primary measure of the correctness of the WaveHoltz solution is a comparison
to the direct solution of the discretized Helmholtz equations (DHS). The DHS solutions are computed using Overture routines that in turn use direct or 
iterative sparse solvers such as those from PETSc~\cite{PETSc}. 
With corrections for the time-discretizations, as discussed in Section~\ref{sec:timeCorrections}, the WaveHoltz solution converges 
to the DHS solution to near machine precision (although this \textsl{exact} match is no longer true when using deflation).
For the examples using deflation, the numerical eigenmodes are computed using the SLEPSc package~\cite{SLEPc2005}.
We note that normally our overset grid solver for the wave equation uses upwind dissipation for stability~\cite{ssmx2023}.
However, for the WaveHoltz results presented here, no upwind dissipation is used. It appears that the WaveHoltz time filter 
is generally able to filter out any weakly unstable modes for typical use cases. 
WaveHoltz still works when upwind dissipation is included, although to achieve a near exact match with the DHS, an additional
correction to remove the effects of the dissipation is needed. 

The computations presented in subsequent sections all take $c=1$ and use a time harmonic Gaussian source term having the form 
\ba
   f(\xv,t) = \alphag \cos(\omega t) \exp\bigl( -\betag  \| \xv-\xv_0\|^2 \bigr), 
   \label{eq:gaussianSource}
\ea
where $\alphag$ is the amplitude, $\xv_0=(x_0,y_0,z_0)$ denotes the center of the Gaussian, and the exponent coefficient $\betag$ 
is determines the approximate width of the Gaussian. Note that the forcing may be adjusted for time discretization errors, and the value for $\alpha$ in the discrete WaveHoltz time filter is taken to be $\alpha_d$ given in~\eqref{eq:alphad} for all calculations, see Section~\ref{sec:fullyDiscreteConvergence}.  The value for $\alpha_d$ is used for all plots of the adjusted $\beta$ function in the subsequent subsections, and so the dependence on $\alpha$ is suppressed for notational convenience.
In addition to the convergence rates (CRs) of the iterations, we also report the effective convergence rate (ECR),
which is an adjusted CR that removes the effect of~$\Np$ (i.e.~since the cost of each wave-solve is proportional to $\Np$), 
\ba
  &  {\rm ECR} \eqdef {\rm CR}^{1/\Np}, 
\label{eq:ECRdef}
\ea
where $\Np$ is the number of periods over which the wave equation is integrated in time, $\Tbar = \Np T$.
The ECR is a better measure of run-time performance than the CR.



\subsection{Disk}  \label{sec:disk}

Helmholtz solutions are computed for a circular disk domain to demonstrate the
use of the WaveHoltz~scheme with an overset grid in two dimensions. 
The overset grid for the disk of radius $R=1$,
consists of an annular boundary-fitted grid and a background Cartesian grid, see Figure~\ref{fig:diskGridAndSolutionFig}. 
Let $\Gcd^{(j)}$ denote the disk grid with target grid spacing $\ds^{(j)}=1/(10 j)$.
The forcing is the Gaussian source in~\eqref{eq:gaussianSource} with $\alphag=-50$, $(x_0,y_0)=(0.25,0.25)$, and $\betag=10$.
The frequency is taken as $\omega=8.1$. Homogeneous Dirichlet boundary conditions are imposed.

{
\newcommand{\drawContour}[7]{%
\begin{scope}[#1]
\draw(0.0,0) node[anchor=south west,xshift=-4pt,yshift=+0pt] {\trimfiga{#2}{\figWidtha}};
  \draw(.5,.5) node[draw,fill=white,anchor=west,xshift=2pt,yshift=1pt] {\scriptsize #3};
\begin{scope}[xshift=-.2cm,yshift=+0pt]
  \draw (\xcb,\ycb) node[anchor=south west,xshift=0.25cm,yshift=.5cm,rotate=-90] {\trimfigcb{colourBarLines}{\cbWidth}{\cbHeight}};
  \draw (.8,0) node[anchor=north,xshift=+3pt,yshift=+2pt] {\scriptsize $#6$};
  \draw (4.8,0) node[anchor=north,xshift=+0pt,yshift=+2pt] {\scriptsize $#7$};
\end{scope}
\end{scope}
}
\newcommand{\cbWidth}{.2cm}
\newcommand{\cbHeight}{4cm}
\newcommand{\xcb}{.5cm}
\newcommand{\ycb}{-.2cm}
\setlength{\ycbTop}{\ycb+\cbHeight}
\setlength{\ycbMid}{\ycb+\cbHeight*\real{.5}}
\newcommand{\trimfigcb}[3]{\includegraphics[width=#2, height=#3, clip, trim=17cm 2.35cm 1.65cm 2.35cm]{#1}}
\newcommand{\figWidtha}{4.5cm}
\newcommand{\trimfiga}[2]{\trimw{#1}{#2}{.11}{.115}{.11}{.11}}
\begin{figure}[htb]
\begin{center}
\begin{tikzpicture}
   \useasboundingbox (0,.25) rectangle (10,4.5);  

   \begin{scope}[yshift=0cm]
     \figByWidth{   0}{-.2}{sicGridG2}{5cm}[0.1][0.1][0.1][0.1]
     \drawContour{xshift=5.cm,yshift=0.00cm}{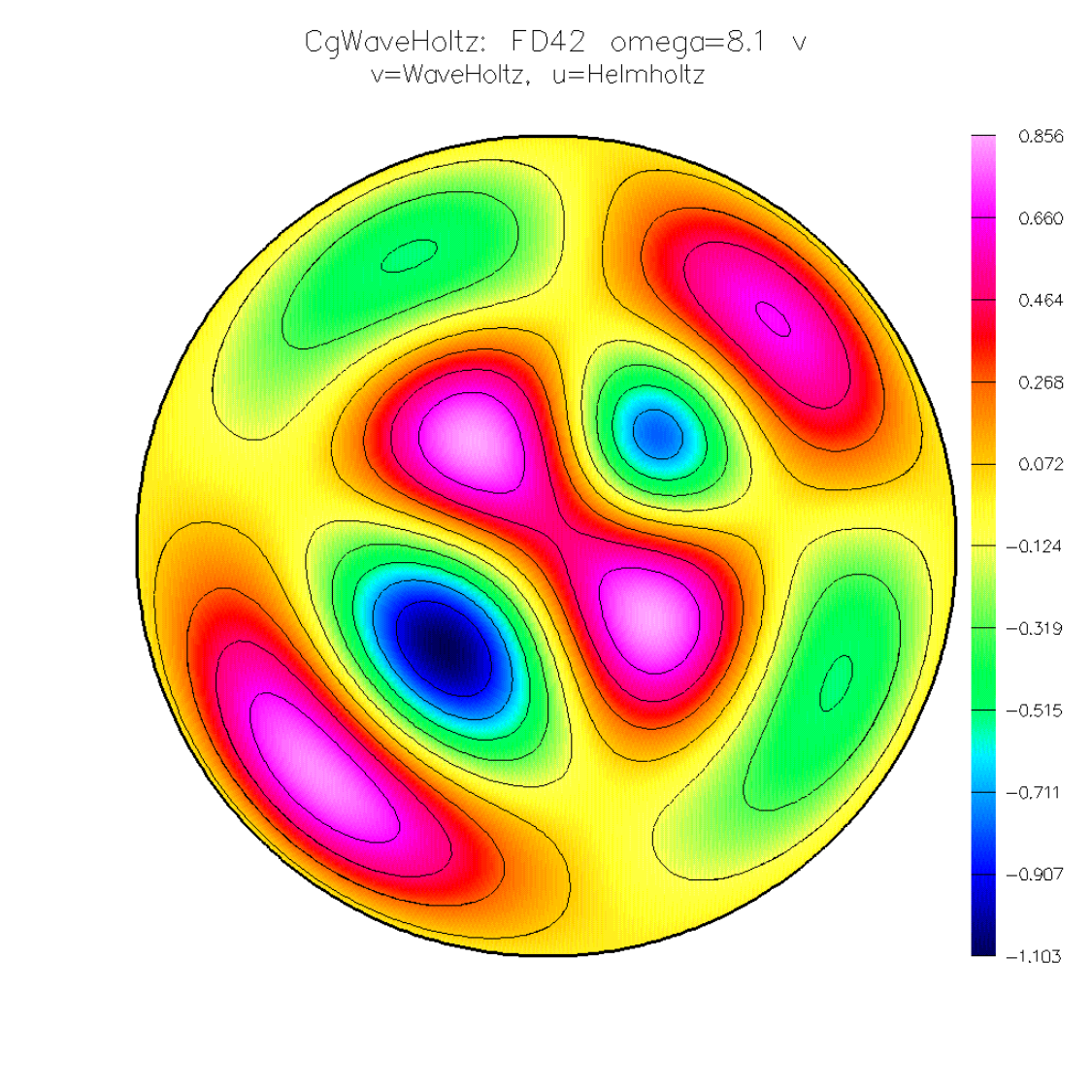}{$v$}{$v$}{$t=0.3$}{$-1.1$}{$.86$};     
   \end{scope}

\end{tikzpicture}
\end{center}
\caption{Left: coarse overset grid $\Gcd^{(2)}$ for the disk. Right: computed WaveHoltz solution, order of accuracy four, on grid $\Gcd^{(8)}$. 
    }
\label{fig:diskGridAndSolutionFig}
\end{figure}
}

{
\newcommand{\figSize}{6.5cm}

\begin{figure}[htb]
\begin{center}
\begin{tikzpicture}[scale=1]
  \useasboundingbox (0,.7) rectangle (14,10.75);  
   \begin{scope}[yshift=5.5cm]
    \figByWidth{  0}{0}{diskG4O2Freq8p1Np1Nits10Exp}{\figSize}[0][0][0][0]
    \figByWidth{7.0}{0}{diskG4O2Freq8p1Np1Nits10Imp}{\figSize}[0][0][0][0]
    \draw (2.8,3) node[draw,fill=white,anchor=west,xshift=0pt,yshift=0pt,inner sep=3pt] {\scriptsize Disk G4, O2, E, $\Np=1$};
    \draw (10,3 ) node[draw,fill=white,anchor=west,xshift=0pt,yshift=0pt,inner sep=3pt] {\scriptsize Disk G4, O2, I, $\Np=1$};
  \end{scope}
   \begin{scope}[yshift=0cm]
    \figByWidth{0.0}{0}{diskG4O4Freq8p1Np1Nits10Exp}{\figSize}[0][0][0][0]
    \figByWidth{7.0}{0}{diskG4O4Freq8p1Np1Nits10Imp}{\figSize}[0][0][0][0]
    \draw (2.8,3) node[draw,fill=white,anchor=west,xshift=0pt,yshift=0pt,inner sep=3pt] {\scriptsize Disk G4, O4, E, $\Np=1$};
    \draw (10,3 ) node[draw,fill=white,anchor=west,xshift=0pt,yshift=0pt,inner sep=3pt] {\scriptsize Disk G4, O4, I, $\Np=1$};
  \end{scope}  

\end{tikzpicture}
\end{center}
\caption{WaveHoltz: disk. Convergence of the FPI and GMRES accelerated WaveHoltz iterations. 
    Top row: order of accuracy four. 
    Bottom row: order of accuracy two.
    Left column: explicit time-stepping. Right column: implicit time-stepping with $10$ time-steps per period.
    } 
\label{fig:diskExplicitImplicit}
\end{figure}
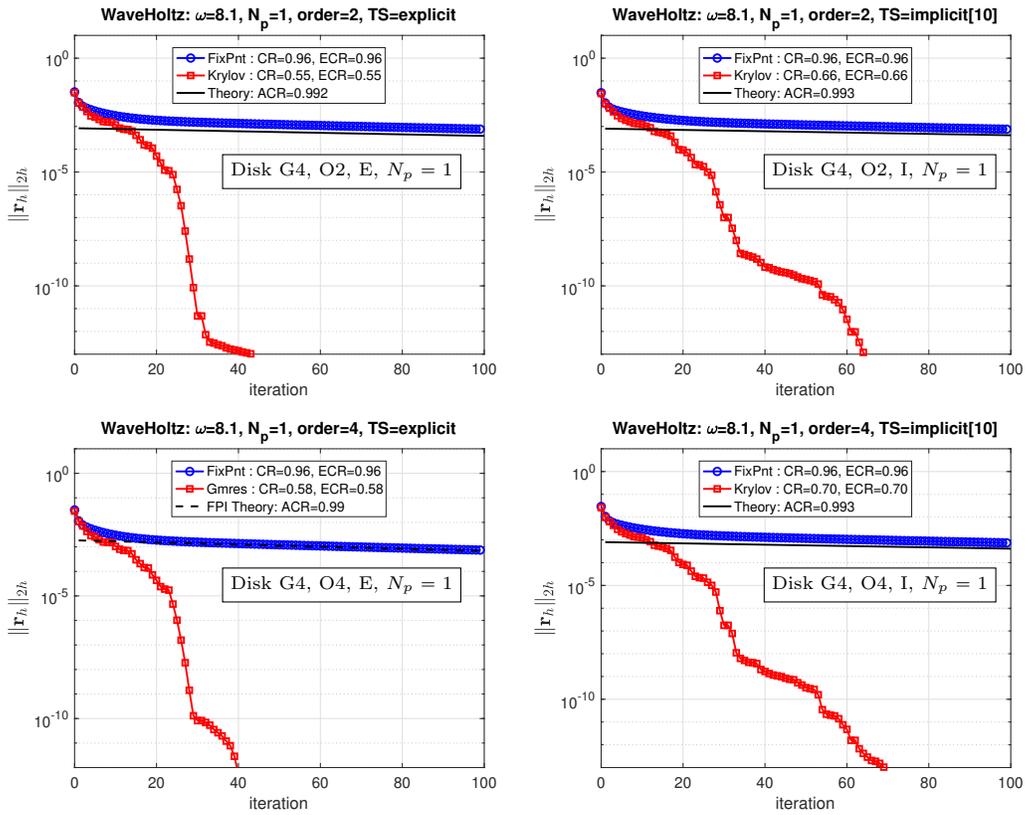
}

{
\newcommand{\figSize}{6.5cm}

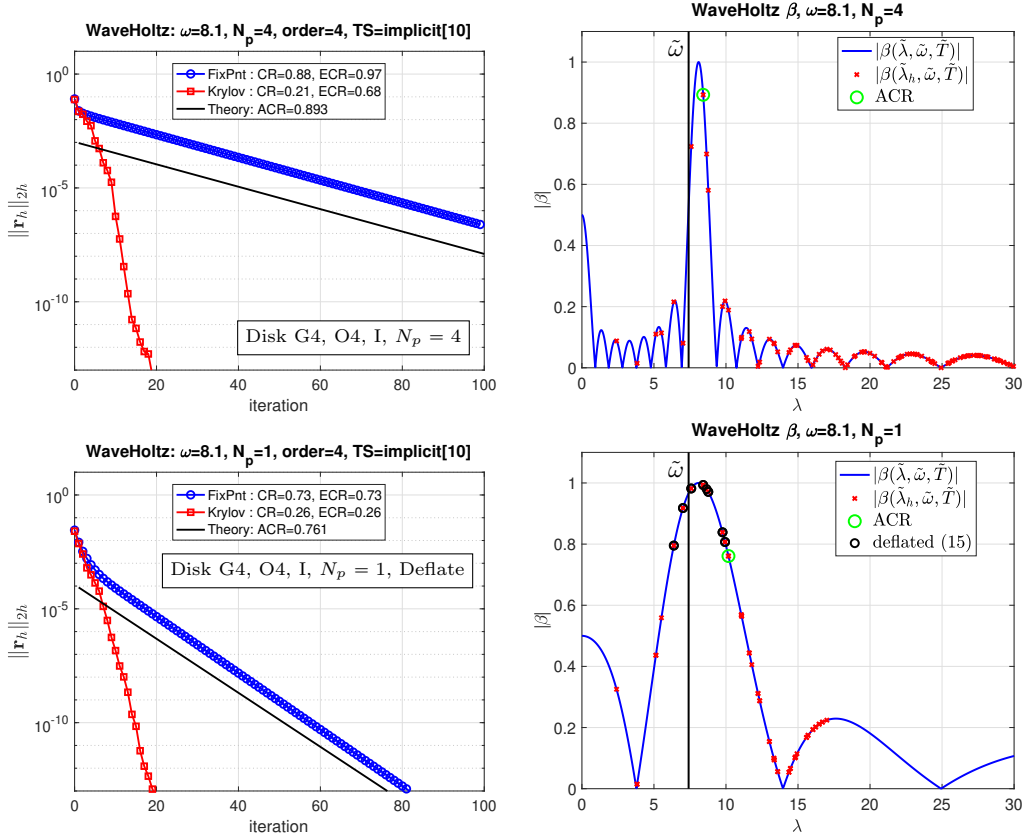
\begin{figure}[htb]
\begin{center}
\begin{tikzpicture}[scale=1]
  \useasboundingbox (0,.7) rectangle (14,10.75);  

   \begin{scope}[yshift=5.6cm]
    \figByWidth{0.0}{0}{diskG4O4Freq8p1Np4Nits10Imp}{\figSize}[0][0][0][0]
    \figByWidth{7.0}{0}{diskG4O4Freq8p1Np4Nits10ImpFixPointMuFunction}{\figSize}[0][0][0][0]
    \draw (3, 1) node[draw,fill=white,anchor=west,xshift=0pt,yshift=0pt,inner sep=3pt,inner sep=3pt] {\scriptsize Disk G4, O4, I, $\Np=4$};
    \draw (9, 5) node[anchor=east,xshift=2pt,yshift=-4pt] {\small $\omegaTilde$};
  \end{scope} 
   \begin{scope}[yshift=0cm]
    \figByWidth{0.0}{0}{diskG4O4Freq8p1Np1Nits10Deflate8Imp}{\figSize}[0][0][0][0]
    \figByWidth{7.0}{0}{diskG4O4Freq8p1Np1Nits10Deflate8ImpFixPointMuFunction}{\figSize}[0][0][0][0]
    \draw (2, 3.5) node[draw,fill=white,anchor=west,xshift=0pt,yshift=0pt,inner sep=3pt,inner sep=3pt] {\scriptsize Disk G4, O4, I, $\Np=1$, Deflate};
    \draw (9, 5) node[anchor=east,xshift=2pt,yshift=-4pt] {\small $\omegaTilde$};
  \end{scope}     

\end{tikzpicture}
\end{center}
\caption{WaveHoltz: disk. 
  Top left: implicit time-stepping with $\Np=4$ periods per time-step.
  Top right: Magnitude of the WaveHoltz filter function $|\beta(\lambda)|$ with eigenvalues marked with red x's.
  Bottom left: implicit time-stepping with deflation. 
  Bottom right: plot of $|\beta(\lambdaTilde,\omegaTilde,\Ttilde)|$  versus $\lambda$, with eigenvalues marked with red x's and deflated eigenvalues marked with black circles.
  The black vertical lines on the right graphs indicate the values of the adjusted frequencies $\omegaTilde$ used to correct for time discretization errors. 
    } 
\label{fig:diskImplicitDeflate}
\end{figure}
}

Figure~\ref{fig:diskExplicitImplicit} shows results for implicit and explicit time-stepping and for second- and fourth-order
accurate approximations using $\Np=1$ periods per time interval. 
The graphs show the scaled $L_2$-norm of the discrete residual vector $\rv_h^{(k)}$,  
\ba 
  \| \rv_h^{(k)} \|_{2h} \eqdef \f{1}{\sqrt{N}} \, \| \rv_h^{(k)} \|_2, 
\ea
versus iteration $k$, where $N$ is the total number of grid points.
For the fixed-point iteration, $\rv_h^{(k)} = \vv_h^{(k)}-\vv_h^{(k-1)}$, while for the GMRES algorithm 
$\rv_h^{(k)} = \bv_h - M \vv_h^{(k)}$ (see equation~\eqref{eq:WaveHoltzLinearSystem}).
The implicit time-stepping results in Figure~\ref{fig:diskExplicitImplicit} used $10$ time-steps per period which corresponds to a grid CFL number of about $60$ with respect to the smallest grid cell. 
Convergence rates for the FPI and GMRES accelerated schemes are shown along with the theoretical estimate for the asymptotic convergence rate (ACR).
The theoretical ACR is determined using the true eigenvalues of the disk. 
The FPI convergence rate is seen to agree well with the theory. 
The Krylov convergence rates using GMRES show good accelerations compared to the FPI.
The convergence rates are roughly the same for the second- and fourth-order accurate approximations.
In all cases the converged WaveHoltz solution agrees with the direct solution of the discrete Helmholtz BVP to a relative difference of about $10^{-12}$.

The top graphs in Figure~\ref{fig:diskImplicitDeflate} shows results using $\Np=4$ periods per time interval with implicit time-stepping and
fourth-order accuracy. The top-right graph shows the absolute value of the WaveHoltz filter function
 $\beta(\lambdaTilde\;\,\omegaTilde,\Ttilde,\alpha_d)$ versus $\lambda$, together with the locations
of the eignvalues (red~x's). A green circle marks the value of $\beta$ that determines the asymptotic convergence rate (ACR). This value of about $0.89$
approximately matches the FPI convergence shown in the top-left graph.
The black vertical lines on the right graphs indicate the values of the adjusted frequencies~$\omegaTilde$
used in the WaveHoltz forcing to correct for time discretization errors (see Section~\ref{sec:timeCorrections}).
 The FPI theoretical ACR is seen to improve from about $0.99$ ($\Np=1$) to about $0.89$ ($\Np=4$), the computed CRs are somewhat better, 
 while the GMRES converges quite a bit faster, ${\rm CR}\approx0.21$. 
The bottom graphs in Figure~\ref{fig:diskImplicitDeflate} show the effect of deflating $15$ eigenmodes ($15$ includes counting any multiple eigenvalues). 
The convergence rates of the FPI and GMRES with deflation are much improved with the ACR for the FPI matching the theoretical value.

\newcommand{\Gcs}{\Gc_s}
\subsection{Solid sphere}  \label{sec:solidSphere}

{
\newcommand{\drawContour}[7]{%
\begin{scope}[#1]
\draw(0.0,0) node[anchor=south west,xshift=-4pt,yshift=+0pt] {\trimfiga{#2}{\figWidtha}};
  \draw(.5,.5) node[draw,fill=white,anchor=west,xshift=2pt,yshift=1pt] {\scriptsize #3};
\begin{scope}[xshift=-.2cm,yshift=+0pt]
  \draw (\xcb,\ycb) node[anchor=south west,xshift=0.25cm,yshift=.5cm,rotate=-90] {\trimfigcb{colourBarLines}{\cbWidth}{\cbHeight}};
  \draw (.8,0) node[anchor=north,xshift=+3pt,yshift=+2pt] {\scriptsize $#6$};
  \draw (4.8,0) node[anchor=north,xshift=+0pt,yshift=+2pt] {\scriptsize $#7$};
\end{scope}
\end{scope}
}
\newcommand{\cbWidth}{.2cm}
\newcommand{\cbHeight}{4cm}
\newcommand{\xcb}{.5cm}
\newcommand{\ycb}{-.2cm}
\setlength{\ycbTop}{\ycb+\cbHeight}
\setlength{\ycbMid}{\ycb+\cbHeight*\real{.5}}
\newcommand{\trimfigcb}[3]{\includegraphics[width=#2, height=#3, clip, trim=17cm 2.35cm 1.65cm 2.35cm]{#1}}
\newcommand{\figWidtha}{4.5cm}
\newcommand{\trimfiga}[2]{\trimw{#1}{#2}{.11}{.115}{.11}{.11}}
\begin{figure}[htb]
\begin{center}
\begin{tikzpicture}
   \useasboundingbox (0,.45) rectangle (11,4.75);  

   \begin{scope}[yshift=0cm]
     \figByWidth{   0}{-.2}{sphereGridsExploded}{5.25cm}[0.][0.][0.][0.]
     \drawContour{xshift=6.cm,yshift=0.00cm}{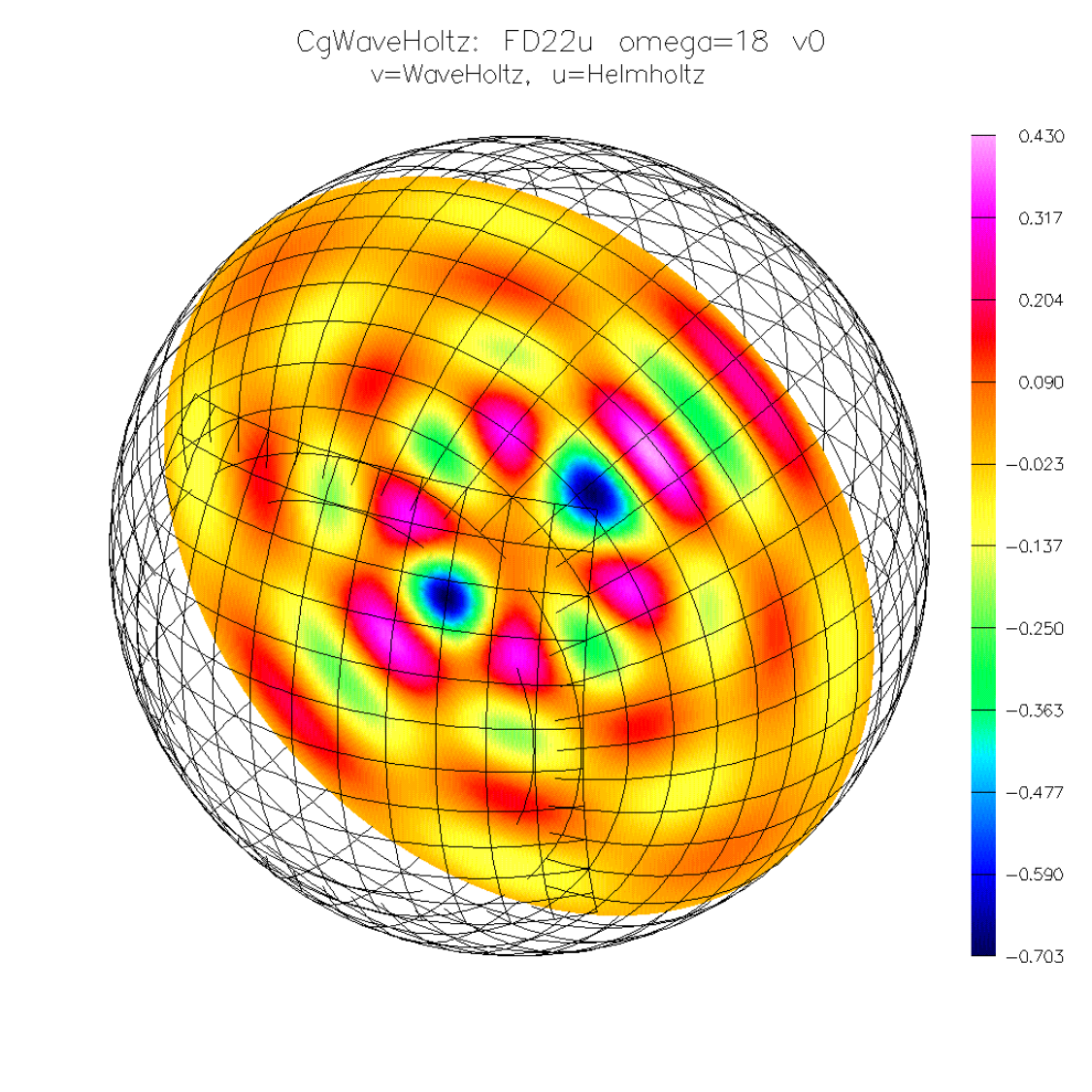}{$v^{(1)}$}{$v$}{$t=0.3$}{$-.70$}{$.43$};     
   \end{scope}

\end{tikzpicture}
\end{center}
\caption{Gaussian source in a solid. Left: exploded view of the overlapping surface patches on the sphere.
Right: Contours of the solution on a cutting plane, with a coarsened version of the grid.
The solution was computed on grid $\Gcs^{(4)}$ to second order of accuracy, and frequency $\omega=18$.
    }
\label{fig:sphereOneFreqSolution}
\end{figure}
}

In this section the WaveHoltz scheme is used to solve Helmholtz problems for a solid spherical domain of radius one.
The overset grid for the domain, denoted by $\Gcs^{(j)}$ with target grid spacing
$\ds^{(j)}=1/(10 j)$, consists of four component grids as shown in Figure~\ref{fig:sphereOneFreqSolution}.  There is
a background Cartesian grid covering the central portion of the solid (not visible in the figure) together with three surface-fitted grids to represent the sphere surface.
The problem is forced by a Gaussian source with frequency $\omega=8.5$, 
strength $\alphag=150$ and exponent $\betag=10$.  The source is centered at $\xv_0=(0.1,0.1,0.1)$.
The boundary conditions are taken to be of Dirichlet type, similar to the previous disk problem.

{
\newcommand{\figSize}{6.5cm}

\begin{figure}[htb]
\begin{center}
\begin{tikzpicture}[scale=1]
  \useasboundingbox (0,.7) rectangle (14,11);  

   \begin{scope}[yshift=5.7cm]
    \figByWidth{0.0}{0}{sphereOneFreqNp4G2Order4Exp}{\figSize}[0][0][0][0]
    \figByWidth{7.0}{0}{sphereOneFreqNp4G2Order4ExpFixPointMuFunction}{\figSize}[0][0][0][0]
    \draw (3.5,1) node[draw,fill=white,anchor=west,xshift=0pt,yshift=0pt,inner sep=3pt] {\scriptsize Sphere, G2, E, O4};
  \end{scope}  
   \begin{scope}[yshift=0cm]
     \figByWidth{0.0}{0}{sphereOneFreqNp4G2Order4Imp}{\figSize}[0][0][0][0]
     \figByWidth{7.0}{0}{sphereOneFreqNp4G2Order4ImpFixPointMuFunction}{\figSize}[0][0][0][0]
     \draw (3.5,1) node[draw,fill=white,anchor=west,xshift=0pt,yshift=0pt,inner sep=3pt] {\scriptsize Sphere, G2, I, O4};
  \end{scope}

\end{tikzpicture}
\end{center}
\caption{Sphere: grid $\Gcs^{(2)}$, order of accuracy 4. Top row: explicit time-stepping.
Bottom row: implicit time-stepping. 
The black vertical lines on the right graphs indicate the values of the adjusted frequencies $\omegaTilde$ used to correct for time discretization errors. 
} 
\label{fig:sphereOneFreq}
\end{figure}
}

Figure~\ref{fig:sphereOneFreq} shows convergence results of the iterations for the sphere
using explicit and implicit time-stepping with $\Np=4$ periods per time interval.
The implicit time-stepping used $\NITS=10$ time-steps per period for a total of $40$ time-steps per wave-solve.
The computed CRs for the fourth-order accurate scheme  (explicit and implicit time-stepping) 
on grid $\Gcs^{(2)}$ are in good agreement with the 
theory (the theoretical ACR is computed using the exact eigenvalues). 
The GMRES accelerated ECRs are very good.
We note that the computation of the direct Helmholtz solution (DHS) for this problem requires too much memory to use a direct sparse solver, and so
an iterative solver is used instead. GMRES with an ILU(100) preconditioner is used.  The large fill-in level of $100$ for ILU is needed
to avoid a failure of the algorithm.


\newcommand{\Gcp}{\Gc_p}
\subsection{Solid pipe}  \label{sec:pipe}

The WaveHoltz scheme is now used to solve Helmholtz problems for a pipe geometry.
The pipe is a solid cylinder of radius $R_b=0.5$ and axial range of $z\in[0,L_z]$ with $L_z=1$.
The overset grid for the domain, denoted by $\Gcp^{(j)}$ with target grid spacing
$\ds^{(j)}=1/(10 j)$, consists of two component grids, 
a background Cartesian grid together with a cylindrical shell near the cylinder surface, as shown in Figure~\ref{fig:pipeGridAndSolution}.
The boundary conditions are taken as periodic in the axial ($z$) direction and Dirichlet on the curved outer boundary of the pipe.
The problem is forced with Gaussian source with frequency $\omega=8.5$, 
strength $\alphag=150$ and exponent $\betag=10$, and it is centered at $\xv_0=(0.1,0.1,0.5)$.

{
\newcommand{\drawContour}[7]{%
\begin{scope}[#1]
\draw(0.0,0) node[anchor=south west,xshift=-4pt,yshift=+0pt] {\trimfiga{#2}{\figWidtha}};
  \draw(.5,.5) node[draw,fill=white,anchor=west,xshift=2pt,yshift=1pt,inner sep=3pt] {\scriptsize #3};
\begin{scope}[xshift=-.2cm,yshift=-2pt]
  \draw (\xcb,\ycb) node[anchor=south west,xshift=0.25cm,yshift=.5cm,rotate=-90] {\trimfigcb{colourBarLines}{\cbWidth}{\cbHeight}};
  \draw (.8,0) node[anchor=north,xshift=+3pt,yshift=+2pt] {\scriptsize $#6$};
  \draw (4.8,0) node[anchor=north,xshift=+0pt,yshift=+2pt] {\scriptsize $#7$};
\end{scope}
\end{scope}
}
\newcommand{\cbWidth}{.2cm}
\newcommand{\cbHeight}{4cm}
\newcommand{\xcb}{.5cm}
\newcommand{\ycb}{-.2cm}
\setlength{\ycbTop}{\ycb+\cbHeight}
\setlength{\ycbMid}{\ycb+\cbHeight*\real{.5}}
\newcommand{\trimfigcb}[3]{\includegraphics[width=#2, height=#3, clip, trim=17cm 2.35cm 1.65cm 2.35cm]{#1}}
\newcommand{\figWidtha}{5cm}
\newcommand{\trimfiga}[2]{\trimw{#1}{#2}{.07}{.115}{.07}{.11}}
\begin{figure}[htb]
\begin{center}
\begin{tikzpicture}
   \useasboundingbox (0,.3) rectangle (11,5);  

   \begin{scope}[yshift=0cm]
     \figByWidth{0.0}{0.15}{pipeGridG2}{5.6cm}[0.][0.][0.05][0.05]
     \drawContour{xshift=6.cm,yshift=0.00cm}{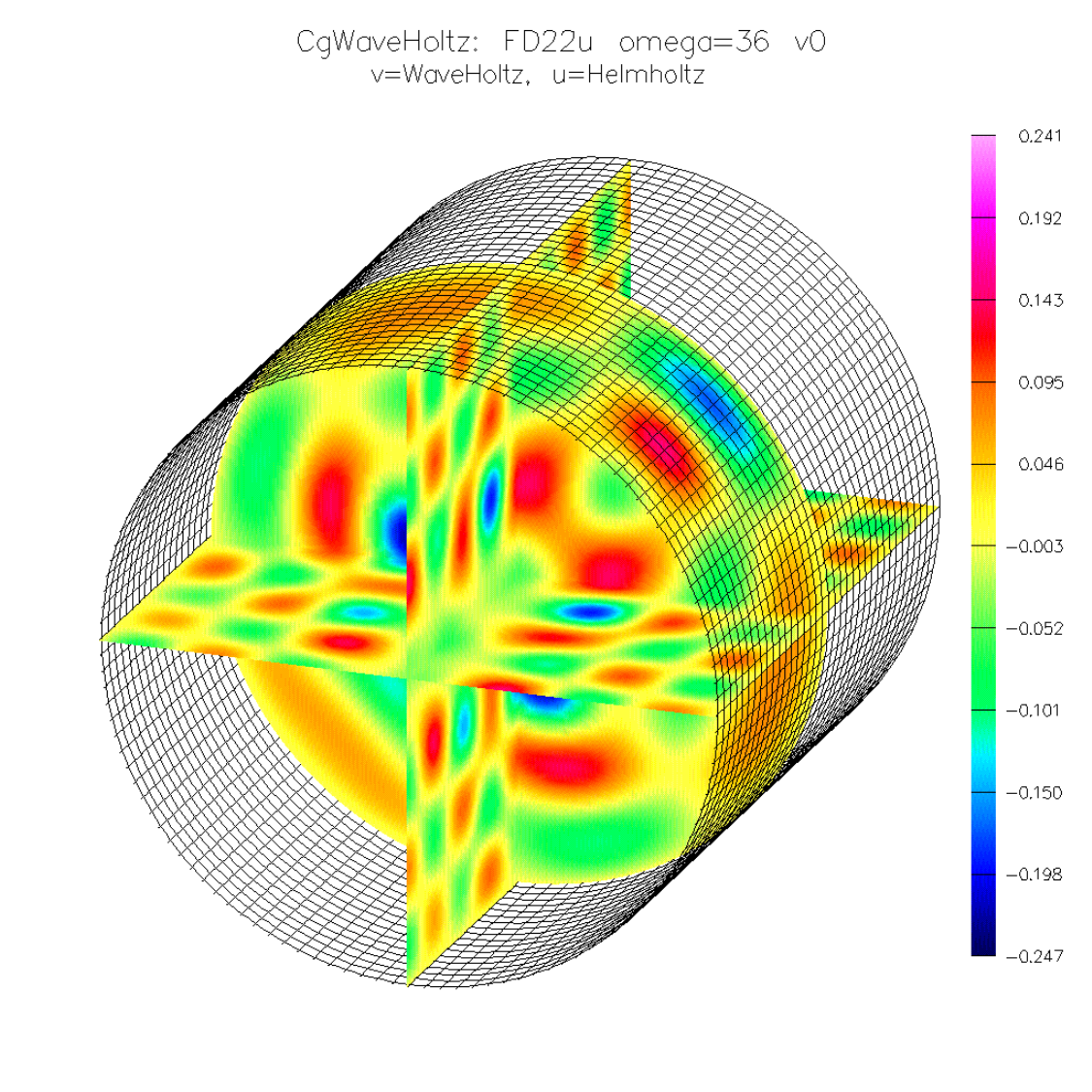}{$\omega=36$}{$v$}{$t=0.3$}{$-.25$}{$.24$}     
   \end{scope}

\end{tikzpicture}
\end{center}
\caption{Gaussian source in a pipe.  
Left: overset grid $\Gcp^{(2)}$ for a solid cylinder.
Right: Contours of the solution on cutting planes, with a coarsened version of the grid.
The solution was computed on grid $\Gcp^{(8)}$ to second order of accuracy, with frequency $\omega=36$.
    }
\label{fig:pipeGridAndSolution}
\end{figure}
}

{
\newcommand{\figSize}{6.5cm}

\begin{figure}[htb]
\begin{center}
\begin{tikzpicture}[scale=1]
  \useasboundingbox (0,.7) rectangle (14,11);  

  \begin{scope}[yshift=5.7cm]
    \figByWidth{  0}{0}{pipeFreq8p5Np4G2O4Exp}{\figSize}[0][0][0][0]
    \figByWidth{7.0}{0}{pipeFreq8p5Np4G2O4ExpFixPointMuFunction}{\figSize}[0][0][0][0]

    \draw (2.5,1.5) node[draw,fill=white,anchor=west,xshift=0pt,yshift=0pt,inner sep=3pt] {\scriptsize Pipe2, E, O4};
  \end{scope}  
   \begin{scope}[yshift=0cm]
    \figByWidth{0.0}{0}{pipeFreq8p5Np4G2O4Imp}{\figSize}[0][0][0][0]
    \figByWidth{7.0}{0}{pipeFreq8p5Np4G2O4ImpFixPointMuFunction}{\figSize}[0][0][0][0]

    \draw (2.5,1.5) node[draw,fill=white,anchor=west,xshift=0pt,yshift=0pt,inner sep=3pt] {\scriptsize Pipe2, I, O4};
  \end{scope}  

\end{tikzpicture}
\end{center}
\caption{Pipe: Grid $\Gcp^{(2)}$, order of accuracy four.
   Top row: explicit time-stepping.
Bottom row: implicit time-stepping.
The black vertical lines on the right graphs indicate the value of the adjusted frequency $\omegaTilde$ used to correct for time discretization errors. 
    } 
\label{fig:pipeOneFreq}
\end{figure}
}

Figure~\ref{fig:pipeOneFreq} shows the iteration convergence results for the pipe when solving 
with $\Np=4$.
The implicit time-stepping uses $\NITS=10$ time-steps per period for a total of $40$ time-steps per wave-solve.
The theoretical convergence rate of the FPI is estimated using the exact eigenvalues. 
The computed CRs for the fourth-order accurate schemes on grid $\Gcs^{(2)}$ are in good agreement with the 
theory.

\newcommand{\Gcde}{\Gc_{\rm de}}
\subsection{Double ellipse (unilluminable room)} \label{sec:doubleEllipse}

As a next example, we solve the Helmholtz problem for the Penrose unilluminable room~\cite{Fukushima2015LightPI}.
The geometry, shown in Figure~\ref{fig:doubleEllipseGridFig}, is designed so that the 
some of the alcoves, two at the top and two at the bottom of the domain, remain dark (or quiet) when there is a light source (or sound source) in the interior.
The design is based on two ellipses of different sizes. 
Two smaller half-ellipses, with semi-axes $(a_1,b_1)=(2,1)$, are located at the top and bottom. 
Two larger half-ellipses, with semi-axes $(a_2,b_2)=(3,6)$, are placed on the left and right.
The left and right ends of the smaller ellipses are located at the foci of the larger ellipses.

{
\newcommand{\drawContour}[7]{%
\begin{scope}[#1]
\draw(0.0,0) node[anchor=south west,xshift=-4pt,yshift=+0pt] {\trimfiga{#2}{\figWidtha}};
\begin{scope}[xshift=.0cm,yshift=0pt]
  \draw (\xcb,\ycb) node[anchor=south west,xshift=0.25cm,yshift=.5cm,rotate=-90] {\trimfigcb{colourBarLines}{\cbWidth}{\cbHeight}};
  \draw (.8,0) node[anchor=north,xshift=+3pt,yshift=+2pt] {\scriptsize $#6$};
  \draw (4.8,0) node[anchor=north,xshift=+0pt,yshift=+2pt] {\scriptsize $#7$};
\end{scope}
\end{scope}
}
\newcommand{\cbWidth}{.2cm}
\newcommand{\cbHeight}{4cm}
\newcommand{\xcb}{.5cm}
\newcommand{\ycb}{-.2cm}
\setlength{\ycbTop}{\ycb+\cbHeight}
\setlength{\ycbMid}{\ycb+\cbHeight*\real{.5}}
\newcommand{\trimfigcb}[3]{\includegraphics[width=#2, height=#3, clip, trim=17cm 2.35cm 1.65cm 2.35cm]{#1}}
\newcommand{\figWidtha}{5.85cm}
\newcommand{\trimfiga}[2]{\trimh{#1}{#2}{.17}{.17}{.1}{.1}}
\newcommand{\figSize}{5.5cm}
\newcommand{\figh}{5.5cm}
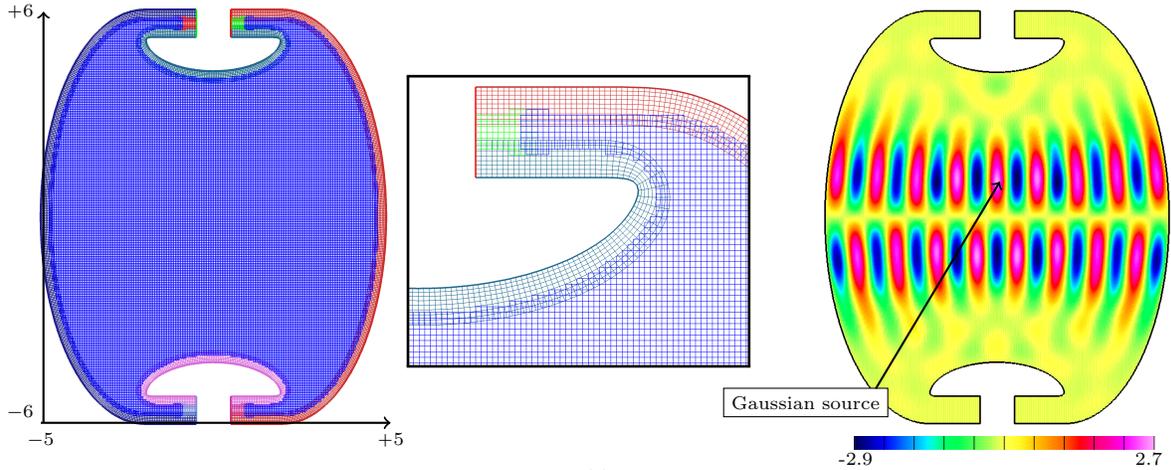
\begin{figure}[htb]
\begin{center}
\begin{tikzpicture}
  \useasboundingbox (0,.1) rectangle (15,5.8);  

  \figByWidth{0.0}{-6pt}{darkCornerRoomGridG2}{4.7cm}[0.1][0.1][0.][0.]
  \draw[thick,black,->,yshift=0pt] (0,0) node[anchor=north] {\scriptsize $-5$} -- (4.65,0.00) node[anchor=north] {\scriptsize $+5$};  
  \draw[thick,black,->,xshift=1pt] (0,0) node[anchor=east,yshift=4pt]  {\scriptsize $-6$} -- (0.00,5.47) node[anchor=east ] {\scriptsize $+6$};  

  \figByWidthb{4.95}{.75}{darkCornerRoomGridG2Zoom}{4.5cm}[0.][0.][0.][0.15]

  \begin{scope}[xshift=10cm,yshift=0]

     \drawContour{xshift=-2pt,yshift=-10pt}{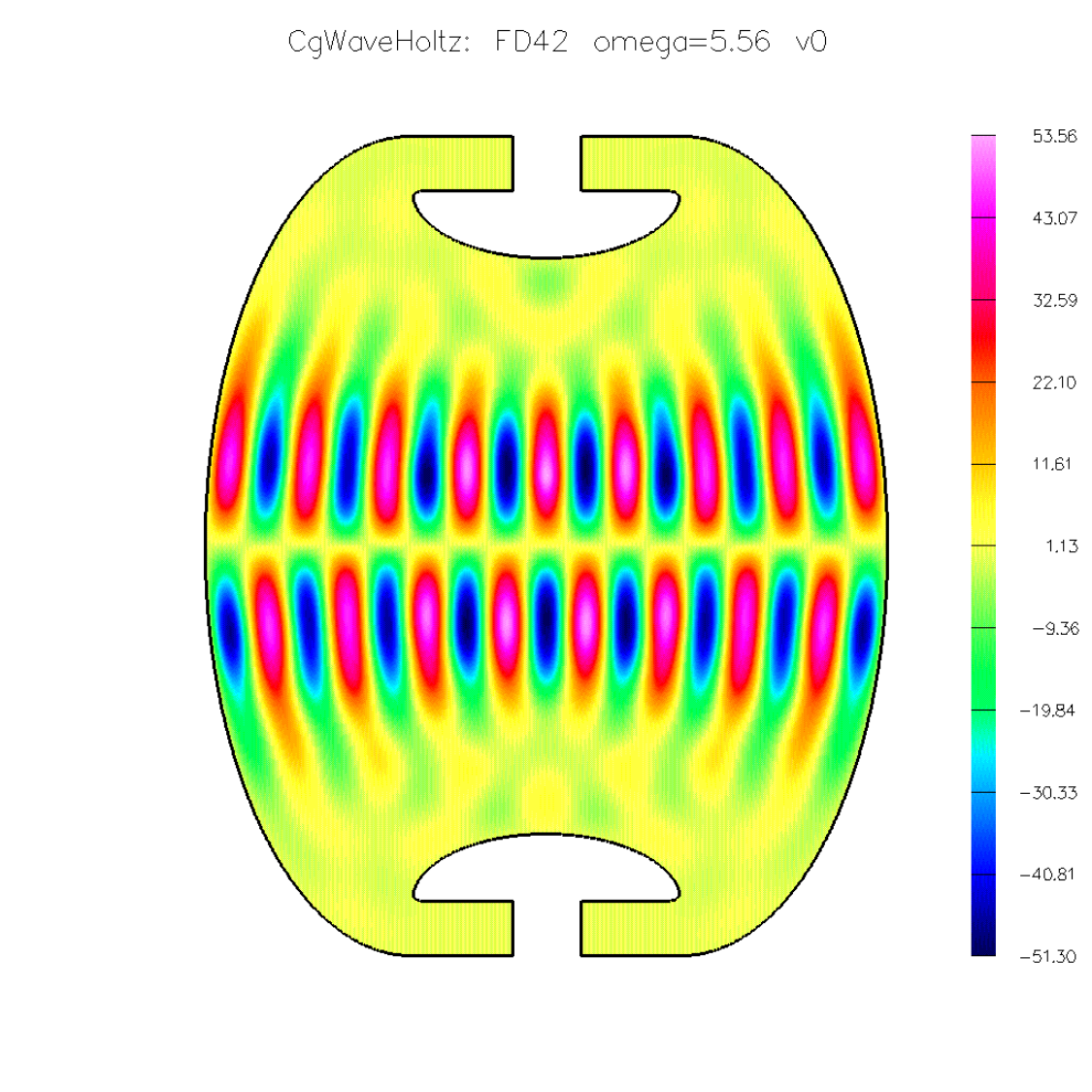}{$v$}{$v$}{$t=0.3$}{$-2.9$}{$2.7$}
     \draw(-1,.2) node[draw,fill=white,anchor=west,xshift=2pt,yshift=2pt,inner sep=3pt] {\scriptsize Gaussian source};  
     \draw[thick,black,->] (1.1,.455) -- (2.75,3.2083);  

  \end{scope} 
\end{tikzpicture}
\end{center}
\caption{Left: double ellipse geometry and overset grid $\Gcde^{(2)}$. Middle: closeup of a portion of the grid.
Right: computed Helmholtz solution for $\omega=5.56$ with a Gaussian source at $(0,1)$.
 }
\label{fig:doubleEllipseGridFig}
\end{figure}
}

The overset grid for the domain is shown in Figure~\ref{fig:doubleEllipseGridFig} (left and middle).
The grid, denoted by $\Gcde^{(j)}$ with target grid spacing $\ds^{(j)}=1/(10 j)$, consists of a total of nine component grids.  Four component grids are placed to fit the curved elliptical boundaries with four small Cartesian grids used to fit the straight portions of the boundaries in the alcoves (see middle image).  The ninth component grid is a large background Cartesian grid covering the bulk of the domain.
Figure~\ref{fig:doubleEllipseGridFig} (right) shows a sample solution computed by the WaveHoltz algorithm, with  homogeneous Dirichlet boundary conditions, for a
Gaussian source~\eqref{eq:gaussianSource} with $\omega=5.56$, $\alphag=400$, $\betag=10$, and $\xv_0=(0,1)$.
The forcing excites a harmonic mode that is active primarily near the center of the domain.

The subsequent calculations for this problem are used, in part, to assess the rule-of-thumb estimates for the PPW given in Recipe~\ref{recipe:PPW}.
To this end, Table~\ref{tab:doubleEllipsePPW} lists PPW data for a range of values for the frequency $\omega$ and the grid resolution given by the index~$j$.
The rule-of-thumb estimates given by $\PPW_p$ for second-order ($p=2$) and fourth-order ($p=4$) accurate approximations use a domain size of $L=12$, i.e.~the height of the domain shown in Figure~\ref{fig:doubleEllipseGridFig}, and a relative error of $\epsilon=10\sp{-2}$.  In the table, the wavelength is computed as $\Lambda=2\pi c/\omega$ (with $c=1$) while the values in the column titled \textit{Actual $\PPW$} are $\Lambda/\ds$.

\begin{table}[htb]\tableFont 
\begin{center}
\begin{tabular}{|c|c|c|c|c|c|c|c|} \hline 
  \multicolumn{8}{|c|}{Double Ellipse, Points-per-wavelength, $\eps=10^{-2}$, $L=12$}    \\ 
   \multicolumn{5}{|c|}{} & Actual & \multicolumn{2}{c|}{Estimated} \\  \hline
   $\omega$  & $\Lambda$  & $N_\Lambda$  &  $j$     & $\ds$  &  $\PPW$  & $\PPW_2$ & $\PPW_4$   \\ \hline 
  $10$   &  $0.628$  & $ 19.1$ &  $ 4$  &  2.50e-02 & $ 25.1$  &  $140.5$  & $ 18.0$  \\
  $10$   &  $0.628$  & $ 19.1$ &  $ 8$  &  1.25e-02 & $ 50.3$  &  $140.5$  & $ 18.0$  \\
  $10$   &  $0.628$  & $ 19.1$ &  $16$  &  6.25e-03 & $100.5$  &  $140.5$  & $ 18.0$  \\
  $10$   &  $0.628$  & $ 19.1$ &  $32$  &  3.13e-03 & $201.1$  &  $140.5$  & $ 18.0$  \\
\hline
  $15$   &  $0.419$  & $ 28.6$ &  $ 4$  &  2.50e-02 & $ 16.8$  &  $172.1$  & $ 19.9$  \\
  $15$   &  $0.419$  & $ 28.6$ &  $ 8$  &  1.25e-02 & $ 33.5$  &  $172.1$  & $ 19.9$  \\
  $15$   &  $0.419$  & $ 28.6$ &  $16$  &  6.25e-03 & $ 67.0$  &  $172.1$  & $ 19.9$  \\
  $15$   &  $0.419$  & $ 28.6$ &  $32$  &  3.13e-03 & $134.0$  &  $172.1$  & $ 19.9$  \\
\hline
  $40$   &  $0.157$  & $ 76.4$ &  $ 4$  &  2.50e-02 & $  6.3$  &  $281.0$  & $ 25.4$  \\
  $40$   &  $0.157$  & $ 76.4$ &  $ 8$  &  1.25e-02 & $ 12.6$  &  $281.0$  & $ 25.4$  \\
  $40$   &  $0.157$  & $ 76.4$ &  $16$  &  6.25e-03 & $ 25.1$  &  $281.0$  & $ 25.4$  \\
  $40$   &  $0.157$  & $ 76.4$ &  $32$  &  3.13e-03 & $ 50.3$  &  $281.0$  & $ 25.4$  \\
  \hline
\end{tabular}
\caption{Actual and estimated points-per-wavelength for the double ellipse domain as a function of frequency $\omega$ and grid resolution $j$. 
 The column titled $\PPW$ holds the actual points-per-wavelength. The columns labeled $\PPW_2$ and $\PPW_4$ contains the rule-of-thumb estimated values
 for second- and fourth-order accurate schemes, respectively, from Recipe~\ref{recipe:PPW}.
    }
\label{tab:doubleEllipsePPW}
\end{center}
\end{table}

{
\newcommand{\drawContour}[7]{%
\begin{scope}[#1]
\draw(0.0,0) node[anchor=south west,xshift=-4pt,yshift=+0pt] {\trimfiga{#2}{\figWidtha}};
  \draw(4,.2) node[draw,fill=white,anchor=west,xshift=2pt,yshift=2pt,inner sep=3pt] {\scriptsize #3};
\begin{scope}[xshift=0cm,yshift=-5pt]
  \draw (\xcb,\ycb) node[anchor=south west,xshift=0.25cm,yshift=.5cm,rotate=-90] {\trimfigcb{colourBarLines}{\cbWidth}{\cbHeight}};
  \draw (.8,0) node[anchor=north,xshift=+3pt,yshift=+2pt] {\scriptsize $#6$};
  \draw (4.8,0) node[anchor=north,xshift=+0pt,yshift=+2pt] {\scriptsize $#7$};
\end{scope}
\end{scope}
}
\newcommand{\cbWidth}{.2cm}
\newcommand{\cbHeight}{4cm}
\newcommand{\xcb}{.5cm}
\newcommand{\ycb}{-.2cm}
\setlength{\ycbTop}{\ycb+\cbHeight}
\setlength{\ycbMid}{\ycb+\cbHeight*\real{.5}}
\newcommand{\trimfigcb}[3]{\includegraphics[width=#2, height=#3, clip, trim=17cm 2.35cm 1.65cm 2.35cm]{#1}}
\newcommand{\figWidtha}{5cm}
\newcommand{\trimfiga}[2]{\trimw{#1}{#2}{.15}{.15}{.11}{.11}}
\newcommand{\figSize}{5.5cm}
\begin{figure}[htb]

\begin{center}
\begin{tikzpicture}
   \useasboundingbox (0,.25) rectangle (16,6);  

  \begin{scope}[yshift=0.0cm] 
     \drawContour{xshift=-.75cm,yshift=0.00cm}{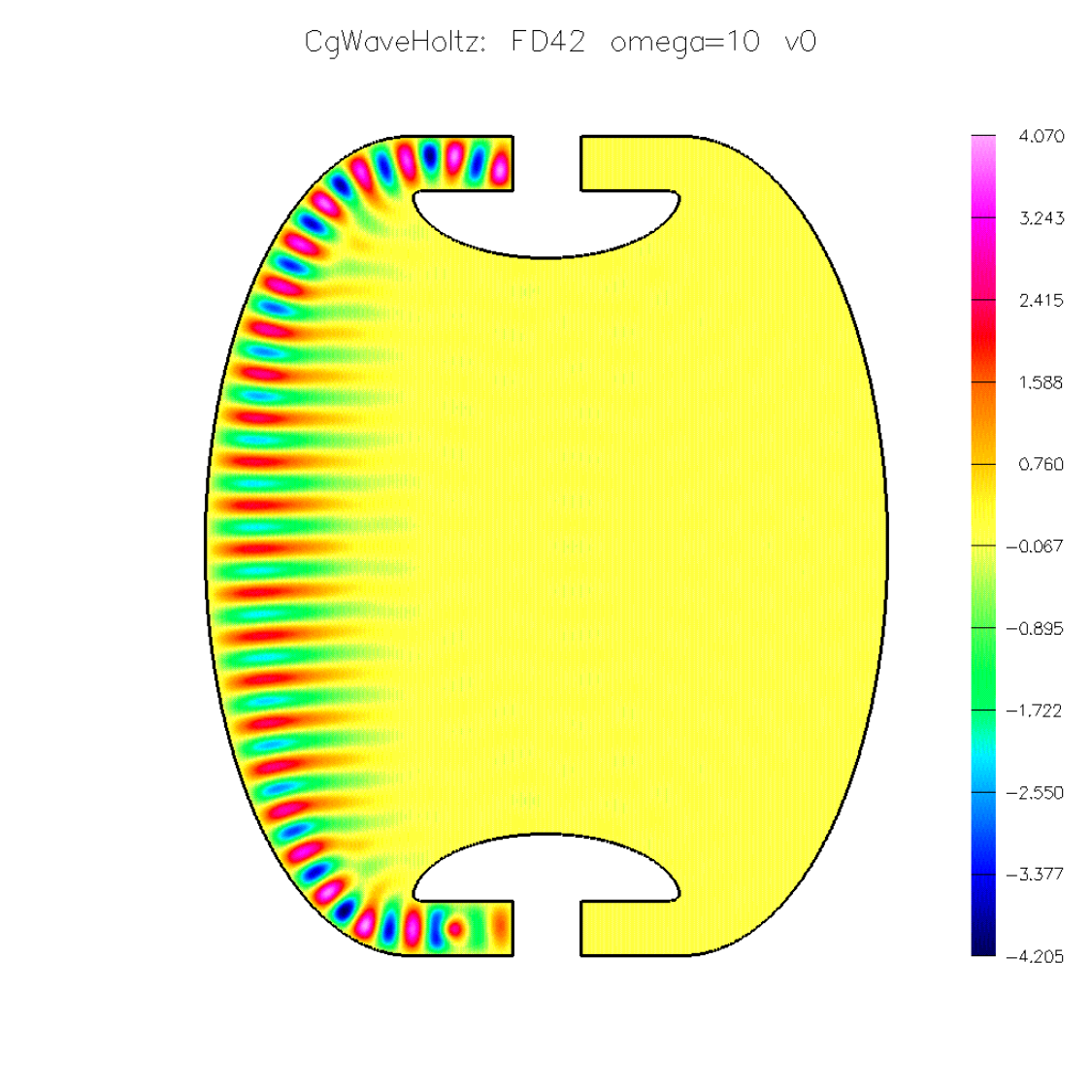}{G4 O4}{$v$}{$t=0.3$}{$-4.2$}{$4.1$}       
   \end{scope}  
   \begin{scope}[xshift=4.75cm,yshift=.7cm]
    \figByWidth{0.00}{0}{darkCornerRoomFreq10Np8Nit10G4Order4Deflate64}{\figSize}[0][0][0][0]
    \figByWidth{5.75}{0}{darkCornerRoomFreq10Np8Nit10G4Order4Deflate64FixPointMuFunction}{\figSize}[0][0][0][0]
    \draw (.75,.9) node[draw,fill=white,anchor=west,xshift=0pt,yshift=0pt,inner sep=3pt] {\scriptsize Double ellipse, O4, $\Gcde^{(4)}$};
  \end{scope}      

\end{tikzpicture}
\end{center}
\caption{Double ellipse Helmholtz solution for $\omega=10$, and Gaussian source located in the lower left alcove at $(-1.4,-5.6)$
The black vertical line on the right graph indicates the values of the adjusted frequencies $\omegaTilde$ used to correct for time discretization errors. 
}
\label{fig:doubleEllipseContoursFreq10}
\end{figure}
\begin{figure}[htb]

\begin{center}
\begin{tikzpicture}
   \useasboundingbox (0,.25) rectangle (16,6);  

  \begin{scope}[yshift=0.0cm]  
    \drawContour{xshift=-.75cm,yshift=0.00cm}{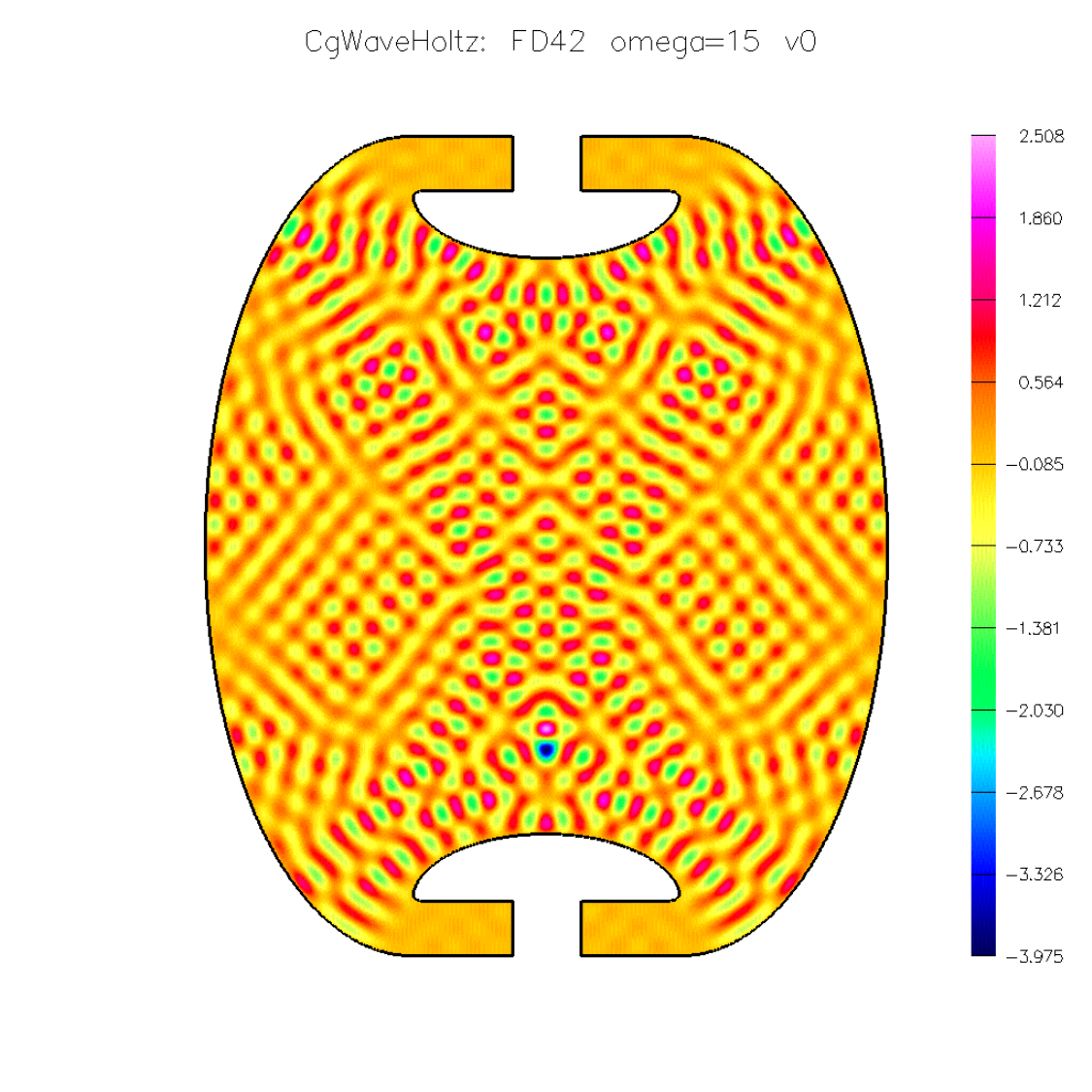}{G4 O4}{$v$}{$t=0.3$}{$-4.0$}{$2.5$}    
   \end{scope}  
   \begin{scope}[xshift=4.75cm,yshift=.7cm]
    \figByWidth{0.0}{0}{darkCornerRoomFreq15Np8Nit10G4Order4Deflate256}{\figSize}[0][0][0][0]
    \figByWidth{5.75}{0}{darkCornerRoomFreq15Np8Nit10G4Order4Deflate256FixPointMuFunction}{\figSize}[0][0][0][0]
    \draw (.75,.9) node[draw,fill=white,anchor=west,xshift=0pt,yshift=0pt,inner sep=3pt] {\scriptsize Double ellipse, O4, $\Gcde^{(4)}$};
  \end{scope}      

\end{tikzpicture}
\end{center}
\caption{Double ellipse Helmholtz solution for $\omega=15$, and Gaussian source located at $(0,-3)$, computed on grid $\Gcde^{(4)}$.
The black vertical line on the right graph indicates the value of the adjusted frequency $\omegaTilde$ used to correct for time discretization errors. 
}
\label{fig:doubleEllipseContoursFreq15}
\end{figure}
}

{
\newcommand{\figSize}{6.5cm}

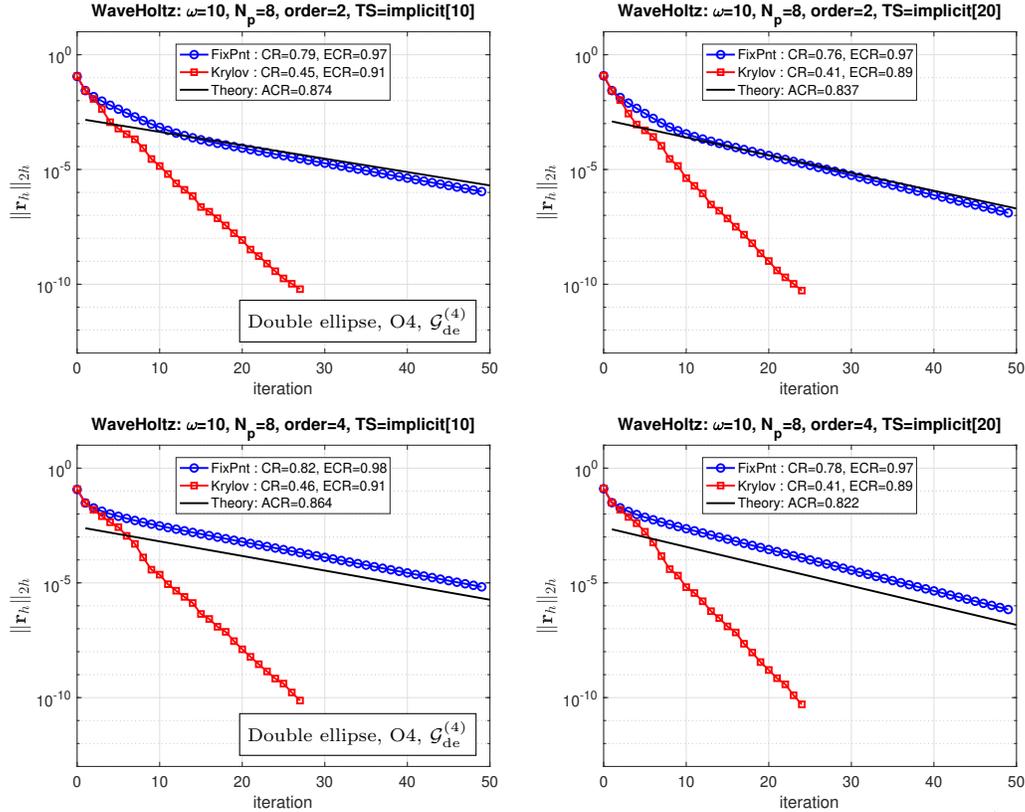
\begin{figure}[htb]
\begin{center}
\begin{tikzpicture}[scale=1]
  \useasboundingbox (0,.8) rectangle (14,11);  
   \begin{scope}[yshift=5.5cm]
    \figByWidth{  0}{0}{darkCornerRoomFreq10Np8Nit10G4Order2Deflate64}{\figSize}[0][0][0][0]
    \figByWidth{7.0}{0}{darkCornerRoomFreq10Np8Nit20G4Order2Deflate64}{\figSize}[0][0][0][0]
    \draw (3,1) node[draw,fill=white,anchor=west,xshift=0pt,yshift=0pt,inner sep=3pt] {\scriptsize Double ellipse, O4, $\Gcde^{(4)}$};
  \end{scope}
   \begin{scope}[yshift=0cm]
    \figByWidth{  0}{0}{darkCornerRoomFreq10Np8Nit10G4Order4Deflate64}{\figSize}[0][0][0][0]
    \figByWidth{7.0}{0}{darkCornerRoomFreq10Np8Nit20G4Order4Deflate64}{\figSize}[0][0][0][0]
    \draw (3,1) node[draw,fill=white,anchor=west,xshift=0pt,yshift=0pt,inner sep=3pt] {\scriptsize Double ellipse, O4, $\Gcde^{(4)}$};

  \end{scope}  

\end{tikzpicture}
\end{center}
\caption{Double ellipse with deflation, implicit time-stepping, $\omega=10$. Comparing results for $\Nits=10$ (left column) and $\Nits=20$ (right column).
         Top: order=2. Bottom: order=4. The convergence is slightly better for $\Nits=20$ but at roughly double the cost in CPU time.
         The black vertical line on the right graph indicates the value of the adjusted frequency $\omegaTilde$ used to correct for time discretization errors. 
    } 
\label{fig:doubleEllipseFreq10Deflate64CompareNits}
\end{figure}
}

{
\newcommand{\drawContour}[7]{%
\begin{scope}[#1]
\draw(0.0,0) node[anchor=south west,xshift=-4pt,yshift=+0pt] {\trimfiga{#2}{\figWidtha}};
  \draw(2,3.0) node[draw,fill=white,anchor=west,xshift=2pt,yshift=2pt,inner sep=2pt] {\scriptsize #3};
  \draw(2,2.5) node[draw,fill=white,anchor=west,xshift=2pt,yshift=2pt,inner sep=2pt] {\scriptsize #4};
\begin{scope}[xshift=0cm,yshift=-5pt]
  \draw (\xcb,\ycb) node[anchor=south west,xshift=0.25cm,yshift=.5cm,rotate=-90] {\trimfigcb{colourBarLines}{\cbWidth}{\cbHeight}};
  \draw (.8,0) node[anchor=north,xshift=+3pt,yshift=+2pt] {\scriptsize $#6$};
  \draw (4.8,0) node[anchor=north,xshift=+0pt,yshift=+2pt] {\scriptsize $#7$};
\end{scope}
\end{scope}
}
\newcommand{\cbWidth}{.2cm}
\newcommand{\cbHeight}{4cm}
\newcommand{\xcb}{.5cm}
\newcommand{\ycb}{-.2cm}
\setlength{\ycbTop}{\ycb+\cbHeight}
\setlength{\ycbMid}{\ycb+\cbHeight*\real{.5}}
\newcommand{\trimfigcb}[3]{\includegraphics[width=#2, height=#3, clip, trim=17cm 2.35cm 1.65cm 2.35cm]{#1}}
\newcommand{\figWidtha}{5cm}
\newcommand{\trimfiga}[2]{\trimw{#1}{#2}{.15}{.15}{.11}{.11}}
\begin{figure}[htb]
\begin{center}
\begin{tikzpicture}
   \useasboundingbox (0,.25) rectangle (16,12);  

  \begin{scope}[yshift=6.25cm]  
    \drawContour{xshift=0.0cm,yshift=0.00cm}{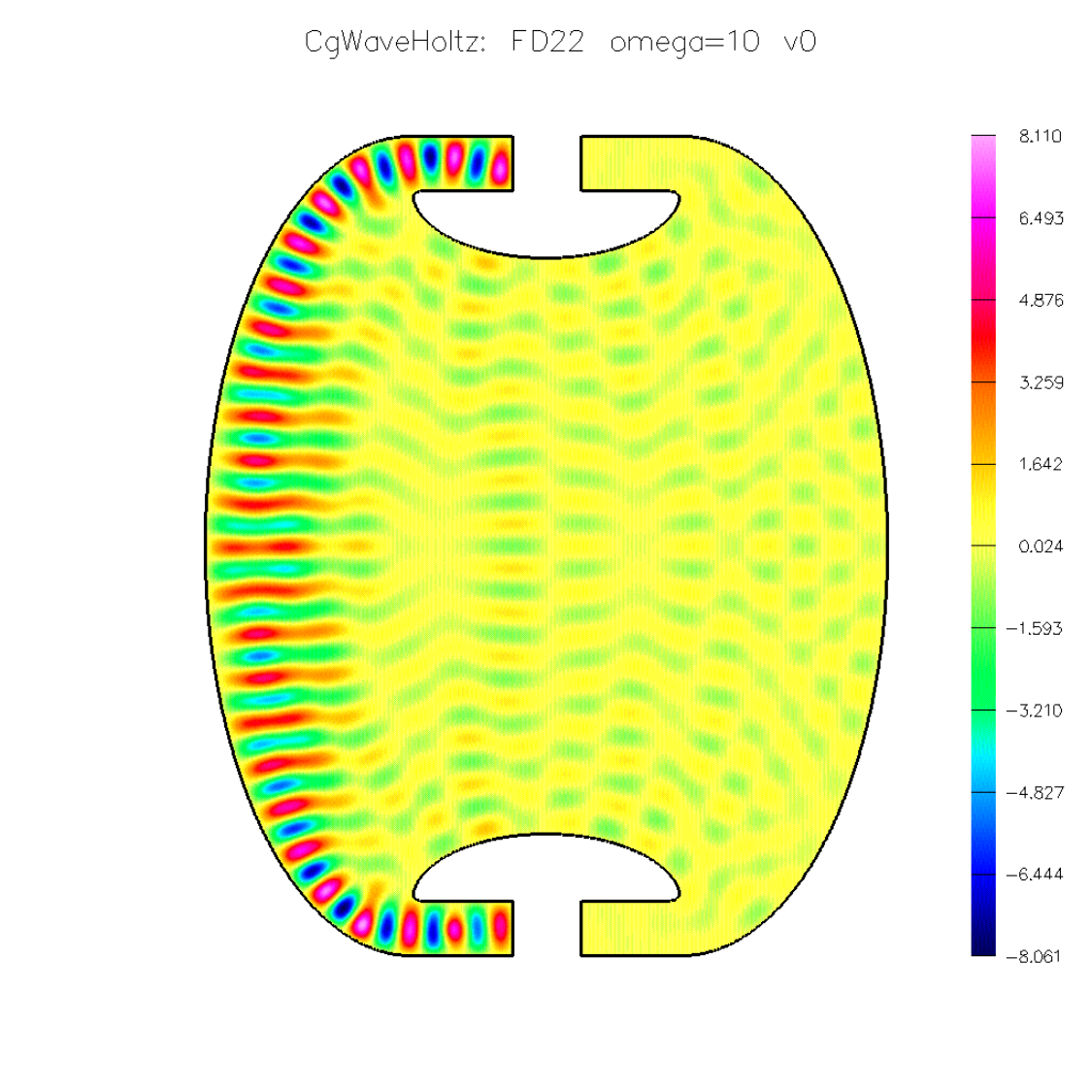}{O2 $\Gcde^{4}$ $\PPW=25$}{$\PPW_2=141$}{$t=0.3$}{$-8.06$}{$8.11$} 
    \drawContour{xshift=5.0cm,yshift=0.00cm}{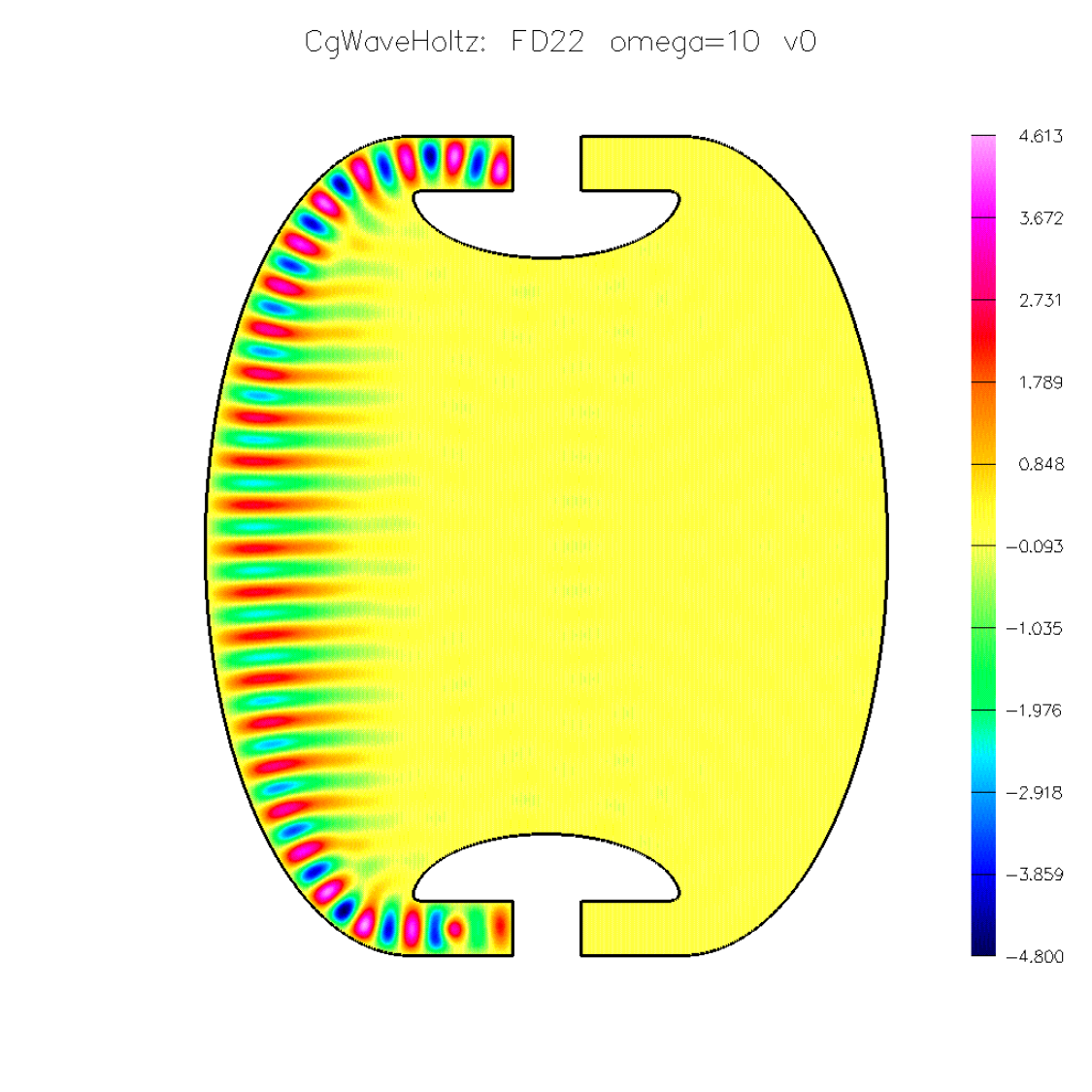}{O2 $\Gcde^{8}$ $\PPW=50$}{$\PPW_2=141$}{$t=0.3$}{$-4.80$}{$4.61$}
    \drawContour{xshift=10.cm,yshift=0.00cm}{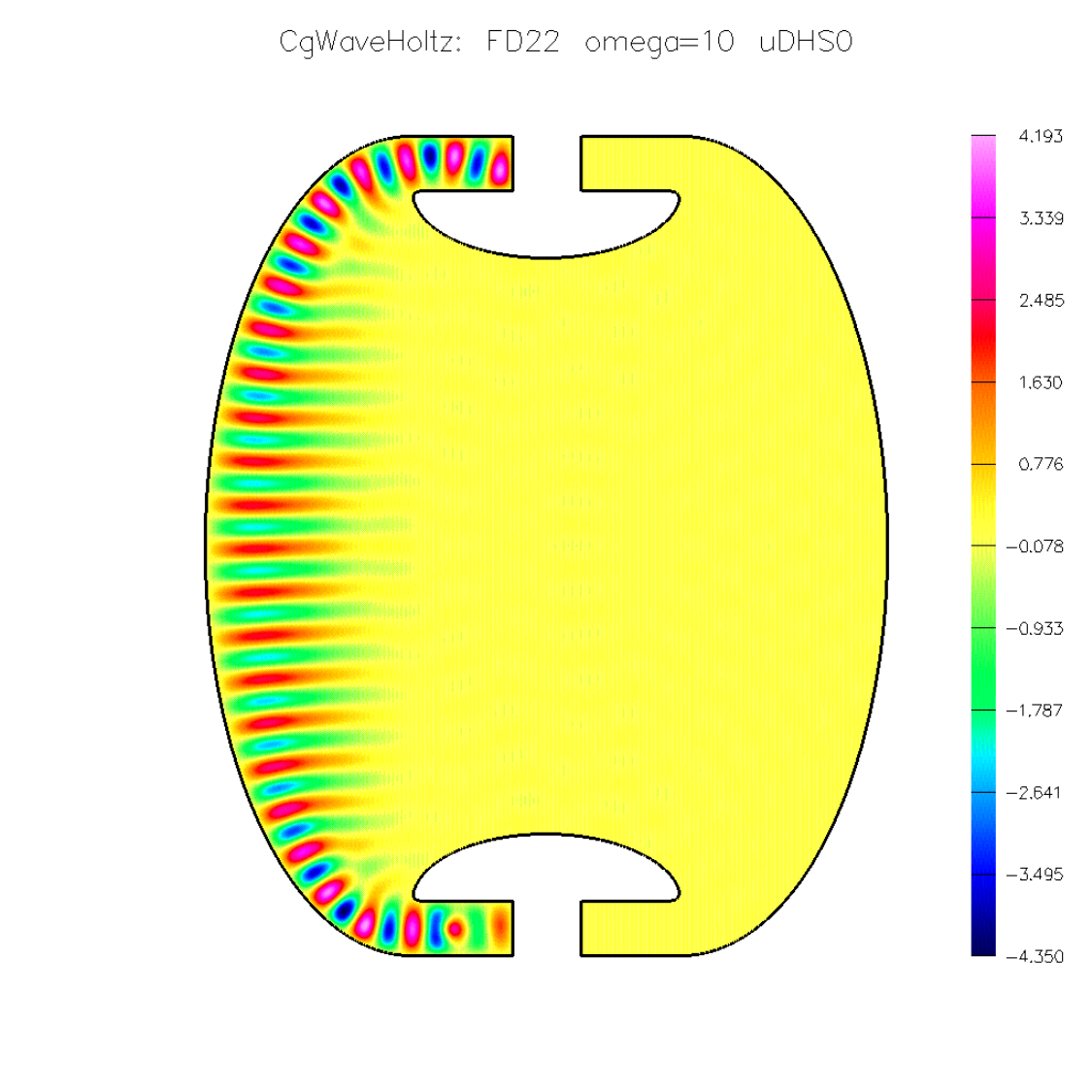}{O2 $\Gcde^{16}$ $\PPW=100$}{$\PPW_2=141$}{$t=0.3$}{$-4.35$}{$4.19$}    
   \end{scope}   
  \begin{scope}[yshift=0.0cm]  
    \drawContour{xshift=0.0cm,yshift=0.00cm}{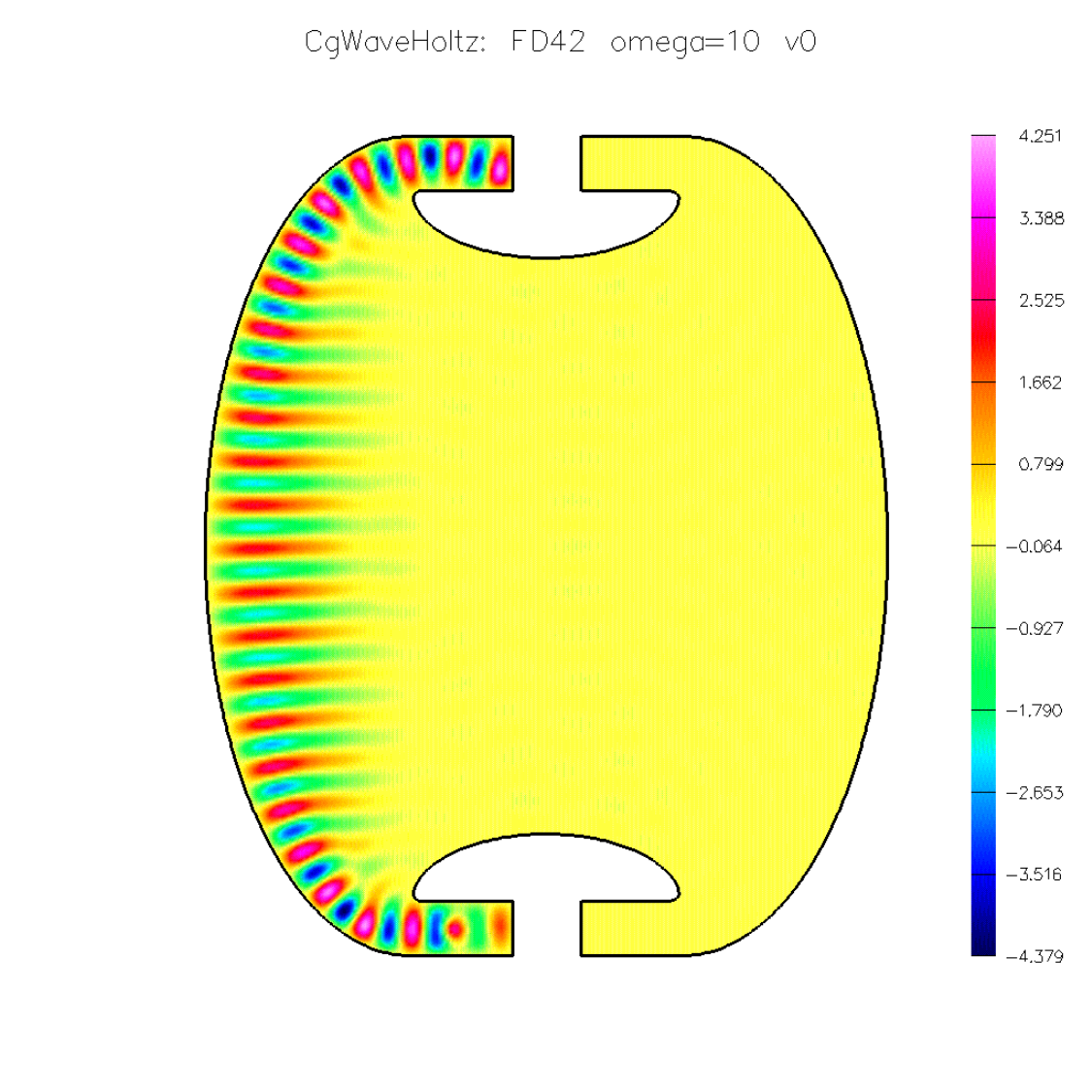}{O4 $\Gcde^{2}$ $\PPW=13$}{$\PPW_4=18$}{$t=0.3$}{$-4.38$}{$4.25$}    
    \drawContour{xshift=5.0cm,yshift=0.00cm}{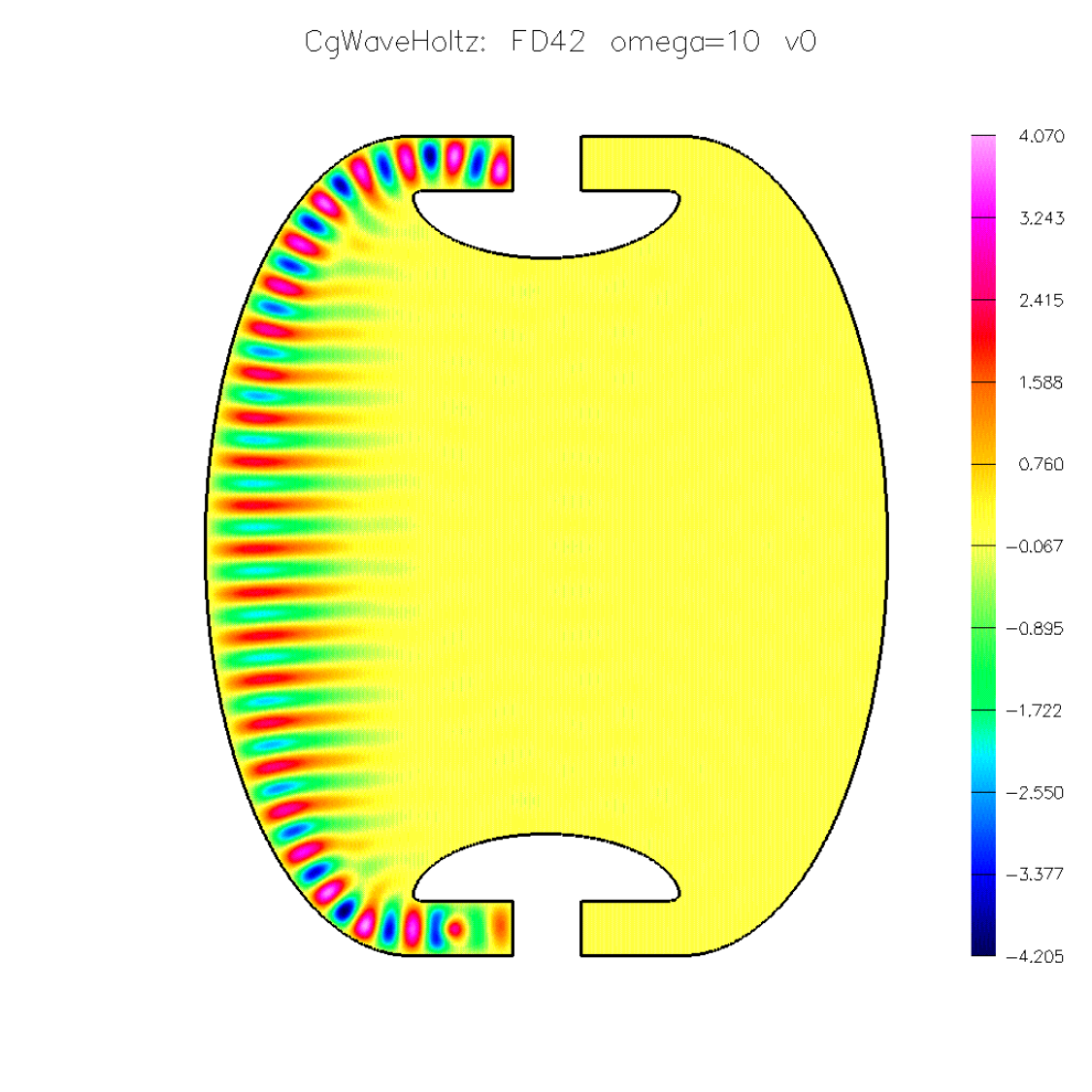}{O4 $\Gcde^{4}$ $\PPW=25$}{$\PPW_4=18$}{$t=0.3$}{$-4.21$}{$4.07$}     
    \drawContour{xshift=10.cm,yshift=0.00cm}{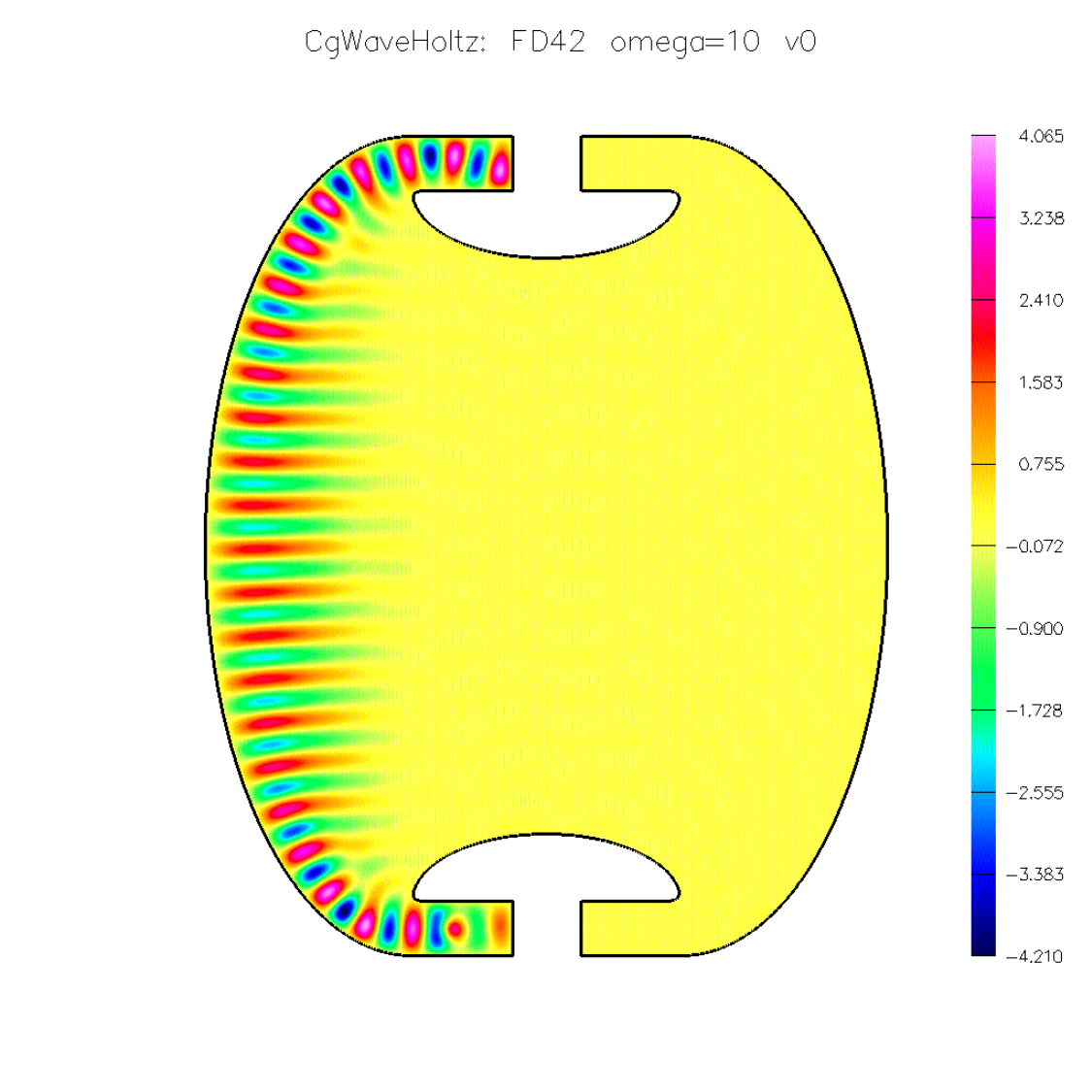}{O4 $\Gcde^{8}$ $\PPW=50$}{$\PPW_4=18$}{$t=0.3$}{$-4.21$}{$4.07$}     
   \end{scope}   

\end{tikzpicture}
\end{center}
\caption{Double ellipse Helmholtz solutions, $\omega=10$. Comparing second-order to fourth-order accurate results. Top row: order two. Bottom row: order four.
  The grid $\Gcde^{(j)}$ has an approximate grid spacing of $\ds^{(j)}=1/(10 j)$. 
  Much finer grids are required for the second-order accurate scheme in order to match the results from the fourth-order accurate scheme (note max and min values on the colour bars).
  $\PPW$ denotes the actual points-per-wavelength used. $\PPW_2$ and $\PPW_4$ are the estimated rule-of-thumb values.
  }
\label{fig:doubleEllipseContoursOrderComparison}
\end{figure}
}

{
\newcommand{\drawContour}[7]{%
\begin{scope}[#1]
\draw(0.0,0) node[anchor=south west,xshift=-4pt,yshift=+0pt] {\trimfiga{#2}{\figWidtha}};
  \draw(2,6.5) node[draw,fill=white,anchor=west,xshift=2pt,yshift=2pt,inner sep=2pt] {\scriptsize #3};
\begin{scope}[xshift=0.3cm,yshift=-5pt]
  \draw (\xcb,\ycb) node[anchor=south west,xshift=0.25cm,yshift=.5cm,rotate=-90] {\trimfigcb{colourBarLines}{\cbWidth}{\cbHeight}};
  \draw (.8,0) node[anchor=north,xshift=+3pt,yshift=+2pt] {\scriptsize $#6$};
  \draw (4.85,0) node[anchor=north,xshift=+0pt,yshift=+2pt] {\scriptsize $#7$};
\end{scope}
\end{scope}
}
\newcommand{\cbWidth}{.2cm}
\newcommand{\cbHeight}{4cm}
\newcommand{\xcb}{.5cm}
\newcommand{\ycb}{-.2cm}
\setlength{\ycbTop}{\ycb+\cbHeight}
\setlength{\ycbMid}{\ycb+\cbHeight*\real{.5}}
\newcommand{\trimfigcb}[3]{\includegraphics[width=#2, height=#3, clip, trim=17cm 2.35cm 1.65cm 2.35cm]{#1}}
\newcommand{\figWidtha}{5.5cm}
\newcommand{\trimfiga}[2]{\trimw{#1}{#2}{.15}{.15}{.11}{.11}}
\newcommand{\figSize}{5.5cm}
\begin{figure}[htb]

\begin{center}
\begin{tikzpicture}
   \useasboundingbox (0,.25) rectangle (11,6.7);  

  \begin{scope}[xshift=-.75cm,yshift=0.0cm] 
     \drawContour{xshift= 0.0cm,yshift=0.00cm}{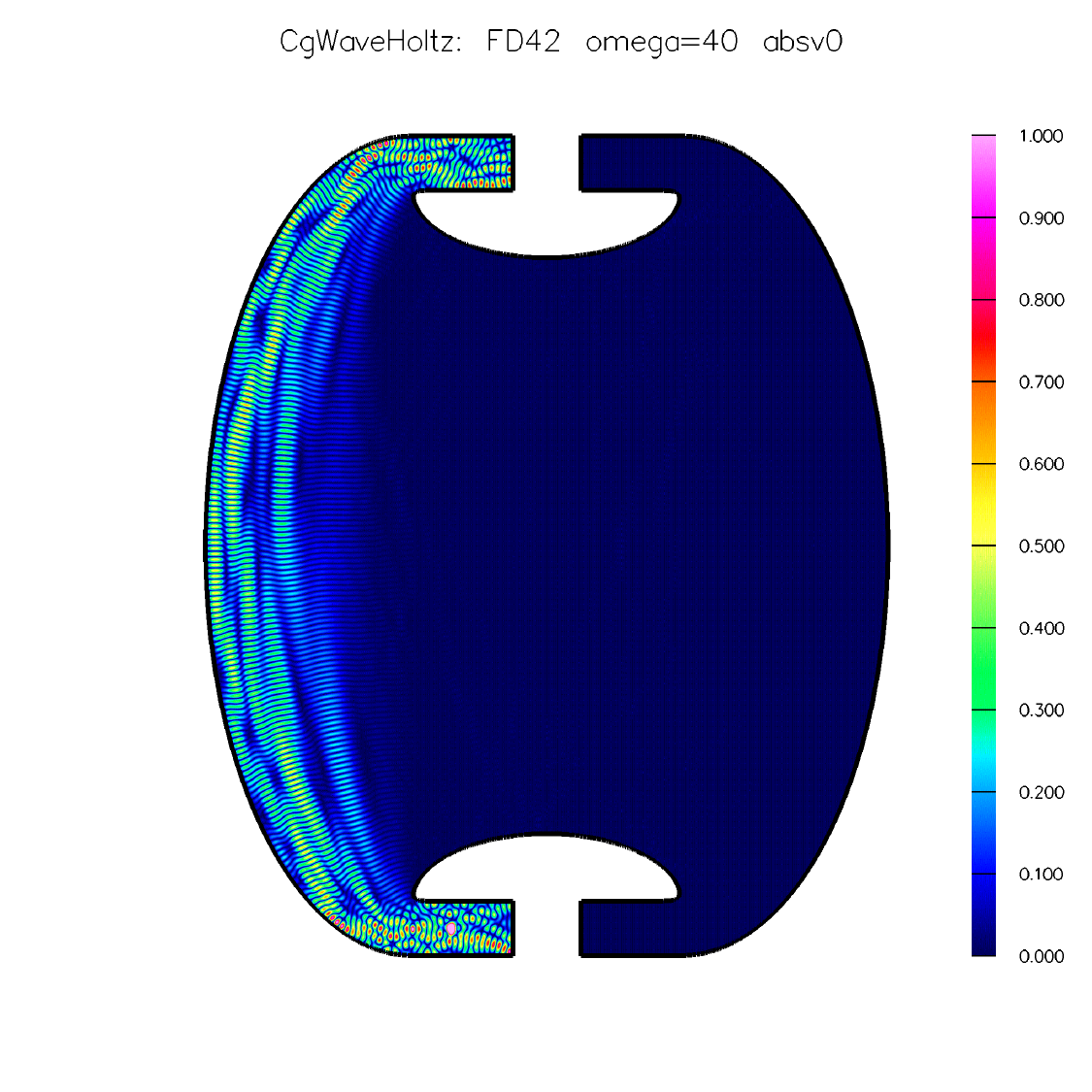}{$|v|$, $\omega=40$, $\Gc^{(16)}$}{$v$}{$t=0.3$}{$0$}{$1$}    
     \drawContour{xshift=5.75cm,yshift=0.00cm}{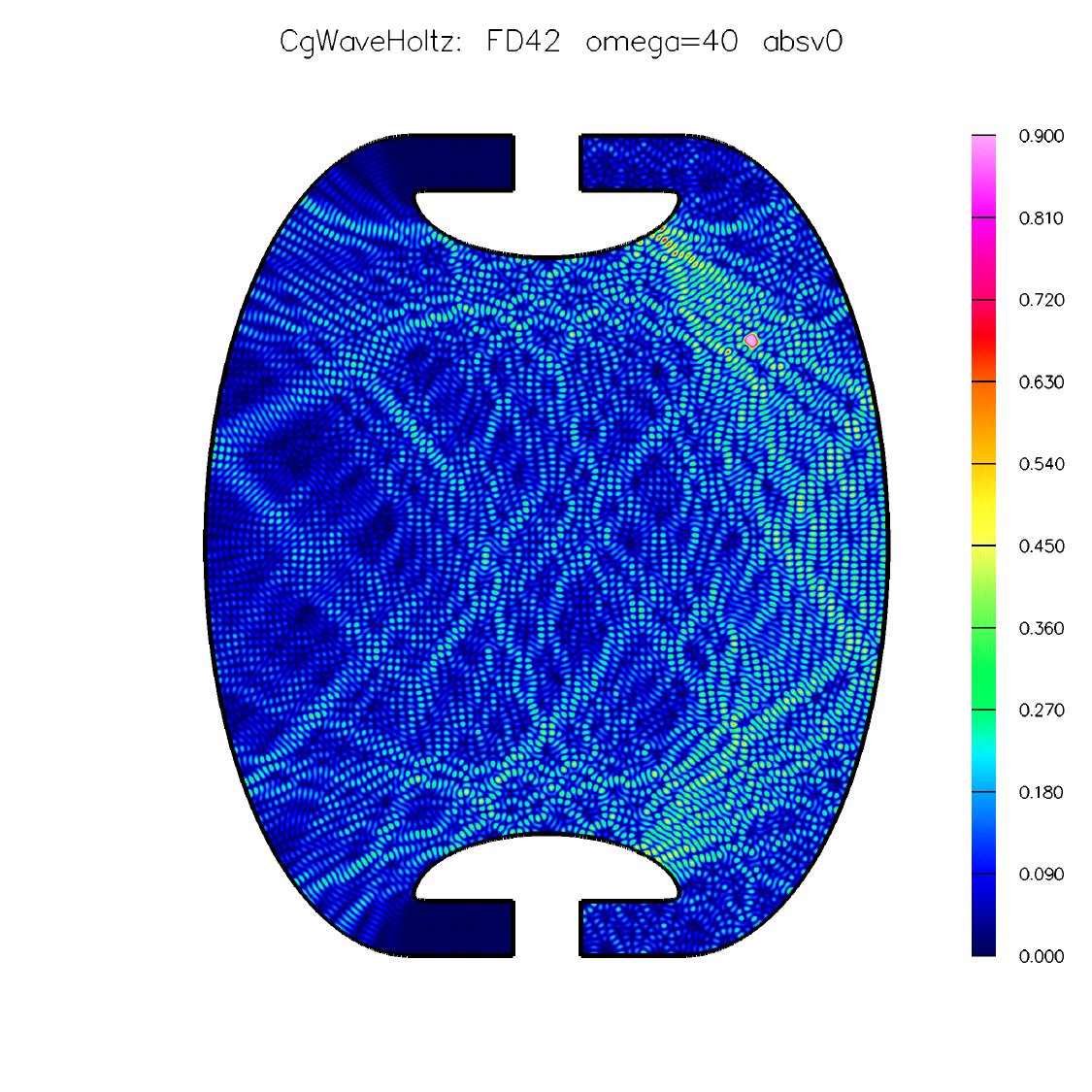}{$|v|$, $\omega=40$, $\Gc^{(16)}$}{$v$}{$t=0.3$}{$0$}{$0.9$}  
   \end{scope}  
\end{tikzpicture}
\end{center}
\caption{Double ellipse Helmholtz solution, $\omega=40$. 
Left: Gaussian source located in the lower left alcove at $(-1.4,-5.6)$ leads to a surface mode on the left edge.
Right: Gaussian source located at $(3,3)$ leads to two \textsl{quiet} alcoves in the top left and bottom right. 
}
\label{fig:doubleEllipseContoursHighFreq}
\end{figure}

}

Figure~\ref{fig:doubleEllipseContoursFreq10} shows results for computations on  grid $\Gcde^{(4)}$
using the Gaussian source with $\omega=10$, $\alphag=400$ and $\betag=10$, and centered at a point $\xv_0=(-1.4,-5.6)$ located
in the lower left alcove. For this choice, a surface mode is generated lying along the left boundary and entering the upper left alcove. 
 The convergence of the WaveHoltz FPI and GMRES iterations
are shown in the graphs in the middle and right of Figure~\ref{fig:doubleEllipseContoursFreq10}.
The results are generated using the fourth-order accurate implicit scheme with $\Np=8$ periods per time interval, $10$ implicit time-steps per period, and $64$ deflated eigenmodes. The FPI convergence rate is seen to match the theory with the GMRES converging rapidly ($\CR \approx 0.27$ and $\ECR\approx 0.85$).
From Table~\ref{tab:doubleEllipsePPW} the actual $\PPW$ is $25$ while the estimated $\PPW_4$ is $18$ which suggests that the computation is resolved (to within a relative error tolerance of
$\eps=10^{-2}$ used for the values in Table~\ref{tab:doubleEllipsePPW}).
The accuracy of the calculations are confirmed using a grid convergence study as discussed below.

Figure~\ref{fig:doubleEllipseContoursFreq15} shows results for a somewhat higher frequency.
The Gaussian source~\eqref{eq:gaussianSource} for this calculation uses $\omega=15$, $\alphag=400$ and $\betag=5$, and a source centered at $\xv_0=(0,-3)$ which is located in the lower center of the interior.
 The convergence of the WaveHoltz FPI and GMRES iterations
are shown in the graphs in the middle and right of Figure~\ref{fig:doubleEllipseContoursFreq15}.
The fourth-order accurate implicit scheme is used with $\Np=8$ periods per time interval, $10$ time-steps per period, and $270$ eigenmodes 
deflated. The FPI convergence rate is seen to match the theory with the GMRES converging rapidly.
For this computation, performed on grid $\Gcde^{(4)}$, the actual $\PPW=17$ while the estimated $\PPW_4=20$ (from Table~\ref{tab:doubleEllipsePPW}), suggesting that the simulation
 is reasonably resolved based on the rule-of-thumb estimates.

Figure~\ref{fig:doubleEllipseFreq10Deflate64CompareNits} compares results using $\Nits=10$ and $\Nits=20$ 
time-steps per-period for both second-order and fourth-order accurate calculations.
The grid for these results is $\Gcde^{(4)}$ and the Gaussian source parameters are the same as those used for the results in Figure~\ref{fig:doubleEllipseContoursFreq10}.  
Using $\Nits=20$ leads to slightly better convergence but at roughly double the cost in CPU time per WaveHoltz iteration. 
Thus using $\Nits=10$ would appear to be the more efficient option in this case. 
We note that the most efficient value for $\Nits$ could depend on the problem being solved.

\begin{table}[hbt]\tableFont 
\begin{center}
\begin{tabular}{|l|c|c|} \hline 
  \multicolumn{3}{|c|}{Double Ellipse Order 2}  \\ \hline
   Grid        & $\| E_v \|_\infty$    &  rate    \\ \hline
 \strutt $\Gcde^{(4)}$  & \num{4.4}{0}  &    \\ \hline
 \strutt $\Gcde^{(8)}$  & \num{6.2}{-1} & 2.8  \\ \hline
 \strutt $\Gcde^{(16)}$ & \num{8.7}{-2} & 2.8  \\ \hline
\end{tabular}
\qquad
\begin{tabular}{|l|c|c|} \hline 
  \multicolumn{3}{|c|}{Double Ellipse Order 4}  \\ \hline
   Grid        & $\| E_v \|_\infty$    &  rate    \\ \hline
 \strutt $\Gcde^{(2)}$ & \num{2.4}{-1} &        \\ \hline
 \strutt $\Gcde^{(4)}$ & \num{9.0}{-3} & 4.8   \\ \hline
 \strutt $\Gcde^{(8)}$ & \num{3.3}{-4} & 4.8   \\ \hline
\end{tabular}
\caption{Double ellipse grid self-convergence. 
  Estimated max-norm errors and convergence rates for the computations shown in Figure~\ref{fig:doubleEllipseContoursOrderComparison}.
}
\label{tab:doubleEllipseSelfConvergence}
\end{center}
\end{table}

Figure~\ref{fig:doubleEllipseContoursOrderComparison} compares results from second-order and fourth-order accurate computations using grids
of varying resolutions. 
It can be seen that much finer grids are required for the second-order accurate scheme in order to match the results from the fourth-order accurate scheme.
The second-order accurate results on grid $\Gcde^{(16)}$ are comparable to the fourth-order accurate results on $\Gcde^{(2)}$, the former grid having a grid spacing $8$ times finer ($64$ times more grid points).
The contour plots also note the values of $\PPW$ (actual points-per-wavelength) along with the estimated values $\PPW_2$ and $\PPW_4$ from Recipe~\ref{recipe:PPW}.
The rule-of-thumb values are seen to provide reasonably good estimates for the required PPW.
According to $\PPW_4=13$ the fourth-order accurate results are nearly resolved on grid $\Gcde^{(2)}$ and well resolved on grid $\Gcde^{(4)}$,
while the second-order accurate results, with $\PPW_2=141$, are just beginning to be resolved on the finest grid $\Gcde^{(16)}$.
The rule-of-thumb values are thus seen to provide good estimates.
Also note, in terms of performance of the the implicit time-stepping scheme, on grid $\Gcde^{(16)}$ with $\Nits=10$, the time-step that was about $230$ 
larger than that for explicit time-stepping.

To estimate the actual errors in the computations shown in Figure~\ref{fig:doubleEllipseContoursOrderComparison}, a grid self-convergence study is performed.
Given computed solutions on three grid resolutions the errors and convergence rates can be estimated using a Richardson extrapolation procedure described in~\cite{pog2008a}.
Table~\ref{tab:doubleEllipseSelfConvergence} shows the estimated max-norm errors, $\| E_v\|_\infty$,  and convergence rates for the second-order and fourth-order accurate schemes.
The estimated errors for second-order accuracy are converging at a rate somewhat better than $2$; however the errors on the coarse and medium resolution grids are relatively large.
Consistent with the previous observations from Figure~\ref{fig:doubleEllipseContoursOrderComparison}, grid $\Gcde^{(16)}$ is still not really fine enough for this second-order accurate computation.
The estimated errors for the fourth-order accurate scheme are seen to be converging at a rate somewhat better than $4$.
At fourth-order accuracy the estimated error of $0.24$ for grid $\Gcde^{(2)}$  indicates that the computed solution is under-resolved,
while the error of $0.009$ for grid $\Gcde^{(4)}$ points to a reasonably resolved calculation. These results are consistent with the 
rule of thumb suggestions for the points-per-wavelength as discussed in the previous paragraph.

Finally, Figure~\ref{fig:doubleEllipseContoursHighFreq} shows results for the higher frequency $\omega=40$ performed using grid $\Gcde^{(16)}$.
The Gaussian source~\eqref{eq:gaussianSource} for these simulations use a weight of $\alphag=2000$ and an exponent of $\betag=10$.
For the results shown on the left, the Gaussian source is located in the lower left alcove at $\xv_0=(-1.4,-5.6)$.
This source generates a surface wave that is located primarily on the left-hand side of the domain.
The simulation shown on the right of Figure~\ref{fig:doubleEllipseContoursHighFreq} is computed with the Gaussian source
located at $\xv_0=(3,3)$. In this case the solution in the lower left and upper left alcoves are relatively quiet with $\vert v\vert$ nearly zero.
From Table~\ref{tab:doubleEllipsePPW} the actual PPW corresponding to the grid $\Gcde^{(16)}$ is approximately $25$ while the rule of thumb estimate
is also approximately $25$ which suggests that these simulations are reasonably resolved.
In contrast, $\PPW_2=281$ for a second-order accurate scheme which would require a grid spacing over $10$ times smaller and a grid with over $100$ times more grid points for this two-dimensional simulation.

\newcommand{\Gcke}{\Gc_{ke}}
\subsection{Knife edge}  \label{sec:knifeEdge}

{
\newcommand{\drawContour}[7]{%
\begin{scope}[#1]
\draw(0.0,0) node[anchor=south west,xshift=-4pt,yshift=+0pt] {\trimfiga{#2}{\figWidtha}};
\draw(3,4.4) node[draw,fill=white,anchor=west,xshift=2pt,yshift=1pt,inner sep=2pt] {\scriptsize #3};
\begin{scope}[xshift=0cm,yshift=-2pt]
  \draw (\xcb,\ycb) node[anchor=south west,xshift=0.25cm,yshift=.5cm,rotate=-90] {\trimfigcb{colourBarLines}{\cbWidth}{\cbHeight}};
  \draw (.8,0) node[anchor=north,xshift=+3pt,yshift=+2pt] {\scriptsize $#6$};
  \draw (6.8,0) node[anchor=north,xshift=+0pt,yshift=+2pt] {\scriptsize $#7$};
\end{scope}
\end{scope}
}
\newcommand{\cbWidth}{.2cm}
\newcommand{\cbHeight}{6cm}
\newcommand{\xcb}{.5cm}
\newcommand{\ycb}{-.2cm}
\setlength{\ycbTop}{\ycb+\cbHeight}
\setlength{\ycbMid}{\ycb+\cbHeight*\real{.5}}
\newcommand{\trimfigcb}[3]{\includegraphics[width=#2, height=#3, clip, trim=17cm 2.35cm 1.65cm 2.35cm]{#1}}
\newcommand{\figWidtha}{7cm}
\newcommand{\trimfiga}[2]{\trimw{#1}{#2}{.12}{.12}{.28}{.28}}
\begin{figure}[htb]
\begin{center}
\begin{tikzpicture}
   \useasboundingbox (-.5,.3) rectangle (15.5,4.75);  

   \figByHeight{0.0}{0}{tipGridSmall64}{4.5cm}[0.03][0.03][0.24][0.24]
   \figByHeightb{.5}{.5}{tipGridSmall64Zoom1}{3.5cm}[0.25][0.25][0.1][0.2]
   \figByHeightb{4.75}{.75}{tipGridSmall64Zoom2}{3.cm}[0.1][0.1][0.1][0.1]

  \draw[thick,black,->,yshift=0pt] (0,0) -- (8.4,0.00);  
  \draw[thick,black,->,xshift=-1pt] (0,0) -- (0.00,4.8);     

   \draw[thick,black,-,xshift=0pt] (0.00,.1) -- (0.00,-.1) node[anchor=north] {\scriptsize $-.5$};
   \draw[thick,black,-,xshift=0pt] (4.07,.1) -- (4.07,-.1) node[anchor=north] {\scriptsize $0$};
   \draw[thick,black,-,xshift=0pt] (8.10,.1) -- (8.10,-.1) node[anchor=north] {\scriptsize $0.5$};

   \draw[thick,black,-,xshift=0pt] (.1,0.0) -- (-.1,0.0) node[anchor=east] {\scriptsize $0$};
   \draw[thick,black,-,xshift=0pt] (.1,4.5) -- (-.1,4.5) node[anchor=east] {\scriptsize $.55$};

   \begin{scope}[xshift=8.cm,yshift=0cm]
     \drawContour{xshift=0.0cm,yshift=0.cm}{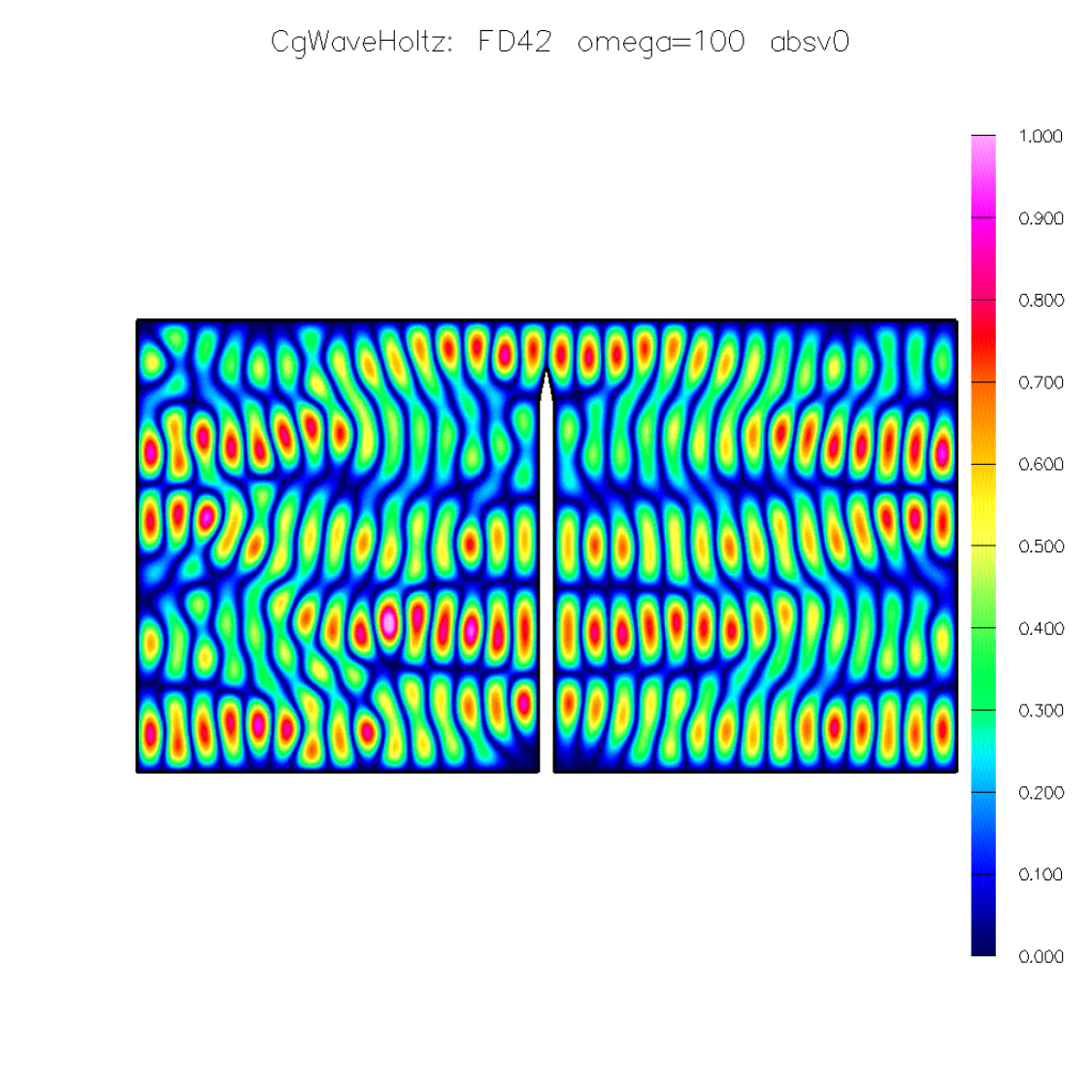}{$|v|$, $\omega=100$}{$v$}{$t=1.0$}{$0.$}{$1.$}     
   \end{scope}   

\end{tikzpicture}
\end{center}
\caption{
    Knife edge. Left: overset grid $\Gcke^{(32)}$ for a knife edge with inserts showing magnified regions near the knife edge and tip. Right: absolute value of the computed Helmholtz solution for $\omega=100$ with a Gaussian source at $(-.2,.2)$.
    }
\label{fig:knifeEdgeGridAndContours}
\end{figure}
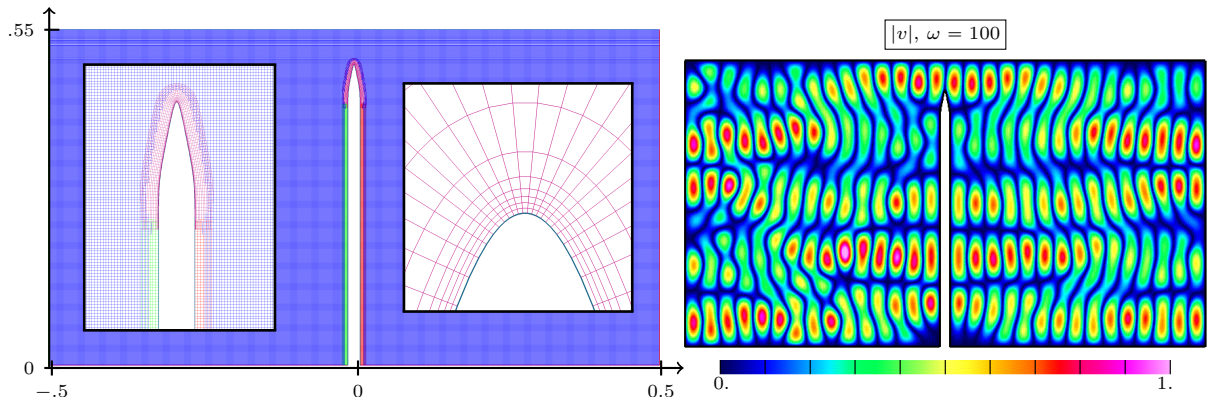
}

{
\newcommand{\figSize}{5.75cm}
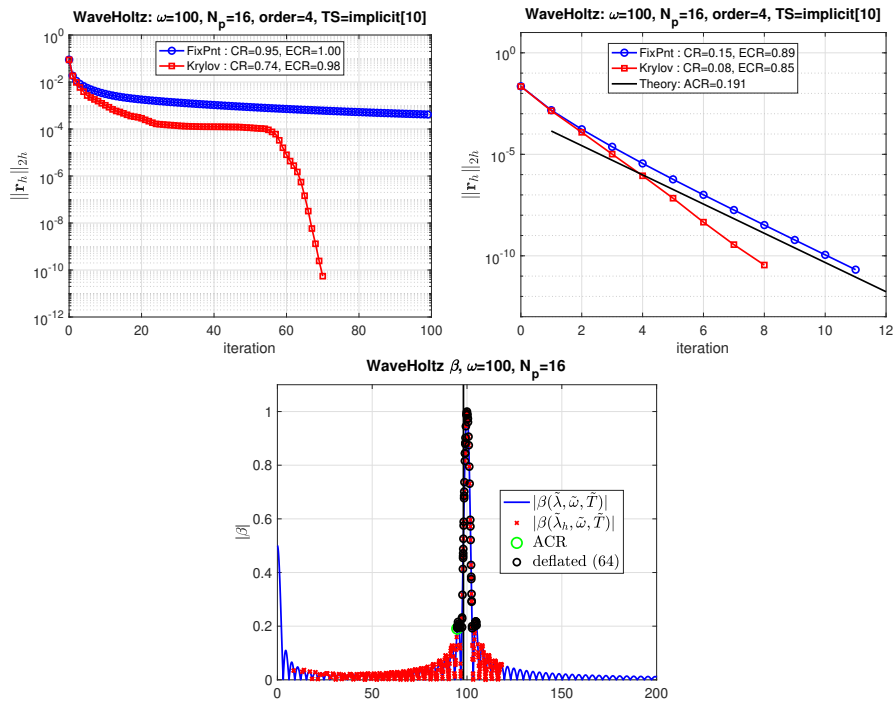
\begin{figure}[htb]

\begin{center}
\begin{tikzpicture}
   \useasboundingbox (0,.8) rectangle (12,9.2);  

   \begin{scope}[xshift=0cm,yshift=0cm]
    \figByWidth{0.0}{4.85}{tipGridSmallG64O4Freq100}{\figSize}[0][0][0][0]
    \figByWidth{6.0}{4.85}{tipGridSmallG64O4Freq100Deflate64}{\figSize}[0][0][0][0]
    \figByWidth{3.0}{0}{tipGridSmallG64O4Freq100Deflate64FixPointMuFunction}{\figSize}[0][0][0][0]
    
  \end{scope}      

\end{tikzpicture}
\end{center}
\caption{Knife edge convergence results for $\omega=100$, on grid $\Gcke^{(64)}$, order four. Top left: convergence history (no deflation). Top right: convergence history when deflating $64$ eigenmodes. Bottom: WaveHoltz filter function $\beta$ function with eigenvalues and deflated eigenvalues.
  }
\label{fig:knifeEdgeConvergenceFreq100}
\end{figure}
}

This example illustrates a Helmholtz problem for which implicit time-stepping is particularly useful.
The geometry contains a thin knife edge that requires a fine grid to resolve the sharp (rounded) tip as illustrated
in Figure~\ref{fig:knifeEdgeGridAndContours}. With a standard explicit time-stepping scheme, the maximum stable time-step is determined
by the smallest cells on the overall grid, and this requires a small global time-step even on grids where such a small time-step is not needed for stability.
One improvement would be to use a local time-stepping method or a locally implicit method such as the one as described in~\cite{carson2024highorderaccurateimplicitexplicittimestepping},
but even then the explicit time step scales at $\dt \sim h/c$ as $h$ goes to zero. In contrast, 
the implicit scheme used here employs a fixed number of time-steps per period, independent of~$h$.

{
\newcommand{\drawContour}[7]{%
\begin{scope}[#1]
\draw(0.0,0) node[anchor=south west,xshift=-4pt,yshift=+0pt] {\trimfiga{#2}{\figWidtha}};
  \draw(4.5,4.) node[draw,fill=white,anchor=west,xshift=2pt,yshift=2pt,inner sep=2pt] {\scriptsize #3};
\begin{scope}[xshift=0cm,yshift=-5pt]
  \draw (\xcb,\ycb) node[anchor=south west,xshift=0.25cm,yshift=.5cm,rotate=-90] {\trimfigcb{colourBarLines}{\cbWidth}{\cbHeight}};
  \draw (.8,0) node[anchor=north,xshift=+3pt,yshift=+2pt] {\scriptsize $#6$};
  \draw (7.1,0) node[anchor=north,xshift=+0pt,yshift=+2pt] {\scriptsize $#7$};
\end{scope}
\end{scope}
}
\newcommand{\cbWidth}{.2cm}
\newcommand{\cbHeight}{6.25cm}
\newcommand{\xcb}{.5cm}
\newcommand{\ycb}{-.2cm}
\setlength{\ycbTop}{\ycb+\cbHeight}
\setlength{\ycbMid}{\ycb+\cbHeight*\real{.5}}
\newcommand{\trimfigcb}[3]{\includegraphics[width=#2, height=#3, clip, trim=17cm 2.35cm 1.65cm 2.35cm]{#1}}
\newcommand{\figWidtha}{7cm}
\newcommand{\trimfiga}[2]{\trimw{#1}{#2}{.12}{.12}{.28}{.28}}
\newcommand{\figSize}{5.5cm}
\begin{figure}[htb]

\begin{center}
\begin{tikzpicture}
   \useasboundingbox (0,.75) rectangle (14,9.3);  

  \begin{scope}[yshift=0.0cm] 
     \drawContour{xshift=-.25cm,yshift=4.65cm}{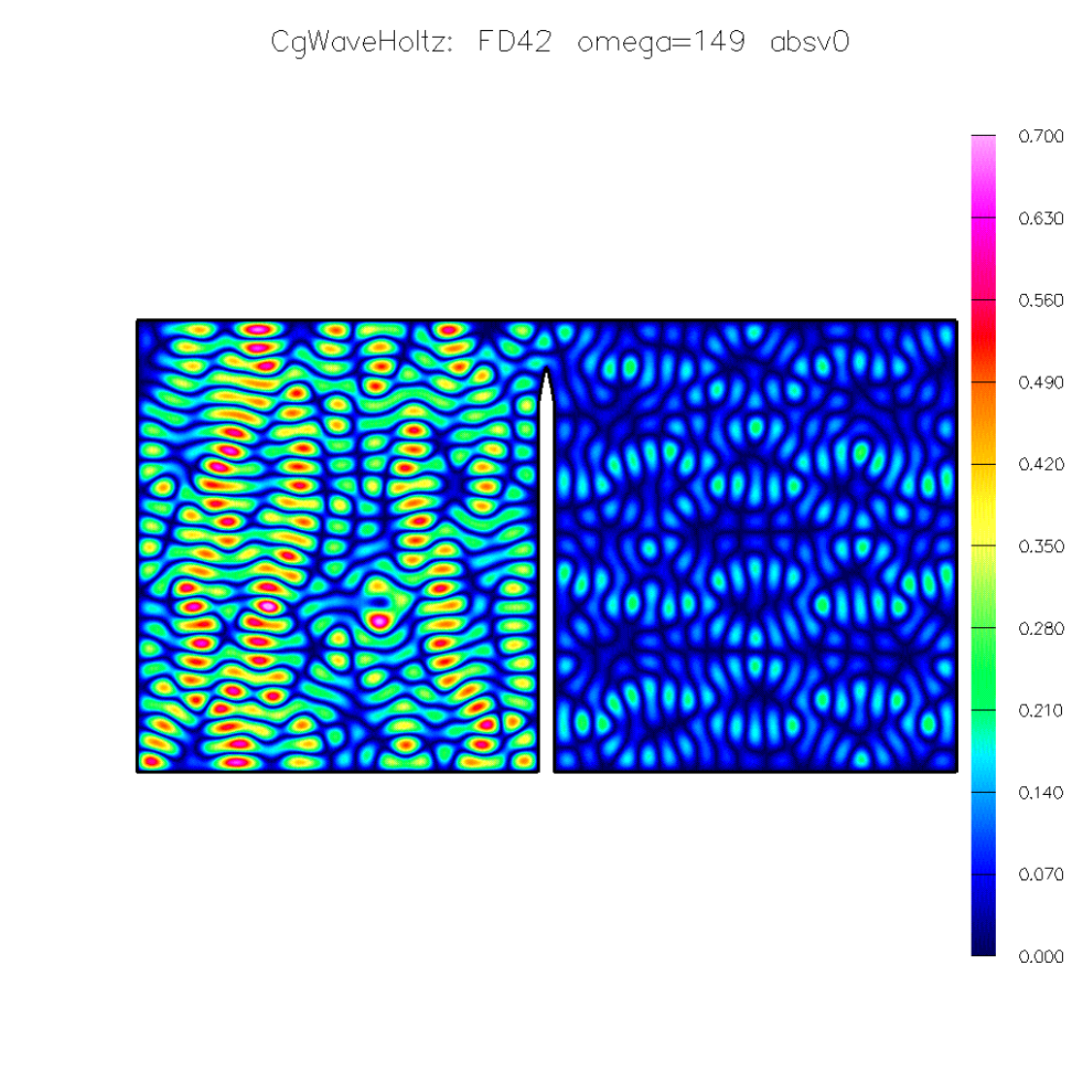}{$|v|$, $\Gcke^{(128)}$, $\PPW=54$}{$|v|$, $\Gcke^{(64)}$}{$t=0.3$}{$0$}{$1$}  
     \drawContour{xshift=-.25cm,yshift=0.50cm}{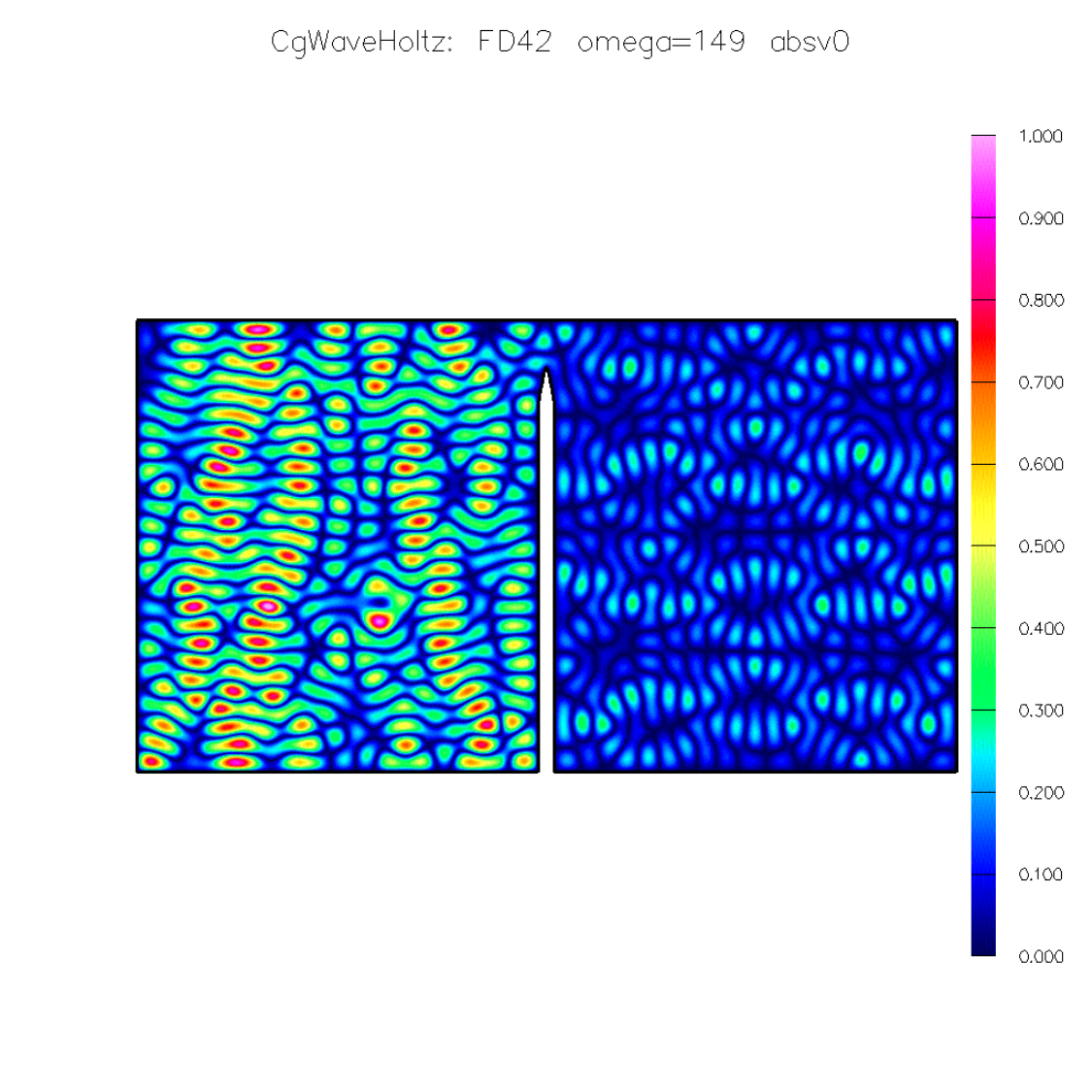}{$|v|$, $\Gcke^{(64)}$, $\PPW=27$}{$v$}{$t=0.3$}{$0$}{$1$}       
   \end{scope}  
   \begin{scope}[xshift=8cm,yshift=0cm]
    \figByWidth{0}{4.75}{tipGridSmallG128O4Freq149Deflate64}{\figSize}[0][0][0][0]
    \figByWidth{0}{0}{tipGridSmallG128O4Freq149Deflate64FixPointMuFunction}{\figSize}[0][0][0][0]
  \end{scope}      

\end{tikzpicture}
\end{center}
\caption{Knife edge Helmholtz solution for $\omega=149$,  Gaussian source located in the lower left at $(-.2,.2)$, 
  fourth-order accurate, deflate $65$ eigenmodes. Left: contours of $|v|$ for the fine grid $\Gcke^{(128)}$ (top) and coarse grid $\Gcke^{(64)}$ (bottom).
  Right: convergence history for the fine grid.
  }
\label{fig:doubleEllipseContoursFreq149}
\end{figure}
}

The overset grid for the knife edge geometry, denoted by $\Gcke^{(j)}$, is shown in Figure~\ref{fig:knifeEdgeGridAndContours},
and consists of four component grids. A background Cartesian grid covers the domain $[-0.5,0.5]\times[0,0.55]$.
Two other Cartesian grids lie adjacent to the lower straight sides of the knife, which has a total height of $0.5$ from its base to tip and a width of $0.02$. 
A curvilinear grid is used
to fit the boundary over the tip of the knife edge.
The nominal grid spacing is $\ds^{(j)}=1/(10 j)$, although the tip grid uses a finer mesh with stretching to resolve
the sharp tip of the knife edge.
Dirichlet boundary conditions are used on all boundaries.

\begin{table}[htb]\tableFont 
\begin{center}
\begin{tabular}{|c|c|c|c|c|c|c|c|} \hline 
  \multicolumn{8}{|c|}{Knife edge, Points-per-wavelength, $\eps=10^{-2}$, $L=1$}   \\  
  \multicolumn{5}{|c|}{} & Actual & \multicolumn{2}{c|}{Estimated} \\  \hline
   $\omega$  & $\Lambda$  & $N_\Lambda$  &  $j$     & $\ds$  &  $\PPW$  & $\PPW_2$ & $\PPW_4$   \\ \hline 
  $100$   &  $0.0628$  & $ 15.9$ &  $32$  &  3.13e-03 & $ 20.1$  &  $128.3$  & $ 17.2$  \\
  $100$   &  $0.0628$  & $ 15.9$ &  $64$  &  1.56e-03 & $ 40.2$  &  $128.3$  & $ 17.2$  \\
  $100$   &  $0.0628$  & $ 15.9$ &  $128$  &  7.81e-04 & $ 80.4$  &  $128.3$  & $ 17.2$  \\
\hline
  $149$   &  $0.0422$  & $ 23.7$ &  $32$  &  3.13e-03 & $ 13.5$  &  $156.6$  & $ 19.0$  \\
  $149$   &  $0.0422$  & $ 23.7$ &  $64$  &  1.56e-03 & $ 27.0$  &  $156.6$  & $ 19.0$  \\
  $149$   &  $0.0422$  & $ 23.7$ &  $128$  &  7.81e-04 & $ 54.0$  &  $156.6$  & $ 19.0$  \\
%
  \hline
\end{tabular}
\caption{Actual and estimated points-per-wavelength for the knife edge domain as a function of frequency $\omega$ and grid resolution $j$. 
 The column titled $\PPW$ holds the actual points-per-wavelength. The columns labeled $\PPW_2$ and $\PPW_4$ contains the rule-of-thumb estimated values
 for second- and fourth-order accurate schemes, respectively, from Recipe~\ref{recipe:PPW}.
    }
\label{tab:knifeEdgePPW}
\end{center}
\end{table}

The right plot of Figure~\ref{fig:knifeEdgeGridAndContours} shows contours of the magnitude of the numerical solution, denoted by $|v|$, from a fourth-order accurate computation
on grid $\Gcke^{(64)}$  using implicit time-stepping. 
The Gaussian source~\eqref{eq:gaussianSource}, with $\omega=100$, $\alphag=7000$ and $\betag=40$, is  located at $\xv_0=(-0.2,0.2)$ which is a point approximately centered in the portion of the domain to the left of the knife.
Using $\Nits=10$ time-steps per period, the implicit scheme has grid CFL numbers of about
$45$ on the Cartesian grids and as high as $866$ on the tip grid. In other words, 
the implicit scheme uses a time-step that is $866$ times larger than that required of an explicit scheme with a global time-step.
From Table~\ref{tab:knifeEdgePPW} this computation on grid $\Gcke^{(64)}$ has $\PPW=40$, while the rule of thumb estimate is $\PPW_4=17$
indicating that the computation should be well resolved (see comments below for the $\omega=149$ simulations).
The convergence history of the WaveHoltz iterations for this case are given in Figure~\ref{fig:knifeEdgeConvergenceFreq100}.
Results with $64$ eigenmodes deflated and without deflation are compared. 
With no deflation the FPI convergence is quite slow while the GMRES accelerated convergence starts out slow but then converges
rapidly. Note that the rapid convergence of GMRES begins at around iteration $60$; at this point GMRES has apparently identified many
of the slowly converging eigenmodes.

Results for a higher-frequency example are given in Figure~\ref{fig:doubleEllipseContoursFreq149}.
In this case the Gaussian source~\eqref{eq:gaussianSource} is defined using $\omega=149$, $\alphag=11000$, $\betag=40$, and $\xv_0=(-0.2,0.2)$.
The solution is computed with $65$ eigenmodes deflated.
The contours of $|v|$ indicate that the forcing has led to a solution that is primarily active in the left half of the domain where the source is approximately centered.
Figure~\ref{fig:doubleEllipseContoursFreq149} shows contours of $|v|$ for a coarse grid $\Gcke^{(64)}$ and fine grid $\Gcke^{(128)}$.
From Table~\ref{tab:knifeEdgePPW} the computation on grid $\Gcke^{(64)}$ had $\PPW=27$ while grid $\Gcke^{(128)}$ had $\PPW=54$.
The rule-of-thumb estimate is $\PPW_4=19$ which suggests that both simulations are well resolved; 
this is confirmed in the computations since the contour plots are nearly indistinguishable.
The right graphs of Figure~\ref{fig:doubleEllipseContoursFreq149} show the convergence history (on the fine grid $\Gcke^{(128)}$) when $65$ eigenmodes are deflated.
The FPI convergence is in good agreement with the theory. The GMRES accelerated convergence is very good.
We note that, as expected, the convergence history on coarse grid $\Gcke^{(64)}$ (not shown) is found to be nearly identical to that for the fine grid.

\section{Optimal $O(N)$ Helmholtz solver at fixed frequency} \label{sec:optimal}

We now describe how the components of the WaveHoltz algorithm can be assembled into an optimal $O(N)$ algorithm to solve Helmholtz problems at a fixed frequency, where $N$ denotes the total number of grid points.
The optimal algorithm has both an $O(N)$ CPU cost as well as an $O(N)$ memory cost.
The WH+ITS+MG+GMRES (WaveHoltz + Implicit-time-stepping + Multigrid + GMRES) algorithm is based on the following ingredients,
\begin{enumerate}
  \item WH: WaveHoltz fixed-point iteration,
  \item ITS: implicit time-stepping of the wave equation using a fixed number of time-steps per period (e.g. $\NITS=10$),
  \item MG: multigrid solution of the implicit time-stepping equations,
  \item GMRES : Krylov accelerated WaveHoltz solver.
\end{enumerate}
The $O(N)$ CPU-time and memory-usage scaling as a function of $N$ is based on the following observations.
\begin{enumerate}
  \item For a fixed frequency, theory and computations show that the convergence of the WaveHoltz FPI is essentially independent of the mesh size (see comments below).
  \item Computations suggest that the convergence of the GMRES accelerated WaveHoltz FPI is also essentially independent of the mesh size 
     (in any case this rate is at least as fast as the FPI convergence rate).
  \item The number of time-steps per WaveHoltz iteration can be fixed, independent of $N$ (e.g.~$10$ time-steps per period).
  \item The cost of the MG solution is $O(N)$ in CPU and $O(N)$ in memory.
  \item The memory use of the GMRES-accelerated WaveHoltz iteration is also $O(N)$ provided the number of GMRES iterations is fixed.
\end{enumerate}
For the results presented here we use the overset grid multigrid solver Ogmg~\cite{automg,multigridWithNonstandardCoarsening2023}.

{
\newcommand{\figSize}{7cm}

\begin{figure}[htb]
\begin{center}
\begin{tikzpicture}[scale=1]
  \useasboundingbox (0,.7) rectangle (7,6);  

  \begin{scope}[yshift=0cm]
    \figByWidth{0.0}{0}{betaWithCoarseFineLambda}{\figSize}[0][0][0][0];
  \end{scope}  

\end{tikzpicture}
\end{center}
\caption{The asymptotic convergence rate of the WaveHoltz filter is normally determined by the discrete eigenvalue closest to $\omega$. 
This figure shows values of $\beta$ evaluated at the coarse grid, fine grid, and continuous eigenvalues of the one-dimensional Laplacian.
As the mesh is refined the ACR approaches the value of $\beta$ at the eigenvalue $\lambda=1$.
    } 
\label{fig:betaWithCoarseAndFineLambda}
\end{figure}
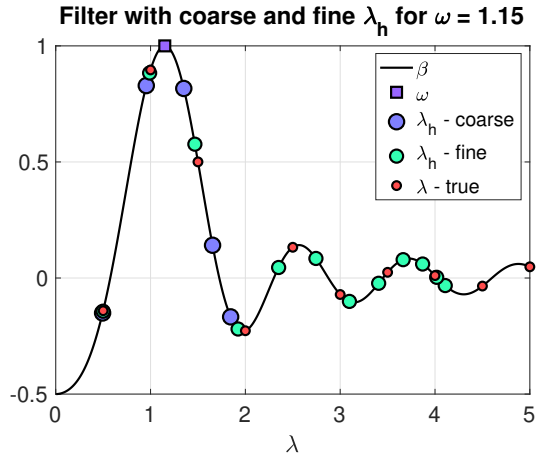
}

Figure~\ref{fig:betaWithCoarseAndFineLambda} illustrates why the convergence rate of the WaveHoltz FPI is essentially
independent of the mesh spacing. (Here we ignore the adjustments to $\lambda$ for finite $\dt$.) 
Consider the  one-dimensional Laplacian on the interval $[0,2\pi]$ with Dirichlet boundary conditions, and with eigenvalues $\lambda_m=m/2$, $m=1,2,3,\ldots$.
The figure shows the filter function $\beta$ evaluated at 
the true eigenvalues together with the eigenvalues of coarse and fine-grid discretizations of it. 
As the mesh is refined, the rate of convergence in WaveHoltz is set by the value of $\beta$ at the true eigenvalue closest to $\omega$. 
As the grid is refined the eigenvalues of the discretized problem converge to the eigenvalues of the continuous problem, and as a result the rate of convergence does not depend significantly on the grid spacing. Further, the poorly resolved eigenvalues are large and far away to the right along the $\lambda$-axis. These are damped rapidly during the WaveHoltz iteration; however when solving the Helmholtz BVP directly these cause ill-conditioning, the condition number then scales as $O (h^{-2})$, where $h$ is a measure of the mesh spacing.

\begin{table}[hbt]\tableFont 
\begin{center}
\begin{tabular}{|c|c|c|c|c|c|c|} \hline 
           &           & \multicolumn{3}{c|}{WaveHoltz} & \multicolumn{2}{c|}{GMRES+ILU(5)} \\ \hline
     grid  & N         &     its &  ECR   &  CPU (s) &       its  &  CPU (s)         \\ \hline 
   square  & $256^2$   &   $14$  &  0.58   &   $4.3$   &  5.3e2   &  1.7e0        \\ \hline
   square  & $512^2$   &   $14$  &  0.58   &   $16$    &  2.8e3   &  5.0e1        \\ \hline
   square  & $1024^2$  &   $14$  &  0.58   &   $63$    &  2.6e4   &  2.0e3        \\ \hline
   square  & $2048^2$  &   $14$  &  0.58   &  $258$    &  1.3e5   &  5.8e4        \\ \hline
\end{tabular}
\caption{Square: solving a Helmholtz problem with WaveHoltz (using multigrid to solve the implicit time-stepping equations) ($\Np=2$, $\NITS=10$) 
and directly with GMRES+ILU(5)+RESTART(50). The WaveHoltz CPU times are seen to approximately scale like a constant times $N$, the total number of grid points.
    }
\label{tab:squareHelmholtzComparison}
\end{center}
\end{table}

In the first example we use the WH+ITS+MG+GMRES scheme to solve the Helmholtz equation on the unit square in two dimensions with Dirichlet boundary conditions.
The frequency is taken as $\omega=11$. The number of periods is $\Np=2$ and there are $\NITS=10$ implicit time-steps per period.
The forcing is the Gaussian source in~\eqref{eq:gaussianSource} with center at $(x_0,y_0)=(0.4,0.4)$, amplitude $\alphag=-100$, and $\betag=20$.
A Cartesian grid is used with equal grid spacings in both directions.
Table~\ref{tab:squareHelmholtzComparison} shows results from a grid refinement study.
The WH+ITS+MG+GMRES algorithm is seen to converge with a fixed number of (GMRES accelerated) iterations (to a fixed tolerance) and the CPU time
is seen to scale linearly with $N$, the total number of grid points.
A graph of the normalized CPU time divided by $N$, versus $N$ is shown in the left plot of Figure~\ref{fig:CPUversusN}; 
the CPU time is normalized so the time for the coarsest grid is $1$.
For comparison, the discretized Helmholtz problem is solved with GMRES and an ILU preconditioner with $5$ levels of fill-in,
 and with a restart length of $50$ (denoted by GMRES+ILU(5)). In this case the number of iterations is seen to increase rapidly with $N$.
The results in Table~\ref{tab:squareHelmholtzComparison} show that the effective convergence rate (ECR), defined in~\eqref{eq:ECRdef}, is independent of $N$.

{
\newcommand{\drawContour}[7]{%
\begin{scope}[#1]
\draw(0.0,0) node[anchor=south west,xshift=-4pt,yshift=+0pt] {\trimfiga{#2}{\figWidtha}};
\begin{scope}[xshift=-.2cm,yshift=-5pt]
  \draw (\xcb,\ycb) node[anchor=south west,xshift=0.25cm,yshift=.5cm,rotate=-90] {\trimfigcb{colourBarLines}{\cbWidth}{\cbHeight}};
  \draw (.8,0) node[anchor=north,xshift=+3pt,yshift=+2pt] {\scriptsize $#6$};
  \draw (4.8,0) node[anchor=north,xshift=+0pt,yshift=+2pt] {\scriptsize $#7$};
\end{scope}
\end{scope}
}
\newcommand{\cbWidth}{.2cm}
\newcommand{\cbHeight}{4cm}
\newcommand{\xcb}{.5cm}
\newcommand{\ycb}{-.2cm}
\setlength{\ycbTop}{\ycb+\cbHeight}
\setlength{\ycbMid}{\ycb+\cbHeight*\real{.5}}
\newcommand{\trimfigcb}[3]{\includegraphics[width=#2, height=#3, clip, trim=17cm 2.35cm 1.65cm 2.35cm]{#1}}
\newcommand{\figWidtha}{4.5cm}
\newcommand{\trimfiga}[2]{\trimw{#1}{#2}{.11}{.115}{.11}{.11}}
\newcommand{\figSize}{5.5cm}
\begin{figure}[htb]

\begin{center}
\begin{tikzpicture}
   \useasboundingbox (0,.25) rectangle (16,4.8);  

  \begin{scope}[yshift=0.0cm] 
     \drawContour{xshift=-.3cm,yshift=0.00cm}{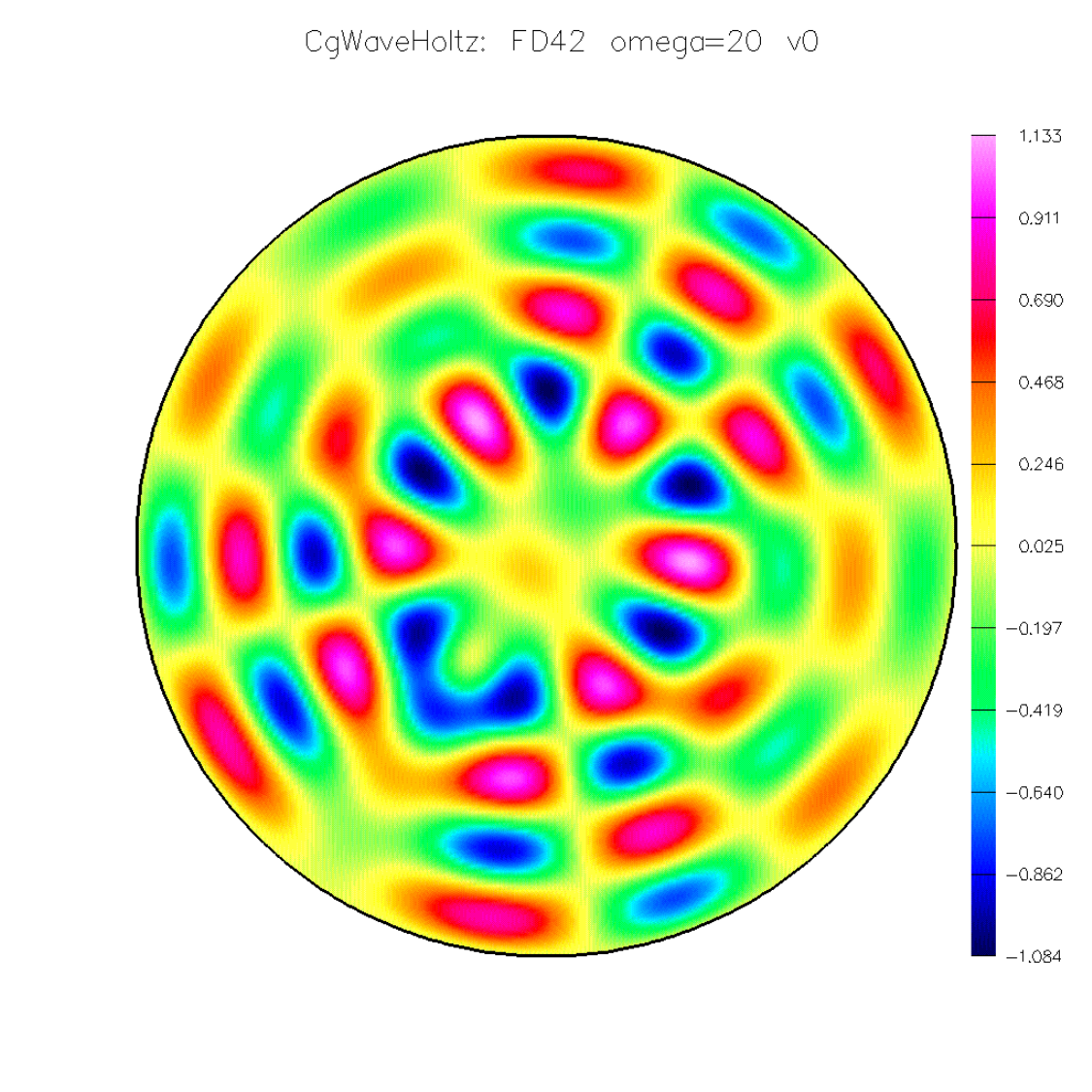}{G8 O4}{$v$}{$t=0.3$}{$-1.1$}{$1.1$}       
   \end{scope}  
   \begin{scope}[xshift=4.85cm,yshift=0cm]
    \figByWidth{0.00}{0}{sice8Freq20Deflate32}{\figSize}[0][0][0][0]
    \figByWidth{5.6}{0}{sice8Freq20Deflate32FixPointMuFunction}{\figSize}[0][0][0][0]
    \draw (.75,.9) node[draw,fill=white,anchor=west,xshift=0pt,yshift=0pt,inner sep=3pt] {\scriptsize Disk, G8, O4};
  \end{scope}      
\end{tikzpicture}
\end{center}
\caption{Disk: Helmholtz solution used for the CPU scaling study. Left: solution on grid $\Gcd^{(8)}$. Middle: WaveHoltz convergence. Right: WaveHoltz filter function. }
\label{fig:diskMultigridSolution}
\end{figure}
}

\begin{table}[hbt]\tableFont 
\begin{center}
\begin{tabular}{|c|c|c|c|c|} \hline 
           &          & \multicolumn{3}{c|}{WaveHoltz}  \\ \hline
     grid  & $\ds$    &     its &  ECR   &  CPU (s)     \\ \hline 
   disk    & $1/40$   &   $19$  &  0.57   &  $4.5$       \\ \hline
   disk    & $1/80$   &   $17$  &  0.54   &  $10.9$      \\ \hline
   disk    & $1/160$  &   $16$  &  0.52   &  $41.1$      \\ \hline
   disk    & $1/320$  &   $16$  &  0.52   &  $160$       \\ \hline
   disk    & $1/640$  &   $16$  &  0.52   &  $583$       \\ \hline
\end{tabular}
\caption{Disk: CPU times versus grid resolution. WaveHoltz plus implicit time-stepping and multigrid leads to an optimal $O(N)$ scheme. The CPU time
approximately doubles as the total number of grid points, $N$, doubles.
   Results for the disk domain, order of accuracy four, $\NITS=10$,  and $\Np=2$.}
\label{tab:diskHelmholtzMultigrid}
\end{center}
\end{table}

{
\newcommand{\figSize}{6cm}

\begin{figure}[htb]
\begin{center}
\begin{tikzpicture}[scale=1]
  \useasboundingbox (0,.7) rectangle (13,5);  

   \begin{scope}[yshift=0cm]
    \figByWidth{0.0}{0}{squareCPUversusN}{\figSize}[0][0][0][0];
    \figByWidth{7.0}{0}{diskCPUversusN}{\figSize}[0][0][0][0];
  \end{scope}  

\end{tikzpicture}
\end{center}
\caption{Normalized values of CPU$(N)$/$N$, versus number of grid points. Left: square. Right: disk. The CPU time is seen to scale linearly with $N$.
    } 
\label{fig:CPUversusN}
\end{figure}
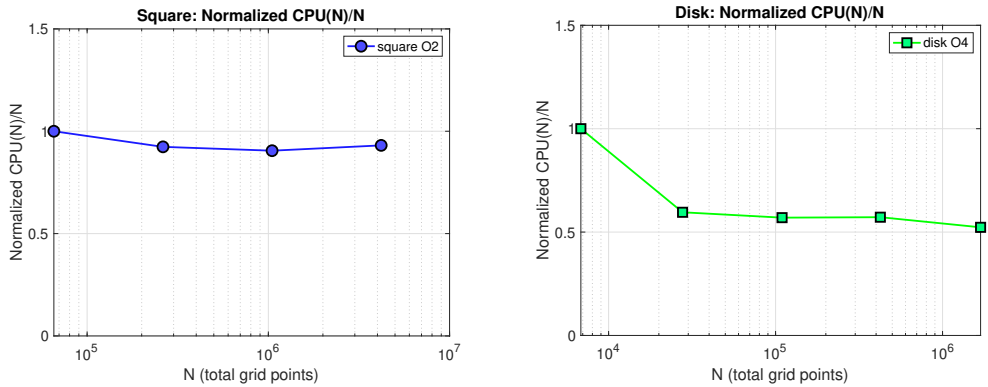
}

We next consider solving a Helmholtz problem on a two-dimensional disk.
The overset grid for this geometry denoted by $\Gcd^{(j)}$ was described in Section~\ref{sec:disk}.
The problem uses a Gaussian forcing with $\omega=20$, $\alphag=400$, $\betag=10$ and $\xv_0=(-0.2,-0.3)$.
The problem is solves using the fourth-order accurate WaveHoltz scheme with implicit time-stepping, $\NITS=10$ time-steps per period, $\Np=2$ periods, and with $58$ eigenmodes deflated.
Figure~\ref{fig:diskMultigridSolution} shows contours of the solution as well as the convergence of the FPI and GMRES accelerated
iterations for grid $\Gcd^{(8)}$.
Table~\ref{tab:diskHelmholtzMultigrid} gives the GMRES accelerated iteration counts and CPU times as the grid is refined.
The results show that the number of GMRES accelerated iterations stays nearly constant as the mesh is refined and that 
the CPU time approximately doubles as the total number of grid points $N$ doubles.
A graph of the CPU time (divided by $N$) versus $N$ is shown in the right plot of Figure~\ref{fig:CPUversusN}.
This gives further evidence of the near optimal CPU time complexity of the algorithm.
We do not show results using GMRES to solve the discrete Helmholtz equations directly since GMRES had great difficulty in solving this problem.
The fill-in level for the ILU preconditioner has to be taken so large as to make the approach almost a sparse direct solver.
For example, ILU$(300)$ is needed for $\Gcd^{(16)}$, and then GMRES converges in two iterations.

\section{Conclusions} \label{sec:conclusions}

We have described an efficient and high-order accurate solver for the Helmholtz equation in complex geometry.
The scheme is based on solving for time-periodic solutions of the related time-domain wave equation using the WaveHoltz algorithm.
WaveHoltz applies a time-filter to the time-dependent solution to remove unwanted frequencies in time.
Complex geometry is treated using overset grids and high-order finite difference schemes.
An optimal $O(N)$ algorithm is described that uses implicit time-stepping to advance the wave equation with a fixed number
of time-steps per period, independent of $N$. The implicit system of equations resulting from implicit time-stepping
can be efficiently solved with multigrid. It is shown how to correct for time discretization errors, even when taking very large time-steps.
GMRES is used to accelerate the basic WaveHoltz fixed-point iteration.
Deflation can be used to remove the slowest converging eigenmodes.
Numerical examples are given in two- and three-space dimensions to illustrate properties of the schemes including the benefits
of using high-order accurate schemes to over-come pollution effects.
A rule-of-thumb for determining the number of points-per-wavelength was derived from the analysis of a model problem and this estimate
 was shown to be useful in computations in complex geometry.
The optimal $O(N)$ behaviour of the algorithm is demonstrated on two examples.
An interesting finding was that no dissipation was needed in the wave equation solver when used with WaveHoltz.
Normally dissipation is needed for stability on overset grids but the WaveHoltz time-filter has apparently filtered out any unstable modes.
In future work we will consider problems with radiation boundary conditions where the Helmholtz solution is complex valued.
The WaveHoltz algorithm can be adjusted to solve for the complex valued solution, \cite{WaveHoltzSemi}.  
Helmholtz problems for systems of equations such as Maxwell's equations of electromagnetics have shown to be efficient in, \cite{peng2021emwaveholtz}, and a natural extension of this work is to the linear dispersive models for electromagnetics in~\cite{adegdm2019}.

\appendix

\section{Corrections for time-discretization errors} \label{sec:timeCorrections}

In this section we describe corrections to the time-stepping scheme
to adjust for time-discretization errors so that 
the solution obtained from the WaveHoltz algorithm matches the solution
to the discretized Helmholtz problem.
With the ability to remove the time-discretization errors, it is no longer necessary 
to solve the wave equation to high-order accuracy in time.
This means that second-order accurate schemes in time can be used, and these are generally more efficient.
Note that these corrections are especially important when using implicit time-stepping
with large time-steps, since the time-discretization errors can be large.

\newcommand{\Nstep}{N_t}
\subsection{Explicit time-stepping correction}  \label{sec:timeCorrectionsExplicit}

Consider the explicit time-stepping scheme (second-order accurate in time and  $p^{\rm th}$-order accurate in space) given by 
\ba
   \Dpt\Dmt W^n_\iv = L_{ph} W_\iv^n - f(\xv_\iv) \cos(\omegaExplicit t), \label{eq:waveExplicitOrder2}
\ea
where the time dependence of the forcing, $\cos(\omegaExplicit t)$, involves a modified frequency $\omegaExplicit$ whose form is yet to be determined.
The time periodic solution to~\eqref{eq:waveExplicitOrder2} is of the form
$W_\iv^n = U_\iv \cos(\omegaExplicit t^n)$, and substituting this into~\eqref{eq:waveExplicitOrder2} leads to 
a discrete Helmholtz equation for $U_\iv$,
\ba
  &  - \omegaHat^2 \, U_\iv =  L_{ph} U_\iv - f(\xv_\iv) ,   \label{eq:waveExplicitTimePeriodic}
\ea
where the frequency $\omegaHat$ is given by
\ba
  & \omegaHat \eqdef \f{\sin(\omegaExplicit\dt/2)}{\dt/2} . 
  \label{eq:omegaHat}
\ea
We actually wish to solve the following discrete Helmholtz problem,
\ba
  &  - \omega^2 U_\iv =  L_{ph} U_\iv - f(\xv_\iv) . 
  \label{eq:HelmholtzDiscrete}
\ea
Comparing~\eqref{eq:waveExplicitTimePeriodic} to~\eqref{eq:HelmholtzDiscrete} 
indicates that we want $\omegaHat=\omega$, and this implies choosing $\omegaExplicit$ to satisfy
\ba
   \f{\sin(\omegaExplicit\dt/2)}{\dt/2} = \omega.
\ea
Note that when using $\omegaExplicit$ in~\eqref{eq:waveExplicitOrder2}, we must solve the wave equation using the new period $\Ttilde=2\pi/\omegaExplicit$,
which in turn changes the time-step. Let $\dt$ denote this new time-step and $\Nstep$ denote the number of time steps. 
Then we require the following relations to hold
\ba
   & \f{\sin(\omegaExplicit\dt/2)}{\dt/2} = \omega,  \quad
   \Ttilde = \f{2\pi}{\omegaExplicit}, \quad
   \dt = \f{\Ttilde}{\Nstep} = \f{2\pi}{\omegaExplicit} \f{1}{\Nstep} .
\ea 

\mni
Here then is the form of time-correction.
\begin{recipe}[Explicit time-stepping correction]
\label{recipe:OmegaTildeExplicit}
First estimate the number of time-steps, $\Nstep$, from the original period $T$ based
on a time-step restriction for stability of the explicit scheme such as that in~\eqref{eq:dtForCartesian}. 
Given $\Nstep$, choose
\ba
  & \dt = \f{2}{\omega}\sin\left(\f{\pi}{\Nstep}\right)   , \qquad
    \omegaExplicit = \f{2\pi}{\Nstep\dt}.
\ea
Check that this new $\dt$ still satisfies the time-step restriction; if not then increase $\Nstep$ until it does.
Solve the wave equation~\eqref{eq:waveExplicitOrder2} using the modified frequency $\omegaExplicit$ and modified period $\Ttilde=2\pi/\omegaExplicit$.
\end{recipe}

\subsection{Implicit time-stepping correction}  \label{sec:timeCorrectionsImplicit}

Consider the implicit time-stepping scheme (second-order accurate in time and  $p^{\rm th}$-order accurate in space) given by 
\ba
   \Dpt\Dmt W^n_\iv = L_{ph}\Big[  \half W_\iv^{n+1} + \half W_\iv^{n-1}  \Big] - f(\xv_\iv)\, \cos(\omegaImplicit t) \cos(\omegaImplicit\dt), 
   \label{eq:waveImplicitOrder2}
\ea
where the time dependence of the forcing is chosen as $\cos(\omegaImplicit t) \cos(\omegaImplicit\dt)$ 
with a modified frequency~$\omegaImplicit$.  The factor $\cos(\omegaImplicit\dt)$ is included in the forcing to enable a convenient choice for $\omegaImplicit$ as shown below.
The time periodic solution to~\eqref{eq:waveImplicitOrder2}, $W_\iv^n = U_\iv \cos(\omegaImplicit t^n)$, satisfies,
\ba
  &  - \omegaHat^2 \, U_\iv =  \cos(\omegaImplicit\dt) L_{ph} U_\iv - f(\xv_\iv)  \cos(\omegaImplicit\dt),   \label{eq:waveImplicitTimePeriodic}
\ea
where $\omegaHat$ is given in~\eqref{eq:omegaHat}. Note that $\omegaHat^2$ can also be written as
\ba
  \omegaHat^2 = \f{2 - 2 \cos(\omegaImplicit\dt)}{\dt^2} .
\ea
Comparing the discrete Helmholtz equation~\eqref{eq:HelmholtzDiscrete}, that we want solve, 
to~\eqref{eq:waveImplicitTimePeriodic} indicates we should choose $\omegaImplicit$ so that
\ba
          &  \f{ 2 - 2\cos(\omegaImplicit\dt) }{\dt^2} \,\f{1}{\cos(\omegaImplicit\dt) } = \omega^2 , \label{eq:omegaTildeImpEqn}
\ea
and solving~\eqref{eq:omegaTildeImpEqn} for $\cos(\omegaImplicit\dt)$ gives a formula that can be used to find $\omegaImplicit$,
\ba          
   \cos(\omegaImplicit\dt) = \f{1}{1 + (\omega\dt)^2/2 }  . \label{eq:cosOmegaTilde}
\ea
We also need to change the period $T$ and $\dt$ so that
\ba
   &  \Nt \, \dt = \Ttilde = \f{2\pi}{\omegaImplicit},   \label{eq:impOmegaTildeEquation}
\ea
where $\Nt$ is the given number of time-steps.
Solving~\eqref{eq:impOmegaTildeEquation} for $\omegaImplicit$ and substituting into~\eqref{eq:cosOmegaTilde} gives
\ba
   &\cos(\omegaImplicit\dt) = \cos\left(\f{2\pi}{N_t}\right) = \f{1}{1 + (\omega\dt)^2/2 }, 
\ea
and thus $\dt$ must satisfy
\ba
   (\omega\dt)^2 = 2\left( \f{1}{\cos\bigl(2\pi/N_t\bigr)} -1 \right)  . \label{eq:implicitCorrectionDtII}
\ea
Note that equation~\eqref{eq:implicitCorrectionDtII} places a minor restriction on the allowable
number of time-steps since the right-hand side must be positive. We require
\ba
         &   \cos\left(\f{2\pi}{N_t}\right)  > 0   \quad
\implies~~  \f{2\pi}{N_t} < \f{\pi}{2}, 
\ea
giving $N_t>4$ and thus we must take at least 5 time-steps per period $T$ (not $\Tbar$),
\ba
  N_t \ge 5 .
\ea

\mni
Here then is the form of time-correction. 
\begin{recipe}[Implicit time-stepping correction]
\label{recipe:OmegaTildeImplicit}
Choose the desired number of time-steps $\Nstep \ge5$ per period $T$.
Compute $\dt$ and $\omegaImplicit$ using
\ba
   & \dt = \f{1}{\omega} \sqrt{ \f{2}{\cos\bigl(2\pi/N_t\bigr)} -2  },\qquad
   \omegaImplicit = \f{2\pi}{N_t\dt},
\label{eq:implicitCorrectionDt}
\ea
Solve the implicit time-stepping equations~\eqref{eq:waveImplicitOrder2} 
using the modified frequency $\omegaImplicit$ and modified period $\Ttilde=2\pi/\omegaImplicit$.
\end{recipe}


\subsection{Discrete filter function} \label{sec:discreteFilterFunction}

In this section we study the form of the discrete filter function $\beta$.
We give a corrected value for the value of $\alpha$ that appears in the $\beta$ function
to account for a potentially large $\dt$ to ensure the discrete $\beta$ function reaches a 
maximum at $\lambda=\omega$.

The discrete filter function is uses a trapezoidal rule quadrature,
\ba
   \beta_d(\lambda,\omega) \eqdef 
       \f{2}{T} \sum_{n=0}^{\Nt} \Big( \cos(\omega t^n) - \f{\alpha}{2}\Big) \, \cos(\lambda t^n) \, \sigma_n \, \dt,
     \label{eq:betad}
\ea
where $T=2\pi/\omega$, $\dt=T/\Nt$, and the quadrature weights are $\sigma_i=1$, $i=1,2,\ldots,\Nt-1$, and
$\sigma_0=\sigma_{\Nt}=\half$. 
It can be shown that $\beta_d$ takes a similar form to equation~\eqref{eq:filterThreeSincs}
 for the continuous~$\beta$, and in particular takes the form
\ba
   \beta_d(\lambda,\omega,T) = \sincd(\omega+\lambda,T) + \sincd(\omega-\lambda,T) -\alpha \, \sincd(\lambda,T),
     \label{eq:discreteBetaThreeSincs}
\ea
where $\sincd$ is an approximate $\sinc$ function defined by
\ba
  \sincd(z,T) = \f{\sin( z \, T)}{T \tan(z\dt/2)/(\dt/2)} .\label{eq:sincd}
\ea
The formula~\eqref{eq:discreteBetaThreeSincs} can be found by replacing 
terms such as $\cos(\lambda t^n)$ in \eqref{eq:betad} by complex exponentials, summing the approproiate geometric series, and then taking the real part.
Note that 
\ba 
  & \sincd(0,T) = 1, \\
  & \sincd( m \omega , T)  = 0, \qquad m=1,2,\ldots, 
\ea
and thus $\beta_d$ is one at $\lambda=\omega$,
\ba
   \beta_d(\omega,T) =1.
\ea
We also want $\beta_d$ to reach a maximum at $\lambda=\omega$.
Now,
\ba
  \f{d\beta_d}{d\lambda}(\omega,T) = 
                 \f{\omega\cos(2\omega T)}{T \tan(\omega\dt)/(\dt/2)}
      -\alpha \, \f{\omega\cos(\omega T)}{T \tan(\omega\dt)/(\dt/2)} , 
\ea
and setting this to zero implies
\ba
    \alpha = \alpha_d \eqdef \f{\tan(\omega\dt/2)}{\tan(\omega\dt)}.   \label{eq:alphadII}
\ea

\mni
\textbf{Summary.} The trapezoidal rule quadrature~\eqref{eq:betad} should use the corrected value of 
$\alpha=\alpha_d$ given in~\eqref{eq:alphadII}.

\section{Discrete dispersion relations for high-order accurate schemes.} \label{sec:discreteDispersion}

In this section we derive some results used in the derivation of the points-per-wavelength rule-of-thumb from Section~\ref{sec:pollutionErrors}.
We prove Theorem~\ref{th:uxxCoeff} which gives the coefficients $b_\mu$ in the discrete approximation to the 
second derivative and then derive the formula~\eqref{eq:kErrorOrderp} for the error in discrete dispersion relation for $\kTilde$.

High-order accurate approximations to the second derivative can be derived from the formal series expansion~\cite{GustafssonKreissOliger95}
\ba
   \p_x^2 u(x) = \Dpx\Dmx \sum_{\mu=0}^{\infty} b_{\mu} (-\dx^2 \Dpx\Dmx)^\mu \, u(x) .   \label{eq:SecondDerivExpansion}
\ea
We now prove Theorem~\ref{th:uxxCoeff} which provides an explicit formula~\eqref{eq:uxxCoeff} for the coefficients $b_\mu$.
\begin{proof}
Substituting $u=e^{i k x}$ into~\eqref{eq:SecondDerivExpansion} gives
\ba
   k^2  =  \f{4 \sin^2(k\dx/2)}{\dx^2} \sum_{\mu=0}^{\infty} b_{\mu} \big(4\sin^2(k\dx/2)\big)^\mu
        =  \f{1}{\dx^2} \sum_{\mu=0}^{\infty} b_{\mu} \big(4\sin^2(k\dx/2)\big)^{\mu+1}
   .  \label{eq:kTildeSymbolInfiniteOrder}
\ea
We proceed formally at this point, assuming the series~\eqref{eq:kTildeSymbolInfiniteOrder} converges. This assumption will be justified when the form for $b_\mu$ is found.
Introduce the normalized wave-number
\ba
   \xi \eqdef k \, \dx.
\ea
Then~\eqref{eq:kTildeSymbolInfiniteOrder} can be written as 
\ba
   \xi^2 = \sum_{\mu=0}^{\infty} b_{\mu} \big(4\sin^2(\xi/2)\big)^{\mu+1} . \label{eq:xiSeries}
\ea
One way to find the coefficients $b_{\mu}$ is to expand the right-hand side of~\eqref{eq:xiSeries} in a Taylor series about~$\xi=0$ and then equate coefficients of powers of~$\xi\sp2$.
This gives an expression for $b_\mu$ in terms of previous values $b_m$, $m=0,1,\ldots,\mu-1$.
To find a closed form expression for $b_\mu$, we instead proceed as follows. Rather than working with powers of $\sin(\xi/2)$ on the right-hand side of~\eqref{eq:xiSeries}, 
we follow~\cite{BanksBucknerHagstrom2022} and introduce 
\ba
   &  \eta \eqdef \sin(\xi/2), 
\ea
which gives $\xi = 2 \arcsin(\eta)$ and then~\eqref{eq:xiSeries} becomes
\ba
    4 \arcsin^2(\eta) = \sum_{\mu=0}^{\infty} b_{\mu} \big(4 \eta^2 \big)^{\mu+1} .
\ea
Substituting the Taylor series for $\arcsin^2(\eta)$
\ba
   \arcsin^2(\eta) = \half \sum_{n=1}^\infty \f{1}{n^2 \, {2n \choose n}} (2\eta)^{2n},  
\ea 
implies
\ba
   2 \sum_{n=1}^\infty \f{4^n}{n^2 \, {2n \choose n}} \eta^{2n} = \sum_{\mu=0}^{\infty} b_{\mu} \big(4 \eta^2 \big)^{\mu+1} 
     = \sum_{\mu=0}^{\infty} b_{\mu} 2^{2\mu+2} \eta^{2\mu+2} 
   = \sum_{n=1}^{\infty} b_{n-1} 4^{n} \eta^{2n} .
\ea 
Whence, equating powers of $\eta\sp2$ gives
\ba
   b_{n-1} = \f{2}{n^2 \, {2n \choose n}} , \qquad n=1,2,\ldots ~.
\ea
Condition~\eqref{eq:uxxCoeff} is obtained upon setting $n=\mu+1$.
\qed
\end{proof}

\bigskip
Now let us derive the formula~\eqref{eq:kErrorOrderp} for the error in $\kTilde$ at order $p$.
The $p\sp{{\rm th}}$-order accurate approximation to the second derivative uses a truncated version of~\eqref{eq:SecondDerivExpansion} given by
\ba
   \p_x^2 u(x) \approx \Dpx\Dmx \sum_{\mu=0}^{p/2-1} b_{\mu} (-\dx^2 \Dpx\Dmx)^\mu \, u(x) . 
\ea
The discrete form of the Helmholtz equation for the model problem~\eqref{eq:pollutionModelProblem} that uses this $p\sp{{\rm th}}$-order approximation is 
\ba
    \Dpx\Dmx \sum_{\mu=0}^{p/2-1} b_{\mu} (-\dx^2 \Dpx\Dmx)^\mu \, U_j  + k^2 \, U_j = f(x_j).   \label{eq:HelmholtzOrderp}
\ea
Substituting $U_j = e^{i \kTilde x_j}$ into~\eqref{eq:HelmholtzOrderp} with $f(x_j)=0$ gives the discrete dispersion 
relation at order $p$, 
\ba
    k^2 = \f{1}{\dx^2} \sum_{\mu=0}^{p/2-1} b_{\mu} \big(4\sin^2(\kTilde\dx/2)\big)^{\mu+1} . \label{eq:disperionRelationOrderp}
\ea 
Equation~\eqref{eq:disperionRelationOrderp} gives an implicit relation between the discrete wave-number $\kTilde$ and $k$.
We wish to find an expression for the error between $k$ and $\kTilde$. 
To this end note that the infinite series~\eqref{eq:kTildeSymbolInfiniteOrder} holds for any $k$ provided $k\dx \ll 1$, 
and in particular the series holds for $\kTilde$, provided $\kTilde\dx \ll 1$,
\ba
   \kTilde^2  = \f{1}{\dx^2} \sum_{\mu=0}^{\infty} b_{\mu} \big(4\sin^2(\kTilde\dx/2)\big)^{\mu+1} .  \label{eq:kTildeSymbolInfiniteOrderKtilde}
\ea
Taking~\eqref{eq:disperionRelationOrderp} minus~\eqref{eq:kTildeSymbolInfiniteOrderKtilde} gives
\bse
\ba
   k^2  & = \kTilde^2  - \f{1}{\dx^2}  \sum_{\mu=p/2}^{\infty} b_{\mu} \big(4\sin^2(\kTilde\dx/2)\big)^{\mu+1} , \\
        & = \kTilde^2 -  \f{1}{\dx^2} \Big[  b_{p/2} (\kTilde\dx)^{p+2} + O\big( (\kTilde\dx)^{p+4} \big) \Big],   \\
        & = \kTilde^2 \Big[ 1 - b_{p/2} (\kTilde\dx)^{p}  +  O \big(\kTilde\dx)^{p+2}\big) \Big] .
\ea
\ese
Taking the square-root of both sides implies   
\ba     
   k &= \kTilde \Big[ 1 - b_{p/2} (\kTilde\dx)^{p} +  O \big( (\kTilde\dx)^{p+2} \big)  \Big ]^{1/2} .  \label{eq:dispersionSquareRoot}
\ea
But to the same order of approximation we can replace $(\kTilde\dx)^p$ with $(k\dx)^p$ and $(\kTilde\dx)^{p+2}$ with $(k\dx)^{p+2}$ 
in~\eqref{eq:dispersionSquareRoot} and  then solving for $\kTilde$ 
gives
\ba
  \kTilde & = k \, \Big[ 1 - b_{p/2} (k\dx)^{p} +  O \big( (k\dx)^{p+2} \big)  \Big]^{-1/2} . \label{eq:dispersionWithSquareRoot}
\ea
Using the first term in the binomial expansion in~\eqref{eq:dispersionWithSquareRoot} leads to the desired relation for the error in $\kTilde$ 
\shadedBoxWithShadow{align}{orange}{
     \kTilde & = k \, \Big[ 1 + \half b_{p/2} (k\dx)^{p} +  O \big( (k\dx)^{p+2}\big)\, \Big] .
}

%
\bibliographystyle{elsart-num}
\bibliography{journal-ISI,henshaw,henshawPapers,WaveHoltz,appelo,helmholtz}
 
\end{document}